\documentclass[12pt,a4paper,twoside]{article}
\usepackage[T1]{fontenc}
\usepackage{lmodern}
\usepackage[a4paper,margin=22mm]{geometry}
\usepackage{authblk,amsmath,amssymb,amsthm,amscd,mathtools}
\usepackage{graphicx,xcolor,colortbl,fancybox,wrapfig,multirow,bigdelim}
\usepackage{mathrsfs,arydshln,enumitem,tikz}
\usepackage[all]{xy}
\usepackage{microtype}
\let\mathrsfs\mathscr
\usepackage[colorlinks=true,allcolors=black]{hyperref}
\usetikzlibrary{automata}
\newcommand{\HH}{\mathcal H}
\newcommand{\EE}{\mathbb E}
\newcommand{\dd}{\,\mathrm d}
\newcommand{\gp}[2]{[#1]_{#2}^{(\gamma)}}
\theoremstyle{definition}
\newtheorem{theorem}{Theorem}[section]
\newtheorem{prop}[theorem]{Proposition}
\newtheorem{lemma}[theorem]{Lemma}
\newtheorem{lem}[theorem]{Lemma}
\newtheorem{corollary}[theorem]{Corollary}

\newtheorem{remark}[theorem]{Remark}

\newtheorem{conjecture}[theorem]{Conjecture}
\newenvironment{demo}[1]{%
  \trivlist
  \item[\hskip\labelsep
        {\bf #1.}]
}{%
\hfill\qedsymbol
  \endtrivlist
}
 \makeatletter
    
    \@addtoreset{equation}{section}
  \makeatother

\def\Z{\mathbb{Z}}

\def\C{\mathbb{C}}

\def\Real{\mathbb{R}}

\def\Norm{{\mathcal{N}}}
\def\sig{{\boldsymbol{\sigma}}}

\def\bla{\boldsymbol{\lambda}}

\def\i{\boldsymbol{i}}

\def\x{\boldsymbol{x}}
\def\y{\boldsymbol{y}}

\title{
Minor summation formulas for hyperpfaffians\\
and Selberg integrals
}
\author[1]{Masao ISHIKAWA\footnote{Partially supported by JSPS KAKENHI Grant Numbers JP16K05068, JP20K03558.}}
\affil[1]{\small Department of Mathematics,
Okayama University, Tsushima, Okayama 700-8530, Japan,
{\tt mi@math.okayama-u.ac.jp}
}
\author[2,3]{Jiang ZENG}
\affil[2]{\small College of Mathematics and Physics,
Wenzhou University, Wenzhou 325035, China;}
\affil[3]{
 ​Université Claude Bernard Lyon 1, CNRS UMR 5208, ICJ, 69100 Villeurbanne, France,
{\tt zeng@math.univ-lyon1.fr}
}
\date{
\small Revised for SIGMA: 30 October 2026.}
\def\defterm#1{{\sl #1}\/}

\renewcommand\tilde{\widetilde}
\newcommand\Pf{\operatorname{Pf}}

\newcommand\Hf{\operatorname{Hf}}

\newcommand\Sym{{\mathfrak{S}}}

\newcommand\sgn{\operatorname{sgn}}

\newcommand\rdots{\mathinner{\mkern1mu\raise0pt\vbox{\kern7pt\hbox{.}}
     \mkern2mu\raise4pt\hbox{.}\mkern2mu\raise8pt\hbox{.}\mkern1mu}}
\def\covered{\mathinner{\mkern1mu\raise0pt\vbox{\kern7pt\hbox{$<$}}
     \mkern-4mu\raise2pt\hbox{.}\mkern2mu}}
\def\covers{\mathinner{\mkern1mu\raise0pt\vbox{\kern7pt\hbox{$>$}}
     \mkern-12mu\raise2pt\hbox{.}\mkern8mu}}

\def\Cat{{\operatorname{Cat}}}
\def\Mot{{\operatorname{Mot}}}
\def\Del{{\operatorname{Del}}}
\def\Sch{{\operatorname{Sch}}}
\def\CBC{{\operatorname{CBC}}}
\def\CTC{{\operatorname{CTC}}}
\def\bB{B}
\def\cB{{\mathrsfs{B}}}

\def\cA{{\mathrsfs{A}}}
\begin{document}
\maketitle
\begin{abstract}
We establish a hyperpfaffian--hyperdeterminant minor summation formula and
derive hyperpfaffian analogues of de Bruijn's integral identities. A
confluent-Vandermonde specialization expresses Hankel hyperpfaffians of
moment sequences as Selberg-type integrals. We apply this method to
Narayana polynomials of types $A$, $B$ and $D$, obtaining evaluations that
include conjectures from our earlier work with Tagawa. The moment
representations distinguish the atom occurring in type $A$ from the
shifted signed measures needed in type $D$. We also evaluate $q$-Hankel
Pfaffians associated with little $q$-Jacobi and Al-Salam--Carlitz
polynomials. Further applications concern Jacobi ensembles with
$\beta=l^2$: polynomial insertions give modified moment hyperpfaffians,
and Kadell's integral yields finite evaluations for integer powers of
the characteristic polynomial, including a quadratic double sum.
We conclude with product conjectures for the Gessel--Xin sequences.
\end{abstract}
\tableofcontents
\section{Introduction}
Moment determinants and Pfaffians connect orthogonal polynomials with
multiple integrals and combinatorial identities. For determinants, the
Vandermonde factor leads naturally to Selberg-type integrals; for Pfaffians,
de Bruijn's integration formulas provide the corresponding link. Higher
Vandermonde powers call for higher-order alternating analogues of
Pfaffians. The purpose of this paper is to develop minor summation and
integration formulas for these hyperpfaffians and to use them to evaluate
expressions built from classical and $q$-orthogonal polynomial moments.
The applications include Narayana polynomials of types $A$, $B$ and $D$,
little $q$-Jacobi and Al-Salam--Carlitz polynomials, and polynomial
observables in Jacobi ensembles.

Our algebraic starting point is a hyperpfaffian--hyperdeterminant minor
summation formula (Theorem~\ref{th:msf1}). It extends the Pfaffian minor
summation formulas of Ishikawa and Wakayama \cite{IW1,IW2} and includes
formulas of Luque and Thibon \cite[Eq.~(87)]{LT1} and Matsumoto
\cite[Theorem~5.1]{Ma1} as special cases. We use it to obtain hyperpfaffian
analogues of de Bruijn's two integration formulas \cite{Bru}; see
Theorems~\ref{th:msf-integral} and~\ref{th:msf2-integral}. The latter
extends the hyperpfaffian formula of Luque and Thibon \cite{LT1}.
This approach complements their work on Hankel hyperdeterminants
\cite{LT2} and on hyperdeterminantal evaluations of the Selberg and Aomoto
integrals \cite{LT3}, as well as Matsumoto's study of Toeplitz
hyperdeterminants and Jack functions \cite{Ma2}.

The main link with moment sequences can be stated explicitly. Following
the Hankel terminology of \cite{LT1,LT2,LT3}, we call a skew-symmetric
hypermatrix of the form
\[
 A(i_1,\ldots,i_l)=
 \prod_{1\le p<q\le l}\bigl(g(i_q)-g(i_p)\bigr)
 f(i_1+\cdots+i_l)
\]
a hypermatrix of Hankel type, and its hyperpfaffian a Hankel hyperpfaffian.
The cases used here have $g(i)=i$ or $g(i)=q^i$.
Let $l$ be a positive even integer, let $n\ge1$, and let $d\psi$ be a finite
positive or signed Borel measure on $[0,R]$, with moments
$\mu_k=\int_0^R x^k\,d\psi(x)$. Put $c_l=\prod_{j=0}^{l-1}j!$.
For an integer $t$ such that $t+\binom l2\ge0$, we have
\begin{align}
&\Pf^{[l,1]}\left(\prod_{p<q}(i_q-i_p)
 \mu_{\sum_p i_p-l+t}\right)_{1\le i_1<\cdots<i_l\le ln}\notag\\
&\qquad=\frac{c_l^n}{n!}\int_{[0,R]^n}
 \prod_{j=1}^n x_j^{t+\binom l2}
 \prod_{i<j}(x_j-x_i)^{l^2}\prod_{j=1}^n d\psi(x_j).
\label{eq:intro-moment}
\end{align}
Here $\Pf^{[l,1]}$ denotes the order-$l$ hyperpfaffian defined in
Section~\ref{sec:Pfaffian}; for $l=2$ it is the ordinary Pfaffian.
Corollary~\ref{cor:Hankel-hyperpfaff} gives a direct proof by a confluent
Vandermonde determinant, together with extensions to other shifts under
absolute integrability assumptions. Formula~\eqref{eq:intro-moment}
explains the exponent $l^2$ and reduces a Hankel hyperpfaffian evaluation
to a multiple integral. Its measure formulation also accommodates atoms
and signed measures, both of which occur in the Narayana applications.
For positive measures, it further shows that varying the nonnegative
shift produces a Stieltjes moment sequence
(Corollary~\ref{cor:stieltjes-preservation}).

The principal evaluations use Selberg's integral \cite{Se},
\begin{align}
S_{n}(\alpha,\beta,\gamma)
&=\int_{[0,1]^n}
\prod_{i=1}^{n}t_{i}^{\alpha-1}\left(1-t_{i}\right)^{\beta-1}
\prod_{1\leq i<j\leq n}\left|t_{i}-t_{j}\right|^{2\gamma}\,d\boldsymbol{t}
\nonumber\\&
=\prod_{j=1}^{n}\frac
{\Gamma\left(\alpha+(j-1)\gamma\right)\Gamma\left(\beta+(j-1)\gamma\right)\Gamma\left(j\gamma+1\right)}
{\Gamma\left(\alpha+\beta+(n+j-2)\gamma\right)\Gamma\left(\gamma+1\right)},
\label{eq:Selberg}
\end{align}
where $d\boldsymbol{t}=dt_1\cdots dt_n$.
For $n>1$, convergence holds when $\Re\alpha,\Re\beta>0$ and
$\Re\gamma>-\min\{1/n,\Re\alpha/(n-1),\Re\beta/(n-1)\}$;
for $n=1$, only the conditions on $\alpha$ and $\beta$ are needed.
Our applications have $\gamma=l^2/2>0$. For background and further
connections, see the survey of Forrester and Warnaar \cite{FW}; discrete
analogues with combinatorial applications appear in \cite{BKW}.
Aomoto's extension \cite{Ao} evaluates, for $0\le k\le n$,
\begin{align}
&\int_{[0,1]^n} \left(\prod_{i=1}^kt_i\right)
\prod_{i=1}^{n}t_{i}^{\alpha-1}\left(1-t_{i}\right)^{\beta-1}
\prod_{1\leq i<j\leq n}\left|t_{i}-t_{j}\right|^{2\gamma}\,d\boldsymbol{t}\nonumber\\&\qquad
=S_n(\alpha, \beta, \gamma) \prod_{j=1}^{k}\frac{\alpha+(n-j)\gamma}
{\alpha+\beta+(2n-j-1)\gamma}.
\label{eq:Aomoto}
\end{align}
By symmetry, this is equivalent to
\begin{align}
&\int_{[0,1]^n} e_{k}(\boldsymbol{t})
\prod_{i=1}^{n}t_{i}^{\alpha-1}\left(1-t_{i}\right)^{\beta-1}
\prod_{1\leq i<j\leq n}\left|t_{i}-t_{j}\right|^{2\gamma}\,d\boldsymbol{t}\nonumber\\&\qquad
=\binom{n}{k}S_n(\alpha, \beta, \gamma) \prod_{j=1}^{k}\frac{\alpha+(n-j)\gamma}
{\alpha+\beta+(2n-j-1)\gamma},
\label{eq:Aomoto2}
\end{align}
where $e_k(\boldsymbol{t})$ is the elementary symmetric polynomial,
defined by
$\prod_{i=1}^n(1+t_i y)=\sum_{k=0}^n e_k(\boldsymbol{t})y^k$.
These two integral evaluations give product and finite-sum formulas for
the moment hyperpfaffians considered below.

One motivation comes from our earlier work with Tagawa \cite{ITZ2} on
Pfaffians associated with orthogonal polynomial moments. That paper
formulated conjectures for Hankel Pfaffians involving Al-Salam--Carlitz
moments, Narayana and related combinatorial sequences, and a sequence
arising from Tamm's Hankel determinant \cite{Ta}.
In the present paper, Theorem~\ref{th:tildeN} gives hyperpfaffian
evaluations for the Narayana polynomials of types $A$, $B$ and $D$.
Their specializations recover and extend the identities in
\cite[Conjecture~6.2]{ITZ2}, previously proved by Ishikawa and Koutschan
\cite{IK} using a computer-assisted application of Zeilberger's holonomic
ansatz. The integral proofs identify the moment measures underlying
these evaluations. In type $A$, an atom at the origin must be taken into
account; in type $D$, the relevant shifted moments require signed
measures. Keeping these distinctions explicit is essential when
applying~\eqref{eq:intro-moment}.

The $q$-applications concern ordinary Pfaffians. A representative
example, proved in \cite[Corollary~3.2]{ITZ2} by a Pfaffian decomposition,
is
\begin{align}
&\Pf\left((q^{i-1}-q^{j-1})\frac{(aq;q)_{i+j+r-2}}{(abq^2;q)_{i+j+r-2}}\right)_{1\leq i<j\leq 2n}
\nonumber\\&=
a^{n(n-1)}q^{n(n-1)(4n+1)/3+n(n-1)r}
\prod_{k=1}^n\frac{(aq;q)_{2k+r-1}(bq;q)_{2(k-1)}(q;q)_{2k-1}}
{(abq^2;q)_{2(k+n)+r-3}}.
\label{eq:ITZ}
\end{align}
Here we use the standard $q$-shifted factorial notation \cite{AAR,GR}:
for $|q|<1$,
\[
 (a;q)_\infty=\prod_{k=0}^\infty(1-aq^k),\qquad
 (a;q)_n=\frac{(a;q)_\infty}{(aq^n;q)_\infty},
\]
with integer $n$ whenever the expression is defined, and
$(a_1,\ldots,a_r;q)_n=\prod_{j=1}^r(a_j;q)_n$ for integer $n$ or
$n=\infty$.
A second proof of~\eqref{eq:ITZ}, based on a quadratic formula for basic
hypergeometric series related to Askey--Wilson polynomials, was given in
\cite[Corollary~3.4]{GITZ2}. Our proof identifies its integral origin:
de Bruijn's formula reduces the Pfaffian to the $k=2$ case of the
Askey--Habsieger--Kadell $q$-Selberg integral
\eqref{eq:ask-conj1}. We also use multivariate Al-Salam--Carlitz
polynomials \cite{AC,BF,Ch} to prove the two identities in
\cite[Conjecture~6.1]{ITZ2} and obtain further shifted evaluations
(Theorem~\ref{conj:03}). Although
Corollary~\ref{cor:q-Hankel-hyperpfaff} provides an integral
representation for $q$-Hankel hyperpfaffians of arbitrary even order,
the explicit $q$-evaluations established here use order two.

Formula~\eqref{eq:intro-moment} also leads to applications in
$\beta$-ensembles with $\beta=l^2$. Hyperpfaffian representations of
partition functions were studied by Sinclair \cite{Sin}, and their
correlation functions are developed by Sinclair and Wells \cite{SW}.
In Section~\ref{sec:applications}, we express polynomial insertions in
Jacobi ensembles as ratios of modified moment hyperpfaffians.
Combining this representation with Kadell's Selberg--Jack integral
\cite{Kad2} gives finite evaluations for every nonnegative integer
power of the characteristic polynomial, with an explicit double sum
for the quadratic case. The underlying Selberg averages are classical;
the contribution here is their explicit realization in terms of Hankel
hyperpfaffians with modified moments. We also record moment formulas
for correlation functions and gap probabilities, and discuss connections
with time-deformed hyperpfaffians \cite{Sinclair2026} and constellation
ensembles \cite{Wolff}.

The paper is organized as follows. Section~\ref{sec:Pfaffian} introduces
the hyperdeterminants and hyperpfaffians used here and proves the minor
summation formula. Section~\ref{sec:Bruijn} develops the integration
formulas and their moment specializations.
Section~\ref{sec:Motzkin} evaluates Hankel hyperpfaffians of Narayana
polynomials by the Selberg and Aomoto integrals.
Section~\ref{sec:Askey} treats little $q$-Jacobi moments and gives the
integral proof of~\eqref{eq:ITZ};
Section~\ref{sec:Al-Salam-Carlitz} treats Rogers--Szeg\"o polynomials and
the related Al-Salam--Carlitz moments.
Section~\ref{sec:applications} develops the ensemble applications.
The final section presents product conjectures for the Gessel--Xin
sequences, including Conjecture~6.3 of~\cite{ITZ2}, which remains open.

\section{Minor summation for hyperpfaffians}\label{sec:Pfaffian}
The aim of this section is to give a hyperpfaffian-hyperdeterminant version of the minor summation formula
\cite[Theorem~1]{IW1}, \cite[Theorem~3.2]{IW2}.
First we give a general form of hyperpfaffian-hyperdeterminant version
in Theorem~\ref{th:msf1},
and then derive corollaries which will be applied 
to obtain formulas linking the so-called Hankel hyperpfaffians to
the Selberg type integrals in the next section.
Some special cases of Theorem~\ref{th:msf1} are obtained in \cite[(87)]{LT1} and \cite[Theorem~5.1]{Ma1}.
The relation between the so-called Hankel Pfaffians and the Selberg integrals is fully understood
if we generalize Pfaffian to hyperpfaffian.
In fact the Selberg integral seems a very general formula 
and the application to evaluations of Pfaffians needs only a special case.
If we pay attention to the hyperpfaffian, it will reveal a full view of the relation.
However, more general arguments of hyperpfaffian may apparently look cumbersome.
Hence this section will provide the complete proof of our theorems in a general form.
\par\smallskip
First of all,
we fix some notation.
Let $\{a_{1},\dots,a_{k}\}_{<}$ denote the set 
$\{a_{1},\dots,a_{k}\}\subseteq \Real$ where $a_{1}<\cdots<a_{k}$.
We use $[n]$ to denote the set $\{1,2,\dots,n\}$ for any positive integer $n$.
Let $\binom{S}{r}$ denote the set of all $r$-element subsets of $S\subseteq \Real$.
If $S=\{s_{1},\dots,s_{n}\}_{<}\subseteq\Real$ is an $n$-element set, 
then we use the symbol $S_{I}$ to denote 
the $r$-element subset $\{s_{i_{1}},\dots,s_{i_{r}}\}_{<}$
for $I=\{i_{1},\dots,i_{r}\}_{<}\in\binom{[n]}{r}$.
Hence we have $\binom{S}{r}=\left\{ S_{I}\mid I \in \binom{[n]}{r}\right\}$.
More generally,
let $l$ and $n$ be positive integers and $S$ any $ln$-element set.
Let $\binom{S}{l,\dots,l}$ ($l$ is repeated $n$ times in the lower line) denote the set of
$n$-tuples $(S_{1},\dots,S_{n})$ such that $S_{j}\in\binom{S}{l}$ for $1\leq j\leq n$
and $S_{1}\uplus\cdots\uplus S_{n}=S$ (disjoint union).
In other wordss we have
$\binom{S}{l,\dots,l}=\left\{(S_{I_{1}},\dots,S_{I_{n}})\mid (I_{1},\dots,I_{n})\in\binom{[ln]}{l,\dots,l}\right\}$.
\par\smallskip
Let $\Sym_{n}$ be the set of permutations of $[n]$,
which form the symmetric group of degree $n$.
We use the one-line notation 
$\sigma=(\sigma(1),\dots,\sigma(n))$
to denote an element $\sigma$ of $\Sym_{n}$.
If $\sigma\in\Sym_{ln}$,
we split $\sigma$ into $n$ blocks of length $l$,
i.e.,
$\displaystyle
\sigma=\left(
\sig_{1},\dots,\sig_{n}
\right)
$
with $\sig_{k}=(\sigma((k-1)l+1),\dots,\sigma(kl))$ for $1\leq k\leq n$.
Let $\Sym_{ln}(l,\dots,l)$ ($n$ times $l$) denote the subset of $\Sym_{ln}$
which is composed of the permutations
$
\sigma=\left(
\sig_{1},\dots,\sig_{n}
\right)
$,
in which the entries of each block are strictly increasing,
i.e.,
\[
\Sym_{ln}(l,\dots,l)=\left\{\,
\sigma=(\sigma(1),\dots,\sigma(ln))\in\Sym_{ln}
\,\big|\,
\sigma((j-1)l+1)<\cdots<\sigma(jl)
\text{ for $1\leq j\leq n$}
\,\right\}.
\]
For example, %$\Sym_{4}(2,2)$ is as follows:
we have
\begin{equation*}
\Sym_{4}(2,2)=
\{(1,2,3,4),
(1,3,2,4),
(1,4,2,3),
(2,3,1,4),
(2,4,1,3),
(3,4,1,2)\}.
\end{equation*}
Let $S \subseteq \Real$ be an $ln$-element set.
There is a natural bijection
from $\Sym_{ln}(l,\dots,l)$ to $\binom{S}{l,\dots,l}$ defined by
\[
\sigma=\left(
\sig_{1},\dots,\sig_{n}
\right)\in \Sym_{ln}(l,\dots,l)
\mapsto
(S_{\sig_{1}},\dots,S_{\sig_{n}})\in \binom{S}{l,\dots,l}.
\]
Hence we can identify an element $I\in\binom{S}{l,\dots,l}$ 
with the corresponding permutation $\sigma\in\Sym_{ln}(l,\dots,l)$,
and use the symbol $\sgn I$ to denote $\sgn\sigma$.
For example,
if $l=n=2$ and $S=\{1,3,5,7\}$,
then we have
\begin{align*}
\binom{S}{2,2}
&=\bigl\{
(\{1,3\},\{5,7\}), (\{1,5\},\{3,7\}), (\{1,7\},\{3,5\}), 
\\&\qquad\qquad
(\{3,5\},\{1,7\}), (\{3,7\},\{1,5\}), (\{5,7\},\{1,3\})
\bigr\}.
\end{align*}
The signs of those elements are $+$, $-$, $+$, $+$, $-$ and $+$, respectively.
\par\smallskip
We divide this section into three subsections.
In Subsection~\ref{sub:hyperdet} we recall the definition of hyperdeterminant, 
and prove some basic properties, 
Proposition~\ref{prop:exterior}, Proposition~\ref{prop:Laplace}.
Especially
Lemma~\ref{lem:hyper-minor} is crucial to prove Theorem~\ref{th:msf1}.
Next,
in Subsection~\ref{sub:hyperpfaff}, 
 we introduce a block-array version of the hyperpfaffian,
which includes Barvinok's original definition \cite{Bar}.
Lemma~\ref{lem:subhyperpfaffian} is a hyperpfaffian generalization of a useful Pfaffian identity,
which will be used to prove the main theorem of the next section,
i.e.,
Theorem~\ref{th:msf2-integral}.
Finally, 
in Subsection~\ref{sub:hyperMSF}, 
we state a hyperpfaffian-hyperdeterminant version of the minor summation formula,
i.e.,
Theorem~\ref{th:msf1}, and two consequences (Corollary~\ref{cor:msf-det} and Theorem~\ref{cor:Pf-Hf}).
\subsection{Hyperdeterminant}\label{sub:hyperdet}
Let $m$ be an even positive integer and $n$ be a positive integer.
An \defterm{$m$-dimensional tensor of order $n$}
is a map $A:[n]^m\to F$, $(i_{1},\dots,i_{m})\mapsto A(i_{1},\dots,i_{m})$,
where $F$ is a field of characteristic $0$.
We express this tensor by
$A=(A(i_{1},\dots,i_{m}))_{1\leq i_{1},\dots,i_{m}\leq n}$,
which we call an \defterm{$m$-dimensional hypermatrix of size $n$}.
We also say the $(i_{1},\dots,i_{m-1})$th row of $A$ is 
$(A(i_{1},\dots,i_{m-1},1),\dots,A(i_{1},\dots,i_{m-1},n))$.
The \defterm{hyperdeterminant} $\det^{[m]}A$ of $A$
is defined to be
\begin{align}
{\det}^{[m]}A
&=\frac1{n!}\sum_{\sigma_{1},\dots,\sigma_{m}\in\Sym_{n}}
\sgn(\sigma_{1}\cdots\sigma_{m})
\prod_{i=1}^{n}A(\sigma_{1}(i),\sigma_{2}(i),\dots,\sigma_{m}(i))
\nonumber\\&
=\sum_{\sigma_{1},\dots,\sigma_{m-1}\in\Sym_{n}}
\sgn(\sigma_{1}\cdots\sigma_{m-1})
\prod_{i=1}^{n}A(\sigma_{1}(i),\dots,\sigma_{m-1}(i),i).
\label{eq:hyperdeterminant}
\end{align}
Especially the $m=2$ case corresponds to the ordinary determinant $\det A$.
For example, if $m=4$ and $n=2$, then
\begin{align*}
{\det}^{[4]}A
&=A(1,1,1,1)A(2,2,2,2)-A(1,1,1,2)A(2,2,2,1)-A(1,1,2,1)A(2,2,1,2)
\\
&+A(1,1,2,2)A(2,2,1,1)
-A(1,2,1,1)A(2,1,2,2)+A(1,2,1,2)A(2,1,2,1)\\
&+A(1,2,2,1)A(2,1,1,2)-A(1,2,2,2)A(2,1,1,1).
\end{align*}
\begin{remark}
\rm
Note that, if $m$ is even,
the second equality in the above definition always holds.
Meanwhile, 
if $m$ is odd and $n>1$, 
then the first sum is equal to zero,
but the second one is not.
In \eqref{eq:hyperdeterminant}, 
the last index is fixed and the first $(m-1)$ indices are permuted
by $\Sym_n$.
This kind of sum varies depending on which index one fixes.
The arguments in this section are consistent if we use the sum in \eqref{eq:hyperdeterminant}
even if $m$ is odd.
\end{remark}
As we define submatrix,
we can define a subhypermatrix of a hypermatrix $A$ 
by restricting the range of each index to a subset of $[n]$.
Fix an $m$-tuple $(I^{(1)},\dots,I^{(m)})$ such that $I^{(1)},\dots, I^{(m)}\in \binom{[n]}r$.
Let $A_{I^{(1)},\dots,I^{(m)}}$ denote the \defterm{subhypermatrix} of $A$ whose $k$th index $i^{(k)}$ moves in $I^{(k)}$
for $1\leq k\leq m$, i.e.,
\[
A_{I^{(1)},\dots,I^{(m)}}=A(i^{(1)}_{\nu_{1}},\dots,i^{(m)}_{\nu_{m}})_{\nu_{1},\dots,\nu_{m}\in [r]},
\]
where $I^{(k)}=\{i^{(k)}_{1},\dots,i^{(k)}_{r}\}_{<}$ for $1\leq k\leq m$.
However, 
since it is cumbersome to use too many sub-indices, 
we simply write %this subhypermatrix as
\[
A_{I^{(1)},\dots,I^{(m)}}=(A(i^{(1)},\dots,i^{(m)}))_{i^{(1)}\in I^{(1)},\dots,i^{(m)}\in I^{(m)}}.
\]
The hyperdeterminant ${\det}^{[m]}A_{I^{(1)},\dots,I^{(m)}}$
is called
the $(I^{(1)},\dots,I^{(m)})$-\defterm{minor} of $A$. 
For instance,
if $m=2$, then $A_{I^{(1)},I^{(2)}}$ denotes the submatrix of $A$
obtained by choosing the rows in $I^{(1)}$ and the columns in $I^{(2)}$.
Let $\{\xi_{i}\}_{i\geq1}$ be anti-commutative symbols,
i.e. $\xi_{j}\xi_{i}=-\xi_{i}\xi_{j}$ for $i,j\geq1$.
The following proposition is elementary.
\begin{prop}
\label{prop:exterior}
Let $m$ and $n$ be positive integers.
If we put
\begin{equation}
\eta_{i}=\sum_{i^{(1)},\dots,i^{(m-1)}\in[n]} A(i^{(1)},\dots,i^{(m-1)},i)
\,\xi_{i^{(1)}}\oplus\cdots\oplus\xi_{i^{(m-1)}}
\label{eq:eta_i}
\end{equation}
for $i\in[n]$,
then
we have
\begin{equation}
\eta_{I}=\eta_{i_{1}}\cdots\eta_{i_{r}}
=\sum_{I^{(1)},\dots,I^{(m-1)}\in\binom{[n]}{r}}
{\det}^{[m]} A_{I^{(1)},\dots,I^{(m-1)},I}\,\,\xi_{I^{(1)}}\oplus\cdots\oplus\xi_{I^{(m-1)}}
\label{eq:eta}
\end{equation}
for any $I=\{i_{1},\dots,i_{r}\}_{<}\in\binom{[n]}{r}$,
where $\xi_{I^{(k)}}=\xi_{i^{(k)}_{1}}\cdots\xi_{i^{(k)}_{r}}$
for $I^{(k)}=\{i^{(k)}_{1},\dots,i^{(k)}_{r}\}_{<}$.
Especially,
if $r=n$ and $I=[n]$,
we obtain
\begin{equation}
\eta_{1}\cdots\eta_{n}={\det}^{[m]}(A)\cdot\underbrace{\xi\oplus\cdots\oplus\xi}_{(m-1) \text{ times}},
\label{eq:eta-full}
\end{equation}
where $\xi=\xi_{1}\cdots\xi_{n}$.
\end{prop}
\begin{demo}{Proof}
By \eqref{eq:eta_i} we have 
\[
\eta_{i_{k}}=\sum_{i^{(1)}_{k},\dots,i^{(m-1)}_{k}\in[n]}
A(i^{(1)}_{k},\dots,i^{(m-1)}_{k},i_{k})
\,\xi_{i^{(1)}_{k}}\oplus\cdots\oplus\xi_{i^{(m-1)}_{k}}.
\]
Multiplying this and using anti-commutativity,
we obtain
\begin{align*}
\eta_{I}&=\eta_{i_{1}}\cdots\eta_{i_{r}}
\\&
=\sum_{i^{(1)}_{1},\dots,i^{(1)}_{r}}\cdots\sum_{i^{(m-1)}_{1},\dots,i^{(m-1)}_{r}}
A(i^{(1)}_{1},\dots,i^{(m-1)}_{1},i_{1})\cdots A(i^{(1)}_{r},\dots,i^{(m-1)}_{r},i_{r})
\\&\qquad\qquad\qquad
\times\xi_{i^{(1)}_{1}}\cdots\xi_{i^{(1)}_{r}}\oplus\cdots\oplus\xi_{i^{(m-1)}_{1}}\cdots\xi_{i^{(m-1)}_{r}}
\\&
=\sum_{I^{(1)}=\{i^{(1)}_{1},\dots,i^{(1)}_{r}\}_{<}\in\binom{[n]}{r}}
\cdots
\sum_{I^{(m-1)}=\{i^{(m-1)}_{1},\dots,i^{(m-1)}_{r}\}_{<}\in\binom{[n]}{r}}\quad
\\&\qquad\qquad
\sum_{\sigma^{(1)},\dots,\sigma^{(m-1)}\in\Sym_{r}}
\sgn\sigma^{(1)}\cdots\sgn\sigma^{(m-1)}
A(i^{(1)}_{\sigma^{(1)}(1)},\dots,i^{(m-1)}_{\sigma^{(m-1)}(1)},i_{1})\cdots
\\&\qquad\qquad
\cdots A(i^{(1)}_{\sigma^{(1)}(r)},\dots,i^{(m-1)}_{\sigma^{(m-1)}(r)},i_{r})
\,\xi_{I^{(1)}}\oplus\cdots\oplus\xi_{I^{(m-1)}}
\\&=
\sum_{I^{(1)},\dots,I^{(m-1)}\in\binom{[n]}{r}}
{\det}^{[m]} A_{I^{(1)},\dots,I^{(m-1)},I}\,\,
\xi_{I^{(1)}}\oplus\cdots\oplus\xi_{I^{(m-1)}},
\end{align*}
which is the desired identity.
\end{demo}
Next, let $I^{(1)},\dots,I^{(m)}\in \binom{[n]}{r}$.
\defterm{The $(I^{(1)},\dots,I^{(m)})$-cofactor of $A$} is defined to be
\begin{equation}
\widetilde a_{I^{(1)},\dots,I^{(m)}}
=(-1)^{|I^{(1)}|+\cdots+|I^{(m)}|}\,{\det}^{[m]}A_{\overline{I^{(1)}},\dots,\overline{I^{(m)}}},
\end{equation}
where $|I^{(k)}|=\sum_{i\in I^{(k)}}i$ and 
$\overline{I^{(k)}}=[n]\setminus I^{(k)}$ stands for the complement of $I^{(k)}$ in $[n]$.
In the following proposition we show that 
we can generalize Proposition~2.1 of
\cite{Ma1} (early version of \cite{Ma2})
whereas we will not use this proposition in this paper.
\begin{prop}
\label{prop:Laplace}
Let $m$ and $n$ be positive integers.
For any $m$-dimensional hypermatrix 
$A=(A(i_{1},\dots,i_{m}))_{1\leq i_{1},\dots,i_{m}\leq n}$
 of size $n$,
we have
\begin{equation}
{\det}^{[m]}A
=(-1)^{m\binom{r+1}2}\sum_{I^{(1)},\dots,I^{(m-1)}\in\binom{[n]}{r}}
a_{I^{(1)},\dots,I^{(m-1)},I}
\,\widetilde a_{I^{(1)},\dots,I^{(m-1)},I}
\end{equation}
for any $I\in\binom{[n]}{r}$.
Here 
we use the symbol $a_{I^{(1)},\dots,I^{(m-1)},I}$
to denote the $(I^{(1)},\dots,I^{(m-1)},I)$-minor of $A$
for simplicity.
\end{prop}
\begin{demo}{Proof}
Let $\eta_i$ be as in \eqref{eq:eta_i}.
Then,
for any $I\in\binom{[n]}{r}$,
we have, from \eqref{eq:eta},
\begin{align*}
\eta&=\eta_{1}\cdots\eta_{n}=(-1)^{|I|-\binom{r+1}{2}}\eta_{I}\eta_{\overline I}
\\&=(-1)^{|I|-\binom{r+1}{2}}
\sum_{I^{(1)},\dots,I^{(m-1)}\in\binom{[n]}{r}}
\sum_{J^{(1)},\dots,J^{(m-1)}\in\binom{[n]}{n-r}}
a_{I^{(1)},\dots,I^{(m-1)},I}\,
a_{J^{(1)},\dots,J^{(m-1)},\overline I}
\\&\qquad\qquad\qquad\qquad\qquad\qquad\qquad\qquad\times
\xi_{I^{(1)}}\xi_{J^{(1)}}\oplus\cdots\oplus\xi_{I^{(m-1)}}\xi_{J^{(m-1)}}.
\end{align*}
Note that 
\[
\xi_{I^{(k)}}\xi_{J^{(k)}}=\begin{cases}
(-1)^{|I^{(k)}|-\binom{r+1}2}\xi &\text{ if $I^{(k)}\cap J^{(k)}=\emptyset$,}\\
0&\text{ otherwise.}
\end{cases}
\]
Hence we obtain
\[
\eta=(-1)^{m\binom{r+1}2}\sum_{I_{1},\dots,I_{m-1}\in\binom{[n]}{r}}
a_{I^{(1)},\dots,I^{(m-1)},I}\,\widetilde a_{I^{(1)},\dots,I^{(m-1)},I}
\,\xi\oplus\cdots\oplus\xi.
\]
Meanwhile, we have already shown that $\eta={\det}^{[m]} A\cdot\xi\oplus\cdots\oplus\xi$.
Hence we obtain the desired identity.
Note that $(-1)^{m\binom{r+1}2}=1$ if $m$ is even.
\end{demo}
The following lemma gives the Laplace expansion 
of a hyperdeterminant along the last index 
and one can state it in more general form.
But this form 
is enough to use
in the proof of Theorem~\ref{th:msf1}.
\begin{lemma}
\label{lem:hyper-minor}
Let $l$, $m$, $n$ and $N$ be positive integers such that $ln\leq N$.
Let  
$A=(A(i_{1},\dots,i_{m}))_{1\leq i_{1},\dots,i_{m}\leq N}$
be an $m$-dimensional hypermatrix of size $N$.
Let $I=\{i_{1},\dots,i_{ln}\}_{<}\in\binom{N}{ln}$
and $I=(I_{1},\dots,I_{n})$ be the $l$-block notation of $I$,
i.e.,
$I_{j}=\{i_{(j-1)l+1},\dots,i_{jl}\}_{<}$ for $1\leq j\leq n$.
Then we have
\begin{align}
{\det}^{[m]}A_{[ln],\dots,[ln],I}
=
\sum_{I^{(1)},\dots,I^{(m-1)}\in\binom{[ln]}{l,\dots,l}}
\sgn I^{(1)}\cdots\sgn I^{(m-1)}
\prod_{j=1}^{n}
{\det}^{[m]}A_{I^{(1)}_{j},\dots,I^{(m-1)}_{j},I_{j}},
\label{eq:hyper-minor}
\end{align}
where $I^{(k)}=(I^{(k)}_{1},\dots,I^{(k)}_{n})$ with $I_{j}^{(k)}\in \binom{[ln]}{l}$
for $1\leq j\leq n$ and $1\leq k\leq m-1$.
\end{lemma}
\begin{demo}{Proof}
Let
\[
\eta_{i}=\sum_{i^{(1)},\dots,i^{(m-1)}\in[ln]}
A(i^{(1)},\dots,i^{(m-1)},i)
\,
\xi_{i^{(1)}}\oplus\cdots\oplus\xi_{i^{(m-1)}}
\]
as in \eqref{eq:eta_i}.
Then, since
$\eta_{I}=\eta_{i_{1}}\cdots\eta_{i_{ln}}=\eta_{I_{1}}\cdots\eta_{I_{n}}$,
we have, by \eqref{eq:eta},
\begin{align*}
\eta_{I}&=\sum_{I^{(1)}_{1},\dots,I^{(m-1)}_{1}\in\binom{[ln]}{l}}\cdots
\sum_{I^{(1)}_{n},\dots,I^{(m-1)}_{n}\in\binom{[ln]}{l}}
{\det}^{[m]}A_{I^{(1)}_{1},\dots,I^{(m-1)}_{1},I_{1}}\cdots
 {\det}^{[m]}A_{I^{(1)}_{n},\dots,I^{(m-1)}_{n},I_{n}}
\\&\qquad\qquad\qquad\qquad\times
\xi_{I^{(1)}_{1}}\cdots\xi_{I^{(1)}_{n}}\oplus\cdots\oplus
\xi_{I^{(m-1)}_{1}}\cdots\xi_{I^{(m-1)}_{n}}.
\intertext{Note that
$\xi_{I^{(k)}_{1}}\cdots\xi_{I^{(k)}_{n}}$ vanishes
unless $[ln]$ is disjoint union of $(I^{(k)}_{1},\dots,I^{(k)}_{n})$
for $1\leq k\leq m-1$.
If we put $I^{(k)}=(I^{(k)}_{1},\dots,I^{(k)}_{n})$ for $1\leq k\leq m-1$,
then 
$\xi_{I^{(k)}_{1}}\cdots\xi_{I^{(k)}_{n}}=\sgn I^{(k)}\,\xi_{[ln]}$ implies}
\eta_{I}&=
\sum_{I^{(1)},\dots,I^{(m-1)}\in\binom{[ln]}{l,\dots,l}}
\sgn I^{(1)}\cdots\sgn I^{(m-1)}\,
{\det}^{[m]}A_{I^{(1)}_{1},\dots,I^{(m-1)}_{1},I_{1}}\cdots 
\\&\qquad\qquad\qquad\qquad\qquad
\cdots {\det}^{[m]}A_{I^{(1)}_{n},\dots,I^{(m-1)}_{n},I_{n}}
\,\xi_{[ln]}\oplus\cdots\oplus\xi_{[ln]}.
\end{align*}
Meanwhile,
by \eqref{eq:eta-full}, we have 
$\eta_{I}={\det}^{[m]}A_{[ln],\dots,[ln],I}\,\xi_{[ln]}\oplus\cdots\oplus\xi_{[ln]}$,
which proves the desired identity.
\end{demo}
It is clear that we can obtain a similar expansion formula along the $k$th index.
\subsection{Hyperpfaffian}\label{sub:hyperpfaff}
Barvinok \cite{Bar} gave the first definition of hyperpfaffian,
and Matsumoto used a variant of hyperpfaffian in \cite{Ma2}.
Here we give a slightly generalized definition,
which unifies the original definitions in \cite{Bar,LT1,Ma2,Sin}.
Ordinary Pfaffian is defined for a skew-symmetric matrix,
that is a matrix $\bB=(\bB(i,j))_{i,j\in[2n]}$ satisfying $\bB(j,i)=-\bB(i,j)$.
Since this condition implies $\bB(i,i)=0$, 
a skew-symmetric matrix is regarded as a map $\cB:\binom{[2n]}{2}\to F$, $\{i,j\}_{<}\mapsto \bB(i,j)$.
For a positive even integer $l$,
Barvinok considered an $l$-alternating tensor $\bB:[ln]^{l}\to F$, 
$\i=(i_{1},\dots,i_{l})\mapsto \bB(\i)=\bB(i_{1},\dots,i_{l})$ of order $ln$
which satisfies $\bB(\tau(\i))=\sgn\tau\cdot \bB(\i)$ for $\tau \in \Sym_{l}$,
where $\tau(\i)$ stands for $\left(i_{\tau(1)},\dots,i_{\tau(l)}\right)$.
The reader can easily see that 
 this $l$-alternating tensor $\bB$ of order $ln$ is 
equivalent to giving a map
$\cB:\binom{[ln]}{l}\to F$
defined by $\cB(\{i_{1},\dots,i_{l}\}_{<})=\bB(i_{1},\dots,i_{l})$. %for $i_{1}<\cdots<i_{l}$.
In this paper we consider
 more general situation.
\par\smallskip
Let $l$, $m$ and $n$ be positive integers,
and let 
$
\bB:\underbrace{[ln]^{l}\times \cdots \times [ln]^{l}}_\text{$m$ times} \to F$, $(\i^{(1)},\dots,\i^{(m)})\mapsto \bB(\i^{(1)},\dots,\i^{(m)})
$
be an $lm$-dimensional tensor of order $ln$ with $\i^{(1)},\dots,\i^{(m)}\in[ln]^{l}$.
We say that this tensor 
$\bB=\left(\bB(\i_1,\dots,\i_{m})\right)_{\i_1,\dots,\i_{m}\in[ln]^{l}}$ 
is \defterm{$l$-alternating} if it satisfies
\[
B(\tau_{1}(\i^{(1)}),\dots,\tau_{m}(\i^{(m)}))
=\sgn(\tau_1)\cdots\sgn(\tau_{m})\,
B(\i^{(1)},\dots,\i^{(m)})
\]
for any $(\tau_1,\dots,\tau_{m})\in(\Sym_l)^{m}$.
This means that, if we permute each index,
then the sign of the permutation is multiplied.
We can associate a map
$\cB:\underbrace{\binom{[ln]}{l}\times \cdots \times \binom{[ln]}{l}}_\text{$m$ times} \rightarrow F$ 
with $\bB$
by setting its value
$\cB(I^{(1)},\dots,I^{(m)})=\bB(\i^{(1)},\dots,\i^{(m)})$, 
where each $l$-tuple
$\i^{(k)}=(i^{(k)}_{1},\dots,i^{(k)}_{l})$ 
has the same entries as
$I^{(k)}=\{i^{(k)}_{1},\dots,i^{(k)}_{l}\}_{<}$ 
 ($1\leq k\leq m$).
This latter map 
$\cB:\underbrace{\binom{[ln]}{l}\times \cdots \times \binom{[ln]}{l}}_\text{$m$ times} \rightarrow F$ is denoted by
\[
\cB=\left(\cB(I^{(1)},\dots,I^{(m)})\right)_{I^{(1)},\dots,I^{(m)}\in\binom{[ln]}{l}},
\]
and called a \defterm{$m$-dimensional $l$-block array of size $ln$}.
An $lm$-dimensional $l$-alternating tensor $\bB$ of order $ln$ and 
an $m$-dimensional $l$-block array $\cB$ of size $ln$ are equivalent 
since we can recover $\bB$ from $\cB$ by the $l$-alternating property.
\par\smallskip
The \defterm{hyperpfaffian} $\Pf^{[l,m]}(\cB)$ of an $m$-dimensional $l$-block array $\cB$ of size $ln$ 
is defined to be
\begin{align}
&\Pf^{[l,m]}(\cB) =
\frac1{n!}\sum_{I^{(1)},\dots,I^{(m)}\in\binom{[ln]}{l,\dots,l}}
\sgn(I^{(1)})\cdots\sgn(I^{(m)})
\prod_{j=1}^{n}
\cB\left(I^{(1)}_{j},\dots, I^{(m)}_{j}\right),
\label{eq:hyperpfaffian}
\end{align}
where $I^{(k)}=(I_{1}^{(k)},\dots,I_{n}^{(k)})$ with $I_{j}^{(k)}\in\binom{[ln]}{l}$ 
for $1\leq j\leq n$ and $1\leq k\leq m$.
We shall write 
$
\Pf^{[l,m]}\left(\cB(I_1,\dots,I_{m})\right)_{I_1,\dots,I_{m}\in\binom{[ln]}{l}}
$
for $\Pf^{[l,m]}(\cB)$.
By removing the sign, 
we define \defterm{hyperhafnian $\Hf^{[l,m]}(\cB)$ of $\cB$} as
\begin{align}
&\Hf^{[l,m]}(\cB) =
\frac1{n!}\sum_{I^{(1)},\dots,I^{(m)}\in\binom{[ln]}{l,\dots,l}}
\prod_{j=1}^{n}
\cB(I^{(1)}_{j},\dots, I^{(m)}_{j}),
\label{eq:hyperhafnian}
\end{align}
which reduces to the ordinary Hafnian when $l=2$ and $m=1$.
We use it only in Theorem~\ref{cor:Pf-Hf}.
Note that, if $l$ and $m$ are both odd and $n>1$
then $\Pf^{[l,m]}(\cB)$ is always zero.
Hence we assume $l$ is always even hereafter.
The $m=1$ case reduces to \cite[(3.3.1)]{Bar}, \cite[(76)]{LT1}, \cite[Sec.1.1]{Sin}
 and the $l=2$ case reduces to \cite[(2.2)]{Ma2}.
Especially the $(l,m)=(2,1)$ case reduces to the ordinary Pfaffian
\[
\Pf^{[2,1]}(\bB) =\frac1{n!}\sum_{{\sigma=(\sigma_{1}(1),\sigma_{1}(2),\dots,\sigma_{n}(1),\sigma_{n}(2))\in\Sym_{2n}}
\atop
{\sigma_{1}(1)<\sigma_{1}(2),\dots,\sigma_{n}(1)<\sigma_{n}(2)}}
\sgn\sigma \cdot \bB(\sigma_{1}(1),\sigma_{1}(2)) \cdots \bB(\sigma_{n}(1),\sigma_{n}(2)),
\]
which we simply write $\Pf(\bB)$.
But, in fact, our new definition of hyperpfaffian $\Pf^{[l,m]}(\cB)$ 
does not essentially generalize  Barvinok's original hyperpfaffian $\Pf^{[l,1]}(\cB)$
since one can write $\Pf^{[l,m]}(\cB)$ in the form of $\Pf^{[lm,1]}(\widetilde \cB)$ by 
taking a %proper 
suitable $\widetilde \cB$
in a similar way as Matsumoto stated in \cite[Proposition~A1]{Ma2}.
We will prove it at the end of this section,
i.e.,
Proposition~\ref{rem:equi-Pf}.
The reader may skip it since it is not used in this paper.
\begin{prop}
Given an $m$-dimensional $l$-block array
\[
\cB=\left(\cB(I^{(1)},\dots,I^{(m)})\right)_{I^{(1)},\dots,I^{(m)}\in \binom{[ln]}{l}}
\]
 of size $ln$,
we put
\begin{equation}
\zeta=\sum_{I^{(1)},\dots,I^{(m)}\in\binom{[ln]}{l}}
\cB(I^{(1)},\dots,I^{(m)})\,\xi_{I^{(1)}}\oplus\cdots\oplus \xi_{I^{(m)}}.
\label{eq:zeta-b}
\end{equation}
Then we have
\begin{equation}
\zeta^n=n!\Pf^{[l,m]}(\cB)\,\underbrace{\xi\oplus\cdots\oplus\xi}_{m\text{ times}},
\end{equation}
where $\xi=\xi_{1}\cdots\xi_{ln}$.
\end{prop}
\begin{demo}{Proof}
From the definition \eqref{eq:zeta-b}, we obtain
\begin{align*}
\zeta^n
&=
\sum_{I^{(1)}_{1},\dots,I^{(m)}_{1}\in\binom{[ln]}{l}}
\cdots
\sum_{I^{(1)}_{n},\dots,I^{(m)}_{n}\in\binom{[ln]}{l}}
\cB(I^{(1)}_{1},\dots,I^{(m)}_{1})\cdots \cB(I^{(1)}_{n},\dots,I^{(m)}_{n})
\\&\qquad\qquad\qquad\qquad\times
\xi_{I^{(1)}_{1}}\cdots\xi_{I^{(1)}_{n}}\oplus\cdots\oplus\xi_{I^{(m)}_{1}}\cdots\xi_{I^{(m)}_{n}}.
\end{align*}
If we set $I^{(k)}=(I^{(k)}_{1},\dots,I^{(k)}_{n})$ for $1\leq k\leq m$,
then we have
\begin{align*}
\zeta^n
&=
\sum_{I^{(1)},\dots,I^{(m)}\in\binom{[ln]}{l,\dots,l}}
\sgn I^{(1)}\cdots\sgn I^{(m)}\,
\cB(I^{(1)}_{1},\dots,I^{(m)}_{1})\cdots \cB(I^{(1)}_{n},\dots,I^{(m)}_{n})
\\&\qquad\qquad\qquad\qquad\times
\xi\oplus\cdots\oplus\xi,
\end{align*}
because $\xi_{I^{(k)}_{1}}\cdots\xi_{I^{(k)}_{n}}=\sgn I^{(k)}\,\xi$ if 
$I^{(k)}_{1}\uplus\cdots\uplus I^{(k)}_{n}=[ln]$ (disjoint union), or $0$ otherwise for $1\leq k\leq m$.
Hence we obtain the desired identity.
\end{demo}
Let $l$, $r$ and $n$ be positive integers such that $l$ is even and $r\leq n$.
Let
\[
\cB=\left(\cB(I_1,\dots,I_{m})\right)_{I_1,\dots,I_{m}\in\binom{[ln]}{l}}
\]
be an $m$-dimensional array of size $ln$.
When $(S^{(1)},\dots,S^{(m)})$ is an $m$-tuple such that $S^{(k)}\in\binom{[ln]}{lr}$
for $1\leq k\leq m$,
we let $\cB_{S^{(1)},\dots,S^{(m)}}$ denote 
the $m$-dimensional $l$-block array
\[
(\cB(I_{1},\dots,I_{m}))_{I_{1}\in \binom{S^{(1)}}{l},\dots,I_{m}\in \binom{S^{(m)}}{l}}
\]
of size $lr$.
We call the hyperpfaffian
\begin{equation}
\Pf^{[l,m]}(\cB_{S^{(1)},\dots,S^{(m)}})
=\frac1{r!}\sum_{I^{(1)}\in\binom{S^{(1)}}{l,\dots,l},...,I^{(m)}\in\binom{S^{(m)}}{l,\dots,l}}
\prod_{k=1}^{m}\sgn(I^{(k)})
\prod_{j=1}^{r}
\cB\left(I^{(1)}_{j},\dots, I^{(m)}_{j}\right),
\label{eq:subhyperpfaffian}
\end{equation}
the \defterm{subhyperpfaffian} of $\cB$,
where we write $I^{(k)}=(I_{1}^{(k)},\dots,I_{r}^{(k)})$.
\par\smallskip
We need the following lemma for later use.
Let $D_{l}(\lambda)=\{l(\lambda-1)+1,l(\lambda-1)+2,\dots,l(\lambda-1)+l\}$ for a positive integer $\lambda$,
and let $D_{l}(\lambda_{1},\dots,\lambda_{n})=D_{l}(\lambda_{1})\uplus\cdots\uplus D_{l}(\lambda_{n})$
for $1 \leq \lambda_{1}<\cdots<\lambda_{n}\leq N$.
We may write $D_{l}(\bla)$ in short for $\bla=\{\lambda_{1},\dots,\lambda_{n}\}_{<}\in\binom{[N]}{n}$.
The following lemma is an easy consequence of \eqref{eq:subhyperpfaffian} and we omit the proof.
\begin{lemma}
\label{lem:subhyperpfaffian}
Let $l$, $m$, $n$ and $N$ be positive integers such that $l$ is even and $n\leq N$.
Let 
$\cB=(\cB(I^{(1)},\dots,I^{(m)}))_{I^{(1)},\dots,I^{(m)}\in\binom{lN}{l}}$
be an $m$-dimensional $l$-block array of size $lN$ defined by
\begin{equation}
\cB(I^{(1)},\dots,I^{(m)})
=\begin{cases}
1
&\text{ if $I^{(k)}=D_{l}(\lambda^{(k)})$ for some $\lambda^{(k)}\in [N]$ ($1\leq k\leq m$),}\\
0&\text{ otherwise.}
\end{cases}
\label{eq:Block}
\end{equation}
Let $S^{(1)},\dots,S^{(m)}\in\binom{[lN]}{ln}$.
Then we have
\begin{equation}
\Pf^{[l,m]}(\cB_{S^{(1)},\dots,S^{(m)}})
=\begin{cases}
1&\text{ if $S^{(k)}=D_{l}(\bla^{(k)})$
for some $\bla^{(k)}\in \binom{[N]}{n}$ ($1\leq k\leq m$),}\\
0&\text{ otherwise.}
\end{cases}
\label{eq:subhyper-deBruijn}
\end{equation}
\end{lemma}
For example,
if $l=2$ and $m=1$,
then $\cB$ is the array whose entry is
\[
\cB(I)=\begin{cases}
1\qquad&\text{ if $I=\{2\lambda-1,2\lambda\}$ for $1\leq \lambda\leq N$,}\\
0\qquad&\text{ otherwise.}
\end{cases}
\]
Assume $S\in\binom{[2N]}{2n}$.
Then $\Pf(\cB_{S})$ equals $1$ if $S$ is of the form 
$S=\{2\lambda_{1}-1,2\lambda_{1}\}\uplus \cdots \uplus \{2\lambda_{n}-1,2\lambda_{n}\}$
for some $1\leq \lambda_{1}<\dots<\lambda_{n}\leq N$,
or $0$ otherwise.
\subsection{Minor summation formulas of hyperpfaffians}\label{sub:hyperMSF}
The following theorem and its corollaries are generalizations of the minor summation
formulas of Pfaffians
(See \cite{IW1,IW2,Ma1,Ma2}).
The original formula corresponds to  the case where $(l, m, r) = (2, 2, 1)$ and $n$ is a positive integer.
That is to say,
if $\cB=(\cB(I))_{I\in\binom{[N]}{2}}=(\cB(\{i,j\}_{<}))_{\{i,j\}_{<}\in\binom{[N]}{2}}$
is a $2$-block array of size $N$
and $H=(H(i,j))_{1\leq i\leq 2n,1\leq j\leq N}$
a $2n\times N$ matrix,
then we have
\begin{equation*}
\sum_{J\in\binom{[N]}{2n}}
\Pf (\cB_{J})
\det H_{[2n],J}
=\Pf(Q)
\end{equation*}
where
\begin{equation*}
Q(I)=Q(\{i,j\}_{<})
=\sum_{K=\{k,l\}_{<}\in\binom{[N]}{2}}
\cB(\{k,l\}_{<})
\det\begin{pmatrix}
H(i,k)&H(i,l)\\
H(j,k)&H(j,l)
\end{pmatrix}
\end{equation*}
for $I=\{i,j\}_{<}\in \binom{[2n]}{2}$
(see \cite[Theorem~1]{IW1}, \cite[Lemma~2.1]{IW2}).
\begin{theorem}
\label{th:msf1}
Let $l$, $m$, $n$, $N$ and $r$ be positive integers such that $m\ge2$ and $l$ is even and $ln\leq N$.
Let 
$H(\nu)=(H(\nu)(i_{1},\dots,i_{m}))_{1\leq i_{1},\dots,i_{m-1}\leq ln,\,1\leq i_{m}\leq N}$
be an $m$-dimensional hypermatrix 
of size $(ln,\dots, ln, N)$ for $1\leq \nu\leq r$,
and let $\cB=(\cB(I^{(1)},\dots,I^{(r)}))_{I^{(1)},\dots,I^{(r)}\in \binom{[N]}{l}}$ 
be an $r$-dimensional $l$-block array of size $N$.
Then we have
\begin{align}
&
\sum_{S^{(1)},\dots,S^{(r)}\in\binom{[N]}{ln}}
\Pf^{[l,r]}(\cB_{S^{(1)},\dots,S^{(r)}})\,
\prod_{\nu=1}^{r}{\det}^{[m]} H(\nu)_{[ln],\dots,[ln],S^{(\nu)}}=\Pf^{[l,(m-1)r]}(Q),
\label{eq:msf}
\end{align}
where
$Q=(Q(I^{(1,1)},\dots,I^{(1,m-1)},\dots,I^{(r,1)},\dots,I^{(r,m-1)}))_{I^{(1,1)},\dots,I^{(r,m-1)}\in \binom{[ln]}{l}}$ 
is the $(m-1)r$-dimensional $l$-block array of size $ln$ defined by
\begin{align}
&
Q(I^{(1,1)},\dots,I^{(1,m-1)},\dots,I^{(r,1)},\dots,I^{(r,m-1)})
\nonumber\\&
=\sum_{K^{(1)},\dots,K^{(r)}\in\binom{[N]}{l}} \cB(K^{(1)},\dots,K^{(r)})\, 
\prod_{\nu=1}^{r}{\det}^{[m]}(H(\nu)_{I^{(\nu,1)},\dots,I^{(\nu,m-1)},K^{(\nu)}}).
\label{eq:msfQ}
\end{align}
\end{theorem}
\begin{demo}{Proof}
From the definition \eqref{eq:hyperpfaffian} of the hyperpfaffian,
we have
\begin{align*}
n!\Pf^{[l,(m-1)r]}(Q)
&=
\sum_{I^{(1,1)},\dots,I^{(1,m-1)},\dots,I^{(r,1)},\dots,I^{(r,m-1)}\in\binom{[ln]}{l,\dots,l}}
\prod_{\nu=1}^{r}\prod_{k=1}^{m-1}\sgn(I^{(\nu,k)})
\\&\times
\prod_{j=1}^{n}
Q(I^{(1.1)}_{j},\dots,I^{(1,m-1)}_{j},\dots,I^{(r,1)}_{j},\dots,I^{(r,m-1)}_{j}),
\end{align*}
where $I^{(\nu,k)}=(I^{(\nu,k)}_{1},\dots,I^{(\nu,k)}_{n})$ with  
$I_{j}^{(\nu,k)}\in\binom{[ln]}{l}$ %and $I^{(s,k)}=\bigcup_{j=1}^{n}I^{(s,k)}_{1}$
for $1\leq j\leq n$, $1\leq \nu\leq r$ and $1\leq k\leq m-1$.
Hence, by \eqref{eq:msfQ}, $n!\Pf^{[l,(m-1)r]}(Q)$ is equal to
\begin{align*}
&
\sum_{I^{(1,1)},\dots,I^{(1,m-1)}\in\binom{[ln]}{l,\dots,l}}
\cdots
\sum_{I^{(r,1)},\dots,I^{(r,m-1)}\in\binom{[ln]}{l,\dots,l}}
\prod_{\nu=1}^{r}\prod_{k=1}^{m-1}\sgn(I^{(\nu,k)})
\\&
\sum_{K^{(1)}_{1},\dots,K^{(r)}_{1}\in\binom{[N]}{l}}
\cdots
\sum_{K^{(1)}_{n},\dots,K^{(r)}_{n}\in\binom{[N]}{l}}
\prod_{j=1}^{n}\left\{
\cB(K^{(1)}_{j},\dots,K^{(r)}_{j})
\prod_{\nu=1}^{r}{\det}^{[m]}(H(\nu)_{I^{(\nu,1)}_{j},\dots,I^{(\nu,m-1)}_{j},K^{(\nu)}_{j}})
\right\}
\\&=
\sum_{K^{(1)}_{1},\dots,K^{(r)}_{1}\in\binom{[N]}{l}}
\cdots
\sum_{K^{(1)}_{n},\dots,K^{(r)}_{n}\in\binom{[N]}{l}}
\cB(K^{(1)}_{1},\dots,K^{(r)}_{1})\cdots \cB(K^{(1)}_{n},\dots,K^{(r)}_{n})
\\&\times
\sum_{I^{(1,1)},\dots,I^{(1,m-1)}\in\binom{[ln]}{l,\dots,l}}
\prod_{k=1}^{m-1}\sgn(I^{(1,k)})\prod_{j=1}^{n}{\det}^{[m]}(H(1)_{I^{(1,1)}_{j},\cdots,I^{(1,m-1)}_{j},K^{(1)}_{j}})
\\&\times\cdots\times
\sum_{I^{(r,1)},\dots,I^{(r,m-1)}\in\binom{[ln]}{l,\dots,l}}
\prod_{k=1}^{m-1}\sgn(I^{(r,k)})
\prod_{j=1}^{n}{\det}^{[m]}(H(r)^{}_{I^{(r,1)}_{j},\dots,I^{(r,m-1)}_{j},K^{(r)}_{j}}).
\end{align*}
By formula \eqref{eq:hyper-minor} in Lemma~\ref{lem:hyper-minor}, it reduces to
\begin{align*}
&\sum_{K^{(1)}_{1},\dots,K^{(r)}_{1}\in\binom{[N]}{l}}
\cdots
\sum_{K^{(1)}_{n},\dots,K^{(r)}_{n}\in\binom{[N]}{l}}
\cB(K^{(1)}_{1},\dots,K^{(r)}_{1})\cdots \cB(K^{(1)}_{n},\dots,K^{(r)}_{n})
\\&\qquad\times
{\det}^{[m]} H(1)_{[ln],\dots,[ln],K^{(1)}}\cdots {\det}^{[m]} H(r)_{[ln],\dots,[ln],K^{(r)}},
\end{align*}
where $K^{(\nu)}=(K^{(\nu)}_{1},\dots,K^{(\nu)}_{n})$ for $1\leq \nu\leq r$.
Since the hyperdeterminant vanishes unless $K^{(\nu)}_{1},\dots,K^{(\nu)}_{n}$ are mutually disjoint for each $\nu$,
we can set $S^{(\nu)}=K^{(\nu)}_{1}\uplus \cdots \uplus K^{(\nu)}_{n}\in\binom{[N]}{ln}$.
Hence we have 
$K^{(\nu)}=(K^{(\nu)}_{1},\dots,K^{(\nu)}_{n})\in\binom{S^{(\nu)}}{l,\dots,l}$
 for each $\nu$.
In this case,
the hyperdeterminant becomes
$
{\det}^{[m]} H(s)_{[ln],\dots,[ln],K^{(\nu)}}
=\sgn K^{(\nu)}
{\det}^{[m]} H(s)_{[ln],\dots,[ln],S^{(\nu)}}
$
and
the above sum becomes 
\begin{align}
&
\sum_{S^{(1)},\dots,S^{(r)}\in\binom{[N]}{ln}}
\sum_{K^{(1)}\in\binom{S^{(1)}}{l,\dots,l}}
\cdots
\sum_{K^{(r)}\in\binom{S^{(r)}}{l,\dots,l}}
\prod_{\nu=1}^{r}\sgn K^{(\nu)}
\nonumber\\&\qquad\qquad\qquad\qquad\times
\prod_{j=1}^{n}\cB\left(K^{(1)}_{j},\dots,K^{(r)}_{j}\right)
\prod_{\nu=1}^{r}{\det}^{[m]} H(\nu)_{[ln],\dots,[ln],S^{(\nu)}}.
\label{eq:minor-sum-primitive}
\end{align}
By \eqref{eq:subhyperpfaffian}, this reduces to
\begin{equation*}
n!\sum_{S^{(1)},\dots,S^{(r)}\in\binom{[N]}{ln}}
\Pf^{[l,r]}(\cB_{S^{(1)},\dots,S^{(r)}})\,
\prod_{\nu=1}^{r}{\det}^{[m]} H(\nu)_{[ln],\dots,[ln],S^{(\nu)}}.
\end{equation*}
This proves the theorem.
\end{demo}
If we put $m=2$ in Theorem~\ref{th:msf1},
then the hyperdeterminant becomes the ordinary determinant and we obtain the following corollary:
\begin{corollary}
\label{cor:msf-det}
Let $l$, $n$, $N$ and $r$ be positive integers such that $l$ is even and $ln\leq N$.
Let $H(\nu)=\left(h_{ij}(\nu)\right)_{1\le i\le ln, 1\le j\le N}$ 
be $ln\times N$ rectangular matrices for $1\leq \nu\leq r$,
and let $\cB=(\cB(I^{(1)},\dots,I^{(r)}))_{I^{(1)},\dots,I^{(r)}\in\binom{[N]}{l}}$ 
be an $r$-dimensional $l$-block array of size $N$.
Then we have
\begin{align}
&
\sum_{S^{(1)},\dots,S^{(r)}\in\binom{[N]}{ln}}
\Pf^{[l,r]}(\cB_{S^{(1)},\dots,S^{(r)}})
\prod_{\nu=1}^{r}\det H(\nu)_{[ln],S^{(\nu)}}
=\Pf^{[l,r]}(Q),
\label{eq:msf-det}
\end{align}
where
$Q=(Q(I^{(1)},\dots,I^{(r)}))_{I^{(1)},\dots,I^{(r)}\in \binom{[ln]}{l}}$ is 
the $r$-dimensional $l$-block array 
of size $ln$ defined by
\begin{equation}
Q(I^{(1)},\dots,I^{(r)})
=\sum_{K^{(1)},\dots,K^{(r)}\in\binom{[N]}{l}} \cB(K^{(1)},\dots,K^{(r)})\, 
\prod_{\nu=1}^{r}\det(H(\nu)_{I^{(s)},K^{(s)}}).
\label{eq:msf-det-Q}
\end{equation}
\end{corollary}
Actually the following theorem is not a direct corollary of Theorem~\ref{th:msf1},
but a consequence of the proof of Theorem~\ref{th:msf1}.
\begin{theorem}
\label{cor:Pf-Hf}
Let $l$, $n$, $N$ and $r$ be positive integers such that $l$ is even and $ln\leq N$.
Let $H(\nu)=(H_{ij}(\nu))_{1\leq i\leq ln,\,1\leq j\leq N}$
be an $ln\times N$ rectangular matrices for $1\leq \nu\leq r$,
and $\cA=(\cA(I))_{I\in \binom{[N]}{l}}$ a $1$-dimensional $l$-block array of size $N$.
If we define the array $Q=(Q(I^{(1)},\dots,I^{(r)}))_{I^{(1)},\dots,I^{(r)}\in \binom{[ln]}{l}}$
 by
\begin{equation}
Q(I^{(1)},\dots,I^{(r)})
=\sum_{K\in\binom{[N]}{l}} \cA(K)\, \prod_{\nu=1}^{r}\det(H(\nu)_{I^{(\nu)},K}),
\end{equation}
then we obtain the following identity:
\begin{equation}
\label{eq:Pf-Hf}
\Pf^{[l,r]}(Q)
=\begin{cases}
\sum_{S\in\binom{[N]}{ln}}
\Pf^{[l,1]}(\cA_{S}) \prod_{\nu=1}^{r}\det\left(H(\nu)_{[ln],S}\right)
&\text{ if $r$ is odd,}\\
\sum_{S\in\binom{[N]}{ln}}
\Hf^{[l,1]}(\cA_{S}) \prod_{\nu=1}^{r}\det\left(H(\nu)_{[ln],S}\right)
&\text{ if $r$ is even.}
\end{cases}
\end{equation}
\end{theorem}
\begin{demo}{Proof}
Recall the proof of Theorem~\ref{th:msf1} when $m=2$.
We put 
\[
\cB(I^{(1)},\dots,I^{(r)})=\cA(I^{(1)})\,\delta_{I^{(1)},I^{(2)}}\cdots\delta_{I^{(1)},I^{(r)}}
\]
for the given $\cA:\binom{[N]}{l}\to F$.
Then \eqref{eq:minor-sum-primitive} becomes 
\begin{equation*}
\sum_{S\in\binom{[N]}{ln}}
\sum_{K\in\binom{S}{l,\dots,l}}
(\sgn K)^{r}
\prod_{j=1}^{n}\cA\left(K_{j}\right)
\prod_{\nu=1}^{r}\det H(\nu)_{[ln],S}.
\end{equation*}
Note that $(\sgn K)^{r}=\sgn K$ if $r$ is odd, or $1$ if $r$ is even.
Hence we obtain the desired result.
\end{demo}
\par\smallskip
At the end of this section, 
we briefly sketch the proof of the following fact:
\begin{prop}\label{rem:equi-Pf}
\rm
Given an $m$-dimensional $l$-block array
\[
\cB=\left(\cB(I^{(1)},\dots,I^{(m)})\right)_{I^{(1)},\dots,I^{(m)}\in\binom{[ln]}{l}}
\]
of size $ln$,
we define the $1$-dimensional $lm$-block array
$
\cA=\left(\cA(J)\right)_{J\in\binom{[lmn]}{lm}}
$
of size $lmn$
by
\[
\cA(J)
=\begin{cases}
\cB(I^{(1)},\dots,I^{(m)})
&\text{ if $J=\bigcup_{k=1}^{m}(I^{(k)}+ln(k-1))$,}\\
0&\text{ otherwise,}
\end{cases}
\]
where $I^{(k)}+ln(k-1)$ stands for 
the subset $\{i+ln(k-1) \mid i\in I^{(k)}\}\in\binom{[lmn]}{l}$ for $1\leq k\leq m$.
Then we have
\begin{equation}
\Pf^{[l,m]}(\cB)=\Pf^{[lm,1]}(\cA).
\label{eq:apendix-matsu}
\end{equation}
\end{prop}
\begin{demo}{Proof}
The proof is similar to that of \cite[Proposition~A1]{Ma2}.
Given an $m$-tuple $(I^{(1)},\dots,I^{(m)})$ 
such that $I^{(k)}=(I^{(k)}_{1},\dots,I^{(k)}_{n})\in\binom{[ln]}{l,\dots,l}$ for $1\leq k\leq m$,
we associate $J=\left(J_{1},\dots,J_{n}\right)\in\binom{[lmn]}{lm,\dots,lm}$ 
with $J_{j}=\bigcup_{k=1}^{m}(I_{j}^{(k)}+ln(k-1))$ for $1\leq j\leq n$.
It is easy to check that $(I^{(1)},\dots,I^{(m)})\mapsto J$ gives a injection 
$\binom{[ln]}{l,\dots,l}^m \to \binom{[lmn]}{lm,\dots,lm}$,
and $\sgn J=\sgn I^{(1)}\cdots \sgn I^{(m)}$.
Further, from the definition of $\cA$, 
we have $\cA(J_{1})\cdots \cA(J_{n})=\prod_{j=1}^{n}\cB(I^{(1)}_{j},\dots,\sigma^{(m)}_{j})$
 if $J$ is in the image of this injection, 
and $0$ otherwise.
This proves the identity \eqref{eq:apendix-matsu}.
\end{demo}

\section{De Bruijn's formula and Hankel hyperpfaffians}
\label{sec:Bruijn}
De Bruijn's original paper contains two major Pfaffian formulas
\cite[(4,7)]{Bru} and \cite[(7.3)]{Bru},
both of which are called by his name.
In this section we establish certain hyperpfaffian analogues of both formulas,
see
Theorem~\ref{th:msf-integral} and Theorem~\ref{th:msf2-integral}.
Especially the latter implies % plays an important role in what follows.
 Corollary~\ref{cor:Hankel-hyperpfaff}, 
which %follows from Theorem~\ref{th:msf2-integral},
 gives a relation
between the Hankel hyperpfaffians and Selberg integrals.
However we are not completely satisfied with Corollary~\ref{cor:q-Hankel-hyperpfaff}
since it is hard to evaluate the integral on the right-hand side of \eqref{eq:gen_formula1}
unless $l=2$.
The main idea of the proof of those two theorems is to
use Theorem~\ref{th:msf1}.
The aim of this section is to reduce certain (hyper)pfaffians and their $q$-analogues,
which we call Hankel (hyper)pfaffians, to Selberg type integrals.
Note that $\Pf^{[l,1]}$ stands for Barvinok's original hyperpfaffian \cite{Bar},
where $l$ is a positive even integer.
\subsection{Jackson integral}
First we recall the definition of the Jackson integral.
Let $f(x)$ be a function defined on an interval $[\alpha,\beta]$.
The Jackson integral from $\alpha$ to $\beta$ is %, usually,
 defined by
\begin{equation*}
\int_{\alpha}^{\beta}f(x)\,d_{q}x
=\int_{0}^{\beta}f(x)\,d_{q}x-\int_{0}^{\alpha}f(x)\,d_{q}x
\end{equation*}
where
\begin{equation*}
\int_{0}^{\alpha}f(x)\,d_{q}x
=(1-q)\alpha\sum_{n=0}^{\infty}f(\alpha q^n)q^n
\end{equation*}
\cite[(1.11.2), (1.11.3)]{GR}.
Here we assume $0<q<1$ unless otherwise stated. All Jackson integrals and multiple sums below are assumed absolutely convergent, so that rearrangement and termwise integration are justified.
If $w$ is the weight function of a measure $d_{q}\omega$,
we have
\begin{equation}
\int_{0}^{\alpha}f(x)\,d_{q}\omega(x)
=\int_{0}^{\alpha}f(x)\,w(x)d_{q}x.
\label{eq:single-integral}
\end{equation}
More generally, we define the multiple $q$-integral by
\begin{equation}
\int\limits_{0\leq x_{1}<\dots<x_{n}\leq \alpha}f(\x)d_{q}\omega(\x)
=(1-q)^{n}\alpha^{n}\sum_{0\leq i_{1}< \cdots < i_{n}}
f(\alpha q^{i_{n}},\dots,\alpha q^{i_{1}})
\prod_{k=1}^{n}w(\alpha q^{i_{k}})
q^{i_{1}+\cdots+i_{n}}.
\label{eq:multiple-integral}
\end{equation}
Let $l$ and $m$ be positive integers,
and let $f:([0,\alpha]^{l})^{m}\rightarrow \C$ 
 be a function.
We say $f$ is \defterm{$l$-alternating} if it satisfies
$f(\sigma^{(1)}(\y^{(1)}),\dots,\sigma^{(m)}(\y^{(m)}))
=\sgn\sigma^{(1)}\cdots \sgn\sigma^{(m)} f(\y^{(1)},\dots,\y^{(m)})$,
where $\y^{(k)}=(y^{(k)}_{1},\dots,y^{(k)}_{l})\in [0,\alpha]^{l}$ 
and $\sigma^{(1)},\dots,\sigma^{(m)}\in\Sym_{l}$ 
for $1\leq k\leq m$.
This is equivalent to giving a function $\binom{[0,\alpha]}{l}^{m}\rightarrow \C$ such that 
$\beta(Y^{(1)},\dots,Y^{(m)})=f(\y^{(1)},\dots,\y^{(m)})$ 
where 
$Y^{(k)}=\{y^{(k)}_{1},\dots,y^{(k)}_{l}\}_{<}$ 
for $1\leq k\leq m$.
This means that it is enough to use only the values of 
an $l$-alternating function
$f(\y^{(1)},\dots,\y^{(m)})$
in the domain $y^{(k)}_{1}<\cdots<y^{(k)}_{l}$ 
($1\leq k\leq m$).
\subsection{The first de Bruijn type formula}
De Bruijn presented two types of formulas in his paper \cite{Bru}.
We establish hyperpfaffian versions of both.
We deduce the following theorem from Theorem~\ref{th:msf1},
which is a generalization of de Bruijn's formula \cite[(4,7)]{Bru}.
\begin{theorem}
\label{th:msf-integral}
Let $l$, $m$, $n$ and $r$ be positive integers such that $m\ge2$ and $l$ is even.
Let $\phi^{(\nu)}_{i_{1},\dots,i_{m-1}}(y)$ be a function on $[0,\alpha]$
for $i_{1},\dots,i_{m-1}\in[ln]$ and $1\leq \nu\leq r$.
Let $f(\y^{(1)},\dots,\y^{(r)})$ be an $l$-alternation function,
where $\y^{(\nu)}=(y^{(\nu)}_{1},\dots,y^{(\nu)}_{l})\in[0,\alpha]^{l}$ for $1\leq \nu\leq r$.
Then we have
\begin{align}
&\int_{0\leq x^{(1)}_{1}<\dots<x^{(1)}_{ln}\leq \alpha}\cdots\int_{0\leq x^{(r)}_{1}<\dots<x^{(r)}_{ln}\leq a}\,
\Pf^{[l,r]}\left(f(\x^{(1)}_{I^{(1)}},\dots,\x^{(r)}_{I^{(r)}})\right)_{I^{(1)},\dots,I^{(r)}\in \binom{[ln]}{l}}
\nonumber\\&\qquad\times
\prod_{\nu=1}^{r}{\det}^{[m]}\left(\phi^{(\nu)}_{j_{1},\dots,j_{m-1}}(x^{(\nu)}_{j_{m}})\right)_{1\leq j_{1},\dots,j_{m}\leq ln}
\,d_{q}\omega(\x^{(s)})
=\Pf^{[l,(m-1)r]}\left(Q\right),
\label{eq:msf-integral}
\end{align}
where
\begin{align}
&
Q(I^{(1,1)},\dots,I^{(1,m-1)},\dots,I^{(r,1)},\dots,I^{(r,m-1)})
=\int_{0\leq x^{(1)}_{1}<\cdots<x^{(1)}_{l}\leq \alpha}
\cdots
\int_{0\leq x^{(r)}_{1}<\cdots<x^{(r)}_{l}\leq \alpha}
\nonumber\\&\qquad
f(\x^{(1)},\dots,\x^{(r)})
\prod_{\nu=1}^{r} 
{\det}^{[m]}(\phi^{(\nu)}_{i^{(\nu,1)},\cdots,i^{(\nu,m-1)}}(x^{(\nu)}_{i_{m}}))_{i^{(\nu,1)}\in I^{(\nu,1)},\dots,i^{(\nu,m-1)}\in I^{(\nu,m-1)},i_{m}\in[l]}
\,d_{q}\omega(\x^{(\nu)})
\label{eq:deBruijn1-Q}
\end{align}
for $I^{(1,1)},\dots,I^{(1,m-1)},\dots,I^{(r,1)},\dots,I^{(r,m-1)}\in\binom{[ln]}{l}$.
Here we use the notation 
$\x^{(\nu)}=(x^{(\nu)}_{1},\\ \dots,x^{(\nu)}_{ln})$
and
$\x^{(\nu)}_{I^{(\nu)}}=(x^{(\nu)}_{i^{(\nu)}_{1}},\dots,x^{(\nu)}_{i^{(\nu)}_{l}})$
for $I^{(\nu)}=\{i^{(\nu)}_{1},\dots,i^{(\nu)}_{l}\}_{<}\in\binom{[ln]}{l}$ ($1\leq \nu\leq r$) in \eqref{eq:msf-integral}.
Meanwhile $\x^{(\nu)}$ stands for $(x^{(\nu)}_{1},\dots,x^{(\nu)}_{l})$
in \eqref{eq:deBruijn1-Q}.
\end{theorem}
\begin{demo}{Proof}
Let $N$ be an integer such that $N\geq ln$.
In this proof we use the notation 
$\alpha q^{I-1}=\{\alpha q^{i_{l}-1},\dots,\alpha q^{i_{1}-1}\}_{<}$
for any $l$-element set $I=\{i_{1},\dots,i_{l}\}_{<}\in \binom{[N]}{l}$.
We take 
\[
\cB(I^{(1)},\dots,I^{(r)})=f(\alpha q^{I^{(l)}-1},\dots,\alpha q^{I^{(r)}-1})
\]
for $I^{(1)},\dots,I^{(r)}\in\binom{[N]}{l}$,
and
\[
H(\nu)(i_{1},\dots,i_{m-1},i_{m})=(1-q)\alpha q^{i_{m}-1}\phi^{(\nu)}_{i_{1},\dots,i_{m-1}}(\alpha q^{i_{m}-1})\omega(\alpha q^{i_{m}-1})
\]
in Theorem~\ref{th:msf1}.
Then, by \eqref{eq:subhyperpfaffian},
the hyperpfaffian on the left-hand side of \eqref{eq:msf} becomes
\begin{align*}
\Pf^{[l,r]}(\cB_{S^{(1)},\dots,S^{(r)}})
&=\frac1{n!}\sum_{I^{(1)}\in\binom{S^{(1)}}{l,\dots,l},\dots,I^{(r)}\in\binom{S^{(r)}}{l,\dots,l}}
\prod_{\nu=1}^{r}\sgn(I^{(\nu)})
\prod_{j=1}^{n}
f\left(\alpha q^{I^{(1)}_{j}-1},\dots,\alpha q^{I^{(r)}_{j}-1}\right),
\end{align*}
where
$I^{(\nu)}$ is in the form $I^{(\nu)}=(I^{(\nu)}_{1},\dots,I^{(\nu)}_{n})$
with $I^{(\nu)}_{j}\in \binom{S^{(\nu)}}{l}$
($1\leq j\leq n$, $1\leq \nu\leq r$).
Meanwhile, by the definition \eqref{eq:hyperdeterminant}, 
the hyperdeterminant on the left-hand side of \eqref{eq:msf} equals
\begin{align*}
&{\det}^{[m]} H(\nu)_{[ln],\dots,[ln],S^{(\nu)}}
\\&
=\sum_{\sigma_{1},\dots,\sigma_{m-1}\in\frak{S}_{ln}}
\sgn(\sigma_{1}\cdots\sigma_{m-1})
\prod_{j=1}^{ln}H(\nu)(\sigma_{1}(j),\sigma_{2}(j),\dots,\sigma_{m-1}(j),s^{(\nu)}_{j})
\\&
=\sum_{\sigma_{1},\dots,\sigma_{m-1}\in\frak{S}_{ln}}
\sgn(\sigma_{1}\cdots\sigma_{m-1})\prod_{j=1}^{ln}
(1-q)\alpha q^{s^{(\nu)}_{j}-1}\phi^{(\nu)}_{\sigma_{1}(j),\dots,\sigma_{m-1}(j)}(\alpha q^{s^{(\nu)}_{j}-1})
\omega(\alpha q^{s^{(\nu)}_{j}-1}),
\end{align*}
where $S^{(\nu)}=\{s^{(\nu)}_{1},\dots,s^{(\nu)}_{ln}\}$.
If $N$ tends to $\infty$,
then it is not hard to see that the left-hand side of \eqref{eq:msf} becomes the left-hand side of \eqref{eq:msf-integral}
because of \eqref{eq:multiple-integral}.
The right-hand side of \eqref{eq:msf-integral} is more easy to see by putting $N \to \infty$ in 
\eqref{eq:msfQ} under this assumption.
This completes the proof of the theorem.
\end{demo}
\subsection{The second de Bruijn type formula and its corollaries}
The following theorem is a hyperpfaffian analogue of another de Bruijn's formula \cite[(7.3)]{Bru},
which will play an important role in this paper.
A special case,
i.e.,
$r=1$,
 is obtained in \cite[(96)]{LT1}.
\begin{theorem}
\label{th:msf2-integral}
Let $l$, $m$, $n$ and $r$ be positive integers such that $m\ge2$ and $l$ is even.
Let $\phi^{(\nu,k)}_{\i'}(x)=\phi^{(\nu,k)}_{i_{1},\dots,i_{m-1}}(x)$ be a function on $[0,\alpha]$
for $\i'=(i_{1},\dots,i_{m-1})\in[ln]^{m-1}$, $\nu\in[r]$ and $k\in[l]$.
Then we have
\begin{align}
&\int_{0\leq x^{(1)}_{1}<\dots<x^{(1)}_{n}\leq \alpha}\cdots\int_{0\leq x^{(r)}_{1}<\dots<x^{(r)}_{n}\leq \alpha}\,
\prod_{\nu=1}^{r}{\det}^{[m]}\left(\phi^{(\nu,1)}_{\i'}(x^{(\nu)}_{j})\Big|\cdots\Big|\phi^{(\nu,l)}_{\i'}(x^{(\nu)}_{j})\right)_{\i'\in[ln]^{m-1},\,j\in[n]}
\nonumber\\&\qquad\qquad
\times \prod_{\nu=1}^{r}d_{q}\omega(\x^{(\nu)})
=\Pf^{[l,(m-1)r]}\left(Q\right),
\label{eq:deBruijin}
\end{align}
where $d_{q}\omega(\x^{(\nu)})=d_{q}\omega(x^{(\nu)}_{1})\cdots d_{q}\omega(x^{(\nu)}_{n})$ and 
$Q$ 
is the $(m-1)r$-dimensional $l$-block array of size $ln$ defined by
\begin{align}
&
Q(I^{(1,1)},\dots,I^{(1,m-1)},\dots,I^{(r,1)},\dots,I^{(r,m-1)})
\nonumber\\&
=\int_{[0,\alpha]^r}
\prod_{\nu=1}^{r} 
{\det}^{[m]}\left(
\phi^{(\nu,k)}_{i^{(\nu,1)},\dots,i^{(\nu,m-1)}}(x^{(\nu)})
\right)_{i^{(\nu,1)}\in I^{(\nu,1)},\dots,i^{(\nu,m-1)}\in I^{(\nu,m-1)},k\in [l]}
\,d_{q}\omega(x^{(\nu)})
\label{eq:deBruijn-det-Q}
\end{align}
for $I^{(1,1)},\dots,I^{(r,m-1)}\in \binom{[ln]}{l}$.
Here
$
\left(\phi^{(\nu,1)}_{\i'}(y_{j})\Big|\cdots\Big|\phi^{(\nu,l)}_{\i'}(y_{j})\right)_{\i'\in[ln]^{m-1},\,j\in[n]}
$
stands for the $m$-dimensional hypermatrix of size $ln$ whose $\i'$th row is given by
\[
\left(\phi^{(\nu,1)}_{\i'}(y_{1}),\cdots,\phi^{(\nu,l)}_{\i'}(y_{1}),\dots,\phi^{(\nu,1)}_{\i'}(y_{n}),\cdots,\phi^{(\nu,l)}_{\i'}(y_{n})\right).
\]
\end{theorem}
\begin{demo}{Proof}
Let $N$ be a positive integer such that $N\geq n$.
We replace $N$ by $lN$ in Theorem~\ref{th:msf1},
and take the $r$-dimensional $l$-block array 
$\cB=(\cB(I^{(1)},\cdots,I^{(r)}))_{I^{(1)},\cdots,I^{(r)}\in\binom{lN}{l}}$
of size $lN$ defined by \eqref{eq:Block}
 in Lemma~\ref{lem:subhyperpfaffian}.
We take the $m$-dimensional hypermatrix 
$H(\nu)=\left(H(\nu)(\i',i_{m})\right)_{\i'\in [ln]^{m-1},i_{m}\in [lN]}$
whose entries are given by
\[
H(\nu)(\i',i_{m})
=\begin{cases}
\phi^{(\nu,\kappa)}_{\i'}\left(\alpha q^{\iota-1}\right)\omega\left(\alpha q^{\iota-1}\right)(1-q)\alpha q^{\iota-1}
&\text{ if $\kappa=1$,}\\
\phi^{(\nu,\kappa)}_{\i'}\left(\alpha q^{\iota-1}\right)
&\text{ if $2\leq \kappa\leq l$,}
\end{cases}
\]
where $\iota$ and $\kappa$ are the uniquely determined integers which satisfy 
$i_{m}=l(\iota-1)+\kappa$, $1\leq \iota \leq N$ and $1\leq \kappa\leq l$.
Then \eqref{eq:msfQ} reads 
\begin{align*}
&
Q(I^{(1,1)},\dots,I^{(1,m-1)},\dots,I^{(r,1)},\dots,I^{(r,m-1)})
\\&
=\sum_{\lambda^{(1)},\dots,\lambda^{(r)}\in[N]} 
\prod_{\nu=1}^{r}{\det}^{[m]}(H(\nu)_{I^{(\nu,1)},\dots,I^{(\nu,m-1)},D_{l}(\lambda^{(\nu)})})
\\&
=\sum_{\lambda^{(1)},\dots,\lambda^{(r)}\in[N]} 
\prod_{\nu=1}^{r}{\det}^{[m]}\left(
\phi^{(\nu,\kappa)}_{i^{(\nu,1)},\dots,i^{(\nu,m-1)}}(\alpha q^{\lambda^{(\nu)}-1})
\right)_{i^{(\nu,1)}\in I^{(\nu,1)},\dots,i^{(\nu,m-1)}\in I^{(\nu,m-1)},\kappa\in[l]}
\\&\qquad\qquad\qquad\qquad\times
\omega\left(\alpha q^{\lambda^{(\nu)}-1}\right)(1-q)\alpha q^{\lambda^{(\nu)}-1}
\end{align*}
for $I^{(1,1)},\dots,I^{(1,m-1)},\dots,I^{(r,1)},\dots,I^{(r,m-1)}\in\binom{[ln]}{l}$.
If we put $N\to\infty$ in this identity, we obtain \eqref{eq:deBruijn-det-Q} from \eqref{eq:single-integral}.
On the other hand, on the left-hand side of \eqref{eq:msf},
the hyperpfaffian vanishes by \eqref{eq:subhyper-deBruijn} 
unless each $S^{(\nu)}$ is in the form 
$S^{(\nu)}=D_{l}(\bla^{(\nu)})$
for some $\bla^{(\nu)}=(\lambda^{(\nu)}_{1},\dots,\lambda^{(\nu)}_{n})\in\binom{[N]}{n}$ ($1\leq \nu\leq r$).
In this case 
we have
\begin{align*}
{\det}^{[m]} H(\nu)_{[ln],\dots,[ln],D_{l}(\bla^{(\nu)})}
&
={\det}^{[m]}\left(\phi^{(\nu,1)}_{\i'}(\alpha q^{\lambda^{(\nu)}_{j}-1})\Big|\cdots\Big|\phi^{(\nu,l)}_{\i'}(\alpha q^{\lambda^{(\nu)}_{j}-1})\right)_{\i'\in[ln]^{m-1},\,j\in[n]}
\\&\times
\prod_{j=1}^{n}\omega \left(\alpha q^{\lambda^{(\nu)}_{j}-1} \right)(1-q)\alpha q^{\lambda^{(\nu)}_{j}-1}.
\end{align*}
Hence we can obtaon the left-hand side of \eqref{eq:deBruijin} from \eqref{eq:multiple-integral}
by putting $N\to\infty$.
Here we use the fact that $l$ is even to rearrange the size of the $\i'$th row as it fits \eqref{eq:multiple-integral}. 
\end{demo}
For example,
if $l=2$, $m=2$ and $r=1$,
then the above formula \eqref{eq:deBruijin} reads
\begin{align*}
\int_{0\leq x_{1}<\cdots<x_{n}\leq \alpha}
\det(\phi_{i}(x_{j})|\psi_{i}(x_{j}))_{i\in[2n],j\in[n]}d_{q}\omega(\x)
=\Pf\left(\int_{0}^{\alpha}\begin{vmatrix} \phi_{i_{1}}(x) & \psi_{i_{1}}(x)\\ \phi_{i_{2}}(x) & \psi_{i_{2}}(x) \end{vmatrix}\,d_{q}\omega(x)\right)_{i_{1},i_{2}\in[2n]},
\end{align*}
where $(\phi_{i}(x_{j})|\psi_{i}(x_{j}))_{i\in[2n],j\in[n]}$ stands for the matrix whose $i$th row equals
\[
\begin{pmatrix}
\phi_{i}(x_{1})&\psi_{i}(x_{1})&\hdots&\phi_{i}(x_{n})&\psi_{i}(x_{n})
\end{pmatrix}
\]
for $i\in[2n]$.
One can prove the following corollary directly from Theorem~\ref{cor:Pf-Hf}
by similar arguments.
The $r=1$ case of \eqref{eq:deBruijin2}
agrees with equation (96) in \cite{LT1}.
\begin{corollary}
\label{cor:deBruijn3}
Let $l$, $n$ and $r$ be positive integers such that $l$ is even.
Let $d_{q}\,\omega(x)=w(x)d_{q}x$ be a measure on $[0,\alpha]$,
and let $\phi^{(s,t)}_{i}(x)$ be a function on $[0,\alpha]$
for $i\in[ln]$ and $t\in[l]$.
Then we have
\begin{align}
&\int_{0\leq x_{1}<\dots<x_{n}\leq \alpha}
\prod_{\nu=1}^{r}\det\left(\phi^{(\nu,1)}_{i}(x_{j})\big|\cdots\big|\phi^{(\nu,l)}_{i}(x_{j})\right)_{i\in[ln],\,j\in[n]}
\,d_{q}\omega(\boldsymbol{x})
\nonumber\\&\qquad\qquad
=\begin{cases}
\Pf^{[l,r]}\left(Q(I^{(1)},\dots, I^{(r)})\right)_{I^{(1)},\dots,I^{(r)}\in\binom{[ln]}{l}},
&\text{ if $r$ is odd,}\\
\Hf^{[l,r]}\left(Q(I^{(1)},\dots,I^{(r)})\right)_{I^{(1)},\dots,I^{(r)}\in\binom{[ln]}{l}},
&\text{ if $r$ is even.}
\end{cases}
\label{eq:deBruijin2}
\end{align}
where
\begin{align}
Q(I^{(1)},\dots,I^{(r)})
=\int_{0}^{\alpha}
\prod_{\nu=1}^{r} 
\det(\phi^{(\nu,\mu)}_{i^{(\nu)}_{\lambda}}(x))_{1\leq \lambda,\mu\leq l}
\,\,d_{q}\omega(x)
\label{eq:deBruijn-det-Q2}
\end{align}
for $I^{(\nu)}=\{i^{(\nu)}_{1},\dots,i^{(\nu)}_{l}\}\in\binom{[ln]}{l}$ ($\nu\in[r]$).
\end{corollary}
Note that equation \eqref{eq:deBruijin2} generalizes 
Luque and Thibon's generalization \cite[(96)]{LT1} of de Bruijn's second formula.
From here we study how to apply this second de Bruijn type formula.
As in \cite{Hab1} we use the symbol
\begin{equation}
\label{eq:Delta1}
\Delta_{k}^{1}(\x)=\Delta_{k}^{1}(x_{1},\dots,x_{n})=\prod_{i<j}\prod_{\nu=0}^{k-1}(x_{j}-q^{\nu}x_{i})(x_{j}-q^{-\nu}x_{i})
\end{equation}
and $\displaystyle (q)_{k}=(q;q)_{k}$ in short.
Recall that $\Pf^{[l,1]}$ stands for Barvinok's original hyperpfaffian in the following corollary.
\begin{corollary}
\label{cor:q-Hankel-hyperpfaff}
Let $l$ and $n$ be positive integers such that $l$ is even.
Let $t$ be an integer,
and let $\mu_{i}=\int_{0}^{\alpha}x^id_{q}\omega(x)$ denote the $i$th moment of the measure $\omega$.
Then we have
\begin{align}
&\Pf^{[l,1]}\Bigl(
\prod_{1\leq j<k\leq l}(q^{i_{j}-1}-q^{i_{k}-1})
\cdot
\mu_{\sum_{k=1}^{l}i_{k}+t-l}
\Bigr)_{1\leq i_{1}<i_{2}<\cdots< i_{l} \leq ln}
\nonumber\\&
=
q^{n\binom{l}3+l\binom{l}2\binom{n}2}
\prod_{k=1}^{l}(q)_{k-1}^{n}
\int_{0\leq x_1<\dots<x_n\leq \alpha}
\prod_{i=1}^{n}x_i^{t+\binom{l}2}
\prod_{1\leq i<j\leq n}(x_{j}-x_{i})^{-l}
\prod_{k=1}^{l}\Delta_{k}^{1}(\x)
\,d_{q}\omega(\x),
\label{eq:gen_formula1}
\end{align}
where $d_{q}\omega(\x)=d_{q}\omega(x_1) \cdots d_{q}\omega(x_n)$.
\end{corollary}
\begin{demo}{Proof}
The proof mainly appeals to
the Vandermonde determinant
\begin{equation}
\Delta_{N}(X)=\det(X_{j}^{i-1})_{1\leq i,j\leq N}=\prod_{1\leq i<j\leq N}(X_{j}-X_{i}).
\label{eq:Vandermonde}
\end{equation}
We put $m=2$ and $r=1$ in Theorem~\ref{th:msf2-integral}
(or $r=1$ in Corollary~\ref{cor:deBruijn3}).
Here we can write $\phi^{(k)}_{i}$ for $\phi^{(1,k)}_{i}$
and $I$ for $I^{(1)}$ in short in \eqref{eq:deBruijin} and \eqref{eq:deBruijn-det-Q}.
We take 
\[
\phi^{(k)}_{i}(x)=\begin{cases}
q^{(l-k)(i-1)}x^{i-1+t}
&\text{ if $k=1$,}\\
q^{(l-k)(i-1)}x^{i-1}
&\text{ if $2\leq k\leq l$,}
\end{cases}
\]
in \eqref{eq:deBruijn-det-Q},
then we obtain
$\det(\phi^{(k)}_{i_{j}}(x))_{j,k\in[l]}
=\det((q^{i_{j}-1})^{l-k})_{j,k\in[l]}\,x^{|I|-l+t}$,
which implies
\begin{equation*}
Q \left( I \right)
=\prod_{1\leq j<k\leq l}
(q^{i_{j}-1}-q^{i_{k}-1})
\int_{0}^{\alpha}x^{|I|-l+t}\,d_{q}\omega(x)
=\prod_{1\leq j<k\leq l}(q^{i_{j}-1}-q^{i_{k}-1})
\,\mu_{|I|-l+t}
\end{equation*}
from \eqref{eq:Vandermonde} with $N=l$.
On the other hand,
if we perform the same substitution as before
on the left-hand side of \eqref{eq:deBruijin}, 
we find that 
\begin{align*}
&\det\left(
\phi^{(1)}_{i}(x_{j})|\cdots|\phi^{(l)}_{i}(x_{j})
\right)_{i\in [ln],\,j\in [n]}
\\&\qquad=\prod_{j=1}^{n}x_{j}^{t}\cdot\det\left(
q^{(l-1)(i-1)}x_{j}^{i-1}\,|\,\cdots\,|\,q^{i-1}x_{j}^{i-1}\,|\,x_{j}^{i-1}
\right)_{i\in [ln],\,j\in [n]},
\end{align*}
where the $i$th row of the matrix equals
\[
(q^{(l-1)(i-1)}x_{1}^{i-1},\dots,q^{i-1}x_{1}^{i-1},x_{1}^{i-1} , \dots\dots , q^{(l-1)(i-1)}x_{n}^{i-1},\dots,q^{i-1}x_{n}^{i-1},x_{n}^{i-1}),
\]
which can be realized as the Vandermonde determinant \eqref{eq:Vandermonde} with $N=ln$
as follows.
Each  integer $i\in [ln]$  can be written in a unique way as
$i=l(j-1)+r$ for $j\in[n]$ and $r\in[l]$. 
If we choose  $X_{i}=q^{l-r}x_{j}$ for $i\in [ln]$ with 
$i=l(j-1)+r$ ($j\in [n]$ and $r \in [l]$), 
then the within-cluster and between-cluster factors give
\begin{align*}
 \Delta_{ln}(X)
 &=q^{n\binom l3+l\binom l2\binom n2}
   \prod_{k=0}^{l-1}(q;q)_k^n\prod_{j=1}^n x_j^{\binom l2}\\
 &\quad\times\prod_{i<j}\left\{(x_j-x_i)^l
       \prod_{r=1}^{l-1}(x_j-q^rx_i)^{l-r}
                         (x_j-q^{-r}x_i)^{l-r}\right\}.
\end{align*}
Combining this factorization with the preceding determinant integral gives
\eqref{eq:gen_formula1}, since the between-cluster product is
$\prod_{i<j}(x_j-x_i)^{-l}\prod_{k=1}^l\Delta_k^1(\x)$.
\end{demo}
The factorization above gives a general integral representation for
$q$-Hankel hyperpfaffians. The product evaluations below use its ordinary
integral counterpart or the case $l=2$. Finding a choice of the functions
$\phi_i^{(k)}$ leading to tractable $q$-integrals for $l>2$ remains a
natural problem.
The ordinary integral counterpart admits a direct proof using a confluent
Vandermonde determinant. This avoids imposing a $q$-dependent family of
measures and also includes measures with atoms.
\begin{corollary}
\label{cor:Hankel-hyperpfaff}
Let $l$ and $n$ be positive integers such that $l$ is even,
 and $t$ an integer.
Let $d\psi$ be a finite signed or positive Borel measure on $[0,\alpha]$,
and write $\mu_i=\int_0^\alpha x^i\,d\psi(x)$. Assume $t+\binom l2\ge0$,
so every moment index below is nonnegative. (The same formula applies to
other integer shifts whenever all terms are absolutely integrable.)
Use $d\psi(\x)=\prod_{j=1}^n d\psi(x_j)$. Then we have
\begin{align}
&\Pf^{[l,1]}\Bigl(
\prod_{1\leq j<k\leq l}(i_{k}-i_{j})\cdot\mu_{\sum_{k=1}^{l}i_{k}-l+t}
\Bigr)_{1\leq i_{1}<\cdots< i_{l} \leq ln}
\nonumber\\&\qquad\qquad
=\frac{\prod_{k=1}^{l}\{(k-1)!\}^{n}}{n!}\int_{[0,\alpha]^n}
\prod_{i}x_i^{t+\binom{l}{2}}\prod_{i<j}(x_{j}-x_{i})^{l^2}
\,d\psi(\x).
\label{eq:hp-psi}
\end{align}
\end{corollary}
\begin{demo}{Proof}
Set $C=\binom l2$, $c_l=\prod_{r=0}^{l-1}r!$ and
$D_r=(r!)^{-1}(d/dx)^r$. For $i_1<\cdots<i_l$,
\[
 \det\bigl(D_{r-1}x^{i_p-1}\bigr)_{p,r=1}^l
 =\frac{\prod_{p<q}(i_q-i_p)}{c_l}
 x^{\sum_p i_p-l-C}.
\]
This follows because the falling factorials of degrees $0,\ldots,l-1$
are monic polynomials in the exponents $i_p-1$.
The full normalized confluent Vandermonde determinant is
\[
 \det\bigl(D_rx_j^{i-1}\bigr)_
 {\substack{1\le i\le ln\\1\le j\le n,\ 0\le r<l}}
 =\prod_{i<j}(x_j-x_i)^{l^2},
\]
where columns are ordered first by $j$ and then by $r$.
Multiply by $\prod_jx_j^{t+C}$ and integrate. Laplace expansion into
blocks of $l$ columns gives a sum over ordered set partitions; by the
hyperpfaffian normalization this sum is $n!c_l^{-n}$ times the
left-hand side of \eqref{eq:hp-psi}. This proves the formula.
Finite expansion and absolute integrability justify each step, including
for signed measures and measures with atoms.
\end{demo}

\begin{corollary}\label{cor:stieltjes-preservation}
Let $d\psi$ be a positive finite measure on $[0,\alpha]$, and fix even
$l$ and $n\ge1$. For $r\ge0$ define
\[
 H_r=\Pf^{[l,1]}\left(\prod_{p<q}(i_q-i_p)
      \mu_{\sum_p i_p-l+r}\right)_{1\le i_1<\cdots<i_l\le ln}.
\]
Then $(H_r)_{r\ge0}$ is a Stieltjes moment sequence. In particular,
$H_rH_{r+2}\ge H_{r+1}^2$ for every $r\ge0$.
\end{corollary}
\begin{demo}{Proof}
By \eqref{eq:hp-psi}, $H_r$ is the $r$th moment of the pushforward,
under $(x_1,\ldots,x_n)\mapsto x_1\cdots x_n$, of the positive measure
\[
 \frac{c_l^n}{n!}\prod_jx_j^{\binom l2}
 \prod_{i<j}(x_j-x_i)^{l^2}\prod_jd\psi(x_j),
 \qquad c_l=\prod_{j=0}^{l-1}j!.
\]
The support of the pushforward lies in $[0,\alpha^n]$.
The inequality follows from Cauchy--Schwarz.
\end{demo}

\def\INT{\mathrsfs{I}}
\def\ConstB#1{\tau_{#1}}
\def\ConstC#1{C_{#1}}
\def\ConstD#1{\eta_{#1}}
\def\ConstH#1{H_{#1}}
\section{Hankel hyperpfaffians for Narayana polynomials}
\label{sec:Motzkin}
This section deals with an application of Corollary~\ref{cor:Hankel-hyperpfaff}
to concrete hyperpfaffian computations.
In \cite[Conjecture 6.2]{ITZ2} we presented a conjecture 
on Pfaffian identities 
involving Motzkin, Delannoy, Schr\"oder numbers and Narayana polynomials.
We define the Narayana polynomials of type $A$, $B$ and $D$,
which provide us a unified treatment of these combinatorial numbers.
The master theorem of this section is Theorem~\ref{th:tildeN}.
The proof of this theorem shows when we can apply \eqref{eq:Selberg} or \eqref{eq:Aomoto}
to the hyperpfaffians involving the Narayana polynomials of type $A$, $B$ and $D$.
As corollaries
we derive remarkable Pfaffian identities from the master theorem.
\subsection{Definitions and main results}
For $n\ge1$ in type $A$, $n\ge0$ in type $B$, and $n\ge2$ in type $D$,
we introduce \defterm{Narayana numbers} 
of type $A$, $B$ and $D$ by
\[
N_{k}(A_n)=\frac1{n}\binom{n}{k}\binom{n}{k-1},
\quad
N_{k}(B_n)=\binom{n}{k}^2,
\quad
N_{k}(D_n)=\binom{n}{k}\left\{\binom{n-1}{k}+\binom{n-2}{k-2}\right\}.
\]
The \defterm{Narayana polynomials} $N(X_{n},a)$ ($n\geq0$) are defined by
\begin{equation}
N(X_{n},a)=\sum_{k=0}^{n}N_{k}(X_{n})a^{k}
\label{eq:Narayana}
\end{equation}
for $X=A,B$ or $D$,
where
we use the convention that $N(A_{0},a)=N(D_{0},a)=1$, $N(D_{1},a)=\frac{a+1}2$
(see \cite[pp.~277--278]{P}).
\par\smallskip
For convenience we introduce the notation 
\begin{equation}
\Phi_{n}(r,s,m)=\prod_{j=1}^{n}
\frac{\binom{2m(j-1)+2r}{m(j-1)+r}\binom{2m(j-1)+2s}{m(j-1)+s}}
{\binom{2m(j-1)+r+s}{m(j-1)+r}}
\prod_{k=2}^{n}
\frac{\binom{mk}{m}}
{\binom{m(2k-3)+r+s}{m(k-1)}}
\label{eq:Phi}
\end{equation}
for nonnegative integers $r,s,m$ and $n$.
This is related to the value of Selberg integral
(see Lemma~\ref{lem:lem-s}).
In this paper we take $\sqrt{a}=\sqrt{r} e^{i\theta/2}$ for a complex number $a=re^{i\theta}$ 
($r\geq0$ and $-\pi<\theta\leq \pi$).
The following is the master theorem of this section,
and all the succeeding Pfaffian identities in this section follow from this theorem.
\begin{theorem}
\label{th:tildeN}
Let $l$ and $n$ be positive integers such that $l$ is even.
Let $r$ be an integer such that $r\geq-\binom{l}2$ (hence the entries of each hyperpfaffian are well-defined),
and 
let $a\in \C\setminus \Real_{\leq0}$. %%be any complex number which is not a negative real number.
Let 
$\ConstB{a}=\frac{(\sqrt{a}-1)^2}{4\sqrt{a}}$ 
and 
$\ConstH{l,n}=\frac{1}{n!}\prod\limits_{k=1}^{l}\{(k-1)!\}^{n}$.
Let $\omega$ be either primitive cube root of unity.
In the $a=1$ branch of \eqref{eq:ThD}, impose the additional condition
$r\ge1-\binom l2$. Then the following identities hold:
\begin{align}
&\Pf^{[l,1]}\Bigl(
\prod_{1\leq j<k\leq l}(i_{k}-i_{j})\cdot 
N\left(A_{|I|+r-l},a\right)
\Bigr)_{I \in \binom{[ln]}{l}}
\nonumber\\&
=\begin{cases}
2^{-n}\ConstH{l,n}
\cdot
\Phi_{n}\left(r+\binom{l}2,1,\frac{l^2}2\right)
&\text{ if $(a,r)=(1,r)$,}\\
2^{-n}a^{n+\frac{l^2}2\binom{n}2}\ConstH{l,n}
\Phi_{n}\left(1,1,\frac{l^2}2\right)
&\text{ if $(a,r)=\left(a,1-\binom{l}{2}\right)$,}\\
2^{n}a^{\frac32n+\frac{l^2}2\binom{n}2}\ConstH{l,n}
\Phi_{n}\left(1,1,\frac{l^2}2\right)
\sum\limits_{k=0}^{n}
\binom{n}{k}
\ConstB{a}^{n-k}
\prod\limits_{j=1}^{k}
\frac{3+l^2(n-j)}{6+l^2(2n-j-1)}
&\text{ if $(a,r)=\left(a,2-\binom{l}{2}\right)$.}
\end{cases}
\label{eq:ThA}
\end{align}
\begin{align}
&\Pf^{[l,1]}\Bigl(
\prod_{1\leq j<k\leq l}(i_{k}-i_{j})\cdot 
N\left(B_{|I|+r-l},a\right)
\Bigr)_{I \in \binom{[ln]}{l}}
\nonumber\\&
=\begin{cases}
\ConstH{l,n}
\Phi_{n}\left(r+\binom{l}2,0,\frac{l^2}2\right)
&\text{ if $(a,r)=(1,r)$,}\\
a^{\frac{l^2}2\binom{n}2}\ConstH{l,n}
\Phi_{n}\left(0,0,\frac{l^2}2\right)
&\text{ if $(a,r)=\left(a,-\binom{l}{2}\right)$,}\\
2^{2n}a^{\frac{n}2+\frac{l^2}2\binom{n}2}\ConstH{l,n}
\Phi_{n}\left(0,0,\frac{l^2}2\right)
\sum\limits_{k=0}^{n}
\binom{n}{k}
\ConstB{a}^{n-k}
\prod\limits_{j=1}^{k}\frac
{1+l^2(n-j)}
{2+l^2(2n-j-1)}
&\text{ if $(a,r)=\left(a,1-\binom{l}{2}\right)$.}
\end{cases}
\label{eq:ThB}
\end{align}
\begin{align}
&\Pf^{[l,1]}\Bigl(
\prod_{1\leq j<k\leq l}(i_{k}-i_{j})\cdot 
N(D_{|I|+r-l},a)
\Bigr)_{I \in \binom{[ln]}{l}}
\nonumber\\
&=\begin{cases}
2^{2n}\ConstH{l,n}
\Phi_{n}\left(r+\binom{l}2-1,0,\frac{l^2}2\right)
\sum\limits_{k=0}^{n}\binom{n}{k}\left(-\frac14\right)^{n-k}
\prod\limits_{j=1}^{k}
\frac
{2r-l-1+l^2(n-j+1)}
{2r-l+l^2(2n-j)}
\\
\hskip300pt\text{ if $(a,r)=(1,r)$,}
\\
2^{2n}\omega^{n+\frac{l^2}2\binom{n}2}\ConstH{l,n}
\Phi_{n}\left(1,0,\frac{l^2}2\right)
\sum\limits_{k=0}^{n}\binom{n}{k}\left(-\frac58\right)^{n-k}
\prod\limits_{j=1}^{k}
\frac{3+l^2(n-j)}
{4+l^2(2n-j-1)}
\\
\hskip270pt\text{ if $(a,r)=\left(\omega,2-\binom{l}2\right)$.}
\end{cases}
\label{eq:ThD}
\end{align}
\end{theorem}
The sums on the right-hand side come from Aomoto integral \eqref{eq:Aomoto2}.
We first note that most of the famous combinatorial numbers are related to
the Narayana polynomials, as detailed below.
\par\medbreak
\label{prop:specialization}
Let $n$ be a nonnegative integer,
and let $\omega=\frac{-1\pm\sqrt{-3}}2$,
which is a primitive cube root of unity.
By specializing $a$ in \eqref{eq:Narayana},
we obtain several classical numbers:
the \defterm{Catalan numbers}
$\Cat(n)=N(A_n,1)=\frac1{2n+1}\binom{2n+1}{n}$,
the \defterm{large Schr\"oder numbers} $\Sch(n)=N(A_n,2)=\sum_{k=0}^{n}\binom{n+k}{2k}\Cat(k)$,
the \defterm{central binomial coefficients} $\CBC(n)=N(B_{n},1)=\binom{2n}{n}$,
the \defterm{central Delannoy numbers} $\Del(n)=N(B_{n},2)=\sum_{k=0}^{n}\binom{n}{k}\binom{n+k}{k}$.
We call the sequence $\Cat^{D}(n)=N(D_{n},1)=(3n-2)\Cat(n-1)$ ($n\geq1$) 
the \defterm{Catalan numbers of type} $D$,
which have several combinatorial meanings (see \cite[12.3]{P}, A051924).
Further we have the \defterm{Motzkin numbers}
$
\Mot(n)=(-1)^{n}\omega^{n+2}
N\left(A_{n+1},\omega\right)
=\sum_{k=0}^{\lfloor n/2 \rfloor}\binom{n}{2k}\Cat(k)
$
for $n\geq0$,
and
the \defterm{central trinomial coefficients}
$
\CTC(n)=(-1)^{n}\omega^{n}
N\left(B_{n},\omega\right)
$
for $n\geq0$,
which is the coefficient of $x^n$ in the expansion of $(1+x+x^2)^n$.
Let ${}_{2}F_{1}(a,b;c;x)=\sum_{n=0}^{\infty}\frac{(a)_{n}(b)_{n}}{n!(c)_{n}}x^n$
denote the hypergeometric series 
with Pochhammer symbol $(x)_{n}=\frac{\Gamma(x+n)}{\Gamma(x)}$.
Finally,
the \defterm{Motzkin numbers of type} $D$
are defined by
\[
\Mot^{D}(n)=(-1)^{n}\omega^{n}
N\left(D_{n},\omega\right)
={}_{2}F_{1}\left(\frac{1-n}2,1-\frac{n}2;1;4\right)
+(n-2){}_{2}F_{1}\left(1-\frac{n}2,\frac32-\frac{n}2;2;4\right)
\]
for $n\geq2$
(A298300).
These specializations are 
roughly summarized
in Table~\ref{tab:specialization}.
\begin{table}[htb]
\begin{center}
\begin{tabular}{c|ccccccc}
\hline
Type & $A$ & $B$ & $D$  \\
\hline
$a=1$      & $\Cat(n)$ & $\CBC(n)$ & $\Cat^{D}(n)$  \\
$a=2$      & $\Sch(n)$ & $\Del(n)$ &   \\ 
$a=\omega$ & $\Mot(n)$ & $\CTC(n)$ & $\Mot^{D}(n)$ \\ 
\hline
\end{tabular}
\caption{Specializations of Narayana polynomials\label{tab:specialization}}
\end{center}
\end{table}
\par\smallskip
By putting $l=2$ and substituting appropriate values into $a$ and $r$ in Theorem~\ref{th:tildeN},
we obtain the following identities conjectured in \cite[Conjecture~6.2]{ITZ2}
for ordinary Pfaffians.
\begin{corollary}
\label{cor:motzkin}
Let $n$ be a positive integer.
Then the following identities hold:
\begin{align}
&\Pf\biggl((j-i)\,\Mot(i+j-3)\biggr)_{1\leq i,j\leq 2n}
=\prod_{k=0}^{n-1}(4k+1),
\label{eq:conj-M_n}
\\
&\Pf\biggl((j-i)\,\Del(i+j-3)\biggr)_{1\leq i,j\leq 2n}
=2^{n^2-1}(2n-1)\prod_{k=1}^{n-1}(4k-1),
\label{eq:conj-D_n}
\\
&\Pf\biggl((j-i)\,\Sch(i+j-2)\biggr)_{1\leq i,j\leq 2n}
=2^{n^2}\prod_{k=0}^{n-1}(4k+1),
\label{eq:conj-S_n}
\\
&\Pf\biggl((j-i)\,N(A_{i+j-2},a)\biggr)_{1\leq i,j\leq 2n}
=a^{n^2}\prod_{k=0}^{n-1}(4k+1).
\label{eq:conj-Nara_n}
\end{align}
\end{corollary}
In fact,
by substituting other values into $a$ and $r$ of Theorem~\ref{th:tildeN}, 
we obtain several more remarkable Pfaffian identities.
In Table~\ref{tab:Derivation} the reader can see how to specialize the parameters
to obtain the identities in the corollaries above and below.
\begin{table}[htb]
\begin{center}
\begin{tabular}{c|ccccccc}
\hline
Identity & Type & Master Identity & $l$ & $a$ & $r\in\Z$  \\
\hline
\eqref{eq:conj-M_n} & $A$ & \eqref{eq:ThA} & $2$ & $\omega$ & $0$ \\
\eqref{eq:conj-D_n} & $B$ & \eqref{eq:ThB} & $2$ & $2$ & $-1$ \\ 
\eqref{eq:conj-S_n} & $A$ & \eqref{eq:ThA} & $2$ & $2$ & $0$ \\ 
\eqref{eq:Mot_sum}  & $A$ & \eqref{eq:ThA} & $2$ & $\omega$ & $1$ \\ 
\eqref{eq:Del_sum}  & $B$ & \eqref{eq:ThB} & $2$ & $2$ & $0$ \\ 
\eqref{eq:Sch_sum}  & $A$ & \eqref{eq:ThA} & $2$ & $2$ & $1$ \\ 
\eqref{eq:MotD_sum} & $D$ & \eqref{eq:ThD} & $2$ & $\omega$ & $1$ \\ 
\eqref{eq:Cat_arb}  & $A$ & \eqref{eq:ThA} & $2$ & $1$ & $r\geq-1$ \\ 
\eqref{eq:CBC_arb}  & $B$ & \eqref{eq:ThB} & $2$ & $1$ & $r\geq-1$ \\ 
\eqref{eq:CatD_arb} & $D$ & \eqref{eq:ThD} & $2$ & $1$ & $r\geq0$ \\ 
\hline
\end{tabular}
\caption{Derivation of the Pfaffian identities\label{tab:Derivation}}
\end{center}
\end{table}
\begin{corollary}
Let $n$ be a positive integer.
Then %we have
\footnotesize
\begin{align}
&\Pf\biggl((j-i)\,\Mot(i+j-2)\biggr)_{1\leq i,j\leq 2n}
=2^{2n}\prod_{k=0}^{n-1}(4k+1)
\sum_{k=0}^{n}\binom{n}{k}\left(-\frac14\right)^{n-k}
\prod_{j=1}^{k}
\frac{3+4(n-j)}{2+4(2n-j)},
\label{eq:Mot_sum}
\\
&\Pf\biggl((j-i)\,\Del(i+j-2)\biggr)_{1\leq i,j\leq 2n}
\nonumber\\&\qquad=2^{n\left(n+\frac52\right)-1}
(2n-1)\prod_{k=1}^{n-1}(4k-1)
\sum_{k=0}^{n}\binom{n}{k}\left(\frac{3\sqrt{2}-4}{8}\right)^{n-k}
\prod_{j=1}^{k}
\frac{1+4(n-j)}{4(2n-j)-2},
\label{eq:Del_sum}
\\
&\Pf\biggl((j-i)\,\Sch(i+j-1)\biggr)_{1\leq i,j\leq 2n}
=2^{n\left(n+\frac52n\right)}
\prod_{k=0}^{n-1}(4k+1)
\sum_{k=0}^{n}\binom{n}{k}\left(\frac{3\sqrt{2}-4}{8}\right)^{n-k}
\prod_{j=1}^{k}
\frac{3+4(n-j)}{2+4(2n-j)},
\label{eq:Sch_sum}
\\
&\Pf\biggl((j-i)\,\Mot^{D}(i+j-1)\biggr)_{1\leq i,j\leq 2n}
=2^{2n}
\prod_{k=0}^{n-1}\frac{(4k+2)!(4k)!}{(2k)!(2k+2n-1)!}
\sum_{k=0}^{n}\binom{n}{k}\left(-\frac{5}{8}\right)^{n-k}
\prod_{j=1}^{k}
\frac{3+4(n-j)}{4(2n-j)}.
\label{eq:MotD_sum}
\end{align}
\normalsize
\end{corollary}
Further, the case of $a=1$ in Theorem~\ref{th:tildeN}
gives us the following formulas:
\begin{corollary}\label{cor:arb-r}
Let $n$ be a positive integer,
and $r\geq-1$ be an integer.
Then %we have
\begin{align}
&\Pf\left((j-i)\,\Cat(i+j+r-2)\right)_{1\leq i,j\leq 2n}
=2^{-n}
\prod_{k=0}^{n-1}
\frac{(4k+2r+2)!(4k+2)!}{(2k+r+1)!(2k+2n+r)!},
\label{eq:Cat_arb}
\\
&\Pf\left((j-i)\,\CBC(i+j+r-2)\right)_{1\leq i,j\leq 2n}
=
\prod_{k=0}^{n-1}
\frac{(4k)!(4k+2r+2)!(2k+1)}{(2k+r+1)!(2k+2n+r-1)!},
\label{eq:CBC_arb}
\\
&\Pf\left((j-i)\Cat^{D}(i+j+r-2)\right)_{1\leq i,j\leq 2n}
\nonumber\\&
=
\frac{2^n}{n!}\prod_{k=0}^{n-1}
\frac{(4k+2r)!(4k)!(2k+2)!}{(2k+r)!(2k)!(2n+2k+r-2)!}
\sum_{k=0}^{n}\binom{n}{k}\left(-\frac14\right)^{n-k}
\prod_{j=1}^{k}\frac{2r+1+4(n-j)}{2r-2+4(2n-j)}.
\label{eq:CatD_arb}
\end{align}
Here we assume $r\geq0$ in \eqref{eq:CatD_arb}
(see the definition of $\Cat^{D}$).
\end{corollary}
The \defterm{Hankel Pfaffian transform}
of a sequence $\{a_{n}\}_{n\geq0}$ is defined to be 
the sequence $\{\mathcal{P}_{n}\}_{n\geq0}$ 
of Pfaffians $\mathcal{P}_{n}=\Pf\left((j-i)a_{i+j-r}\right)_{1\leq i,j\leq 2n}$.
Corollary~\ref{cor:arb-r} gives
the Hankel Pfaffian transforms of
the sequences $\{\Cat(n+r)\}_{n\geq0}$ and $\{\CBC(n+r)\}_{n\geq0}$
for a fixed $r$.
We summarize some identified sequences with OEIS references 
of \eqref{eq:Cat_arb} and \ref{eq:CBC_arb} in Table~\ref{tab:arb_r}.
\begin{table}[htb]
\begin{center}
\begin{tabular}{c|ccccccc}
\hline
Identity & $r=-1$ & $r=0$ & $r=1$ & $r=2$ \\
\hline
\eqref{eq:Cat_arb} & A000407 & A007696 & A001813 & A007696 \\
\eqref{eq:CBC_arb} &         &         & A147626 &         \\ 
\hline
\end{tabular}
\caption{Hankel Pfaffian transforms of Catalan numbers and CBCs\label{tab:arb_r}}
\end{center}
\end{table}
\par\smallskip
As pointed in \cite{ITZ2},
we can  also derive
\eqref{eq:Cat_arb} and \eqref{eq:CBC_arb}
from \eqref{eq:ITZ}.
It is interesting to note that 
\eqref{eq:Cat_arb}
is 
a Pfaffian analogue of Desainte-Catherine and Viennot's determinant identity:
\begin{equation}
\det(\Cat(i+j+r-1))_{1\leq i,j\leq n}=\prod_{1\leq i\leq j\leq r}\frac{i+j+2n}{i+j},
\label{eq:Desainte-Catherine-Viennot}
\end{equation}
see \cite{DV,Ta}.
The reader may think that the product on the right-hand side of \eqref{eq:Cat_arb}
is not so simple as that of \eqref{eq:Desainte-Catherine-Viennot}.
But, if one takes the ratio of the products between $r$ and $r+1$,
it looks very similar to the ratio of \eqref{eq:Desainte-Catherine-Viennot}.
 As  shown in \cite{DV}, the determinant of \eqref{eq:Desainte-Catherine-Viennot}  
 is the number of Young tableaux with entries 
 from $[r]$ satisfying: 
 there are at most $2n$ rows; 
 the rows are strictly increasing; 
 the columns are non-decreasing; 
 every column has an even number of cells. 
It would be a challenging problem to find combinatorial
objects which count some of the above Pfaffians.
\subsection{Proof of Theorem~\ref{th:tildeN}}
Now we proceed to prove Theorem~\ref{th:tildeN}.
First 
we need to find the generating function 
\[
G(X,a,z)=\sum_{n=0}^{\infty}N(X_{n},a)z^{n}
\]
of the Narayana polynomials \eqref{eq:Narayana} 
for $X=A,B$ or $D$.
\begin{lemma}
\label{lem:GF-Narayana}
We have
\begin{align}
G(A,a,z)
&=\frac{1-(a-1)z-\sqrt{(a-1)^2z^2-2(a+1)z+1}}{2z},
\label{eq:gf-narayana-A}
\\%%%%%%%%%%
G(B,a,z)
&=\frac{1}{\sqrt{(a-1)^2z^2-2(a+1)z+1}},
\label{eq:gf-narayana-B}
\\%%%%%%%%%%
G(D,a,z)
&=\frac12\left\{
\sqrt{(a-1)^2z^2-2(a+1)z+1}+
\frac{1+(a+1)z}{\sqrt{(a-1)^2z^2-2(a+1)z+1}}
\right\}.
\label{eq:gf-narayana-D}
\end{align}
\end{lemma}
\begin{demo}{Proof}
Let $y$ be the function at the right-hand side
of \eqref{eq:gf-narayana-A} 
multiplied by $z$.
Then
\[
y^2+\{(a-1)z-1\}y+z=0
\]
and $\lim\limits_{z\to 0}\frac{y}{z}=1$.
This equation can be rewritten as 
$y=z\left(1+\frac{ay}{1-y}\right)$.
By Lagrange inversion, for $n\ge1$,
\begin{align*}
 [z^{n+1}]y
 &=\frac1{n+1}[u^n]\left(1+\frac{au}{1-u}\right)^{n+1}\\
 &=\frac1{n+1}\sum_{k=1}^n\binom{n+1}{k}\binom{n-1}{k-1}a^k
 =\frac1n\sum_{k=1}^n\binom nk\binom n{k-1}a^k.
\end{align*}
For $n=0$ the coefficient is $1$.
This establishes \eqref{eq:gf-narayana-A}.
\par\smallskip
The Legendre polynomials $P_n(x)$ have the explicit formula
 (see \cite[p. 162]{R}):
\[
P_n(x)
=\left(\frac{x-1}{2}\right)^n\sum_{k=0}^n {n\choose k}^2\left(\frac{x+1}{x-1}\right)^k.
\]
Substituting $a=\frac{x+1}{x-1}$,
we obtain
$$
 \sum_{k=0}^n{n\choose k}^2a^k=(a-1)^n P_n\left(\frac{a+1}{a-1}\right).
$$
Equation \eqref{eq:gf-narayana-B} then follows from 
the following generating function
(see \cite[Chapter 10]{R})
\[
(1-2xt+t^2)^{-1/2}=\sum_{n=0}^\infty P_n(x)t^n.
\]
\par
By definition we have
\begin{align*}
N_{k}(D_{n})&=\binom{n}{k}^2-\binom{n}{k}\binom{n-2}{k-1}
=N_{k}(B_{n})-n N_{k}(A_{n-1})
\end{align*}
for $n\geq 2$, %and
$N(A_{0},a)=N(B_{0},a)=N(D_{0},a)=1$, 
$N(B_{1},a)=a+1$ and $N(D_{1},a)=\frac{a+1}2$.
Since
\[
(zG(A,a,z))'=\sum_{n\geq 0} (n+1)N(A_n,a)z^n=\sum_{n\geq 1} n N(A_{n-1},a)z^{n-1},
\]
we have $G(D,a,z)=G(B,a,z)-(zG(A,a,z))'\cdot z-\frac{a-1}2z$,
which immediately implies \eqref{eq:gf-narayana-D} 
by using \eqref{eq:gf-narayana-A} and \eqref{eq:gf-narayana-B}.
\end{demo}
A positive moment representation must not be confused with quasi-definiteness.
Nonvanishing leading Hankel determinants characterize a quasi-definite
moment functional; positive representing measures instead require positive
semidefiniteness of all finite Hankel matrices. We verify the required
representations explicitly. This also identifies an atom in type $A$ and
avoids an invalid unshifted density in type $D$.

\begin{lemma}\label{lem:measure}
Let $a>0$, $b_\pm=(1\pm\sqrt a)^2$, $I_a=[b_-,b_+]$, and
$R_a(x)=\sqrt{(b_+-x)(x-b_-)}$ on $I_a$. The type $A$ and $B$
Narayana polynomials are moments of the positive measures
\begin{align}
 d\mu_A(x)&=(1-a)_+\delta_0(dx)
       +\frac{R_a(x)}{2\pi x}\,\mathbf1_{(b_-,b_+)}(x)\,dx,
       \label{eq:psi-A}\\
 d\mu_B(x)&=\frac{\mathbf1_{(b_-,b_+)}(x)}{\pi R_a(x)}\,dx.
       \label{eq:psi-B}
\end{align}
Here $(1-a)_+=\max(1-a,0)$, and $\int x^k\,d\mu_X=N(X_k,a)$
for $k\ge0$ and $X=A,B$. Put $P_a(x)=2x^2-(a+1)x+(a-1)^2$.
For type $D$, the shifted representation is
\begin{equation}\label{eq:psi-D}
 N(D_{k+2},a)=\frac1{2\pi}\int_{I_a}\frac{x^kP_a(x)}{R_a(x)}\,dx,
 \qquad k\ge0.
\end{equation}
At $a=1$, one may use the lower shift
\begin{equation}\label{eq:D-one-shift}
 N(D_{k+1},1)=\frac1\pi\int_0^4
                  \frac{x^k(x-1)}{\sqrt{x(4-x)}}\,dx,\qquad k\ge0.
\end{equation}
The measures in the two shifted type $D$ formulas may be signed.
\end{lemma}
\begin{demo}{Proof}
Choose $Q(z)=\sqrt{(z-b_-)(z-b_+)}$ analytic off $I_a$ and
normalized by $Q(z)/z\to1$ at infinity. The generating functions give
\[
 \frac1zG(A,a,1/z)=\frac{z-(a-1)-Q(z)}{2z},\qquad
 \frac1zG(B,a,1/z)=\frac1{Q(z)}.
\]
The boundary jumps across $I_a$ give the continuous densities in
\eqref{eq:psi-A}--\eqref{eq:psi-B}. If $a\ne1$, then
$Q(0)=-|a-1|$, so the residue of the first transform at zero is
$(1-a)_+$. The second transform has no atom. At $a=1$ the displayed
densities are integrable and there is no atom. Equivalently, these
Stieltjes transforms determine all coefficients at infinity and hence
all moments.

The generating functions also give, for $k\ge0$,
\begin{equation}\label{eq:D-B-shift}
 N(D_{k+2},a)=N(B_{k+2},a)-\frac{a+1}{2}N(B_{k+1},a)
                  +\frac{(a-1)^2}{2}N(B_k,a).
\end{equation}
Inserting \eqref{eq:psi-B} proves \eqref{eq:psi-D}.
When $a=1$, the identity
$N(D_{k+1},1)=N(B_{k+1},1)-N(B_k,1)$ also holds for $k\ge0$;
this follows from the generating functions, including the convention
$N(D_1,1)=1$. It proves \eqref{eq:D-one-shift}.
\end{demo}

The atom in \eqref{eq:psi-A} contributes only to the zeroth moment.
In type $D$, dividing the density in \eqref{eq:psi-D} by $x^2$
does not in general give a representing finite measure for the unshifted
sequence. At $a=1$ it would be
$(x-1)/(\pi x\sqrt{x(4-x)})$, which is not integrable at zero.
The shifted formulas above are the ones used below.

The following lemma simplifies the values of the Selberg integrals
in Theorem~\ref{th:tildeN}.
\begin{lem}
\label{lem:lem-s}
Let $r,s,m$ and $n$ be nonnegative integers.
Then %we have 
\begin{equation}
S_n\left(r+\frac{1}{2}, s+\frac{1}{2}, {m}\right)
=\frac{\pi^n}{2^{2n\{m(n-1)+r+s\}}}
\cdot
\Phi_{n}(r,s,m)
\label{eq:lem-s}
\end{equation}
\end{lem}
\begin{demo}{Proof} 
We use $\Gamma(x+n)=\Gamma(x)(x)_{n}$ for nonnegative integer $n$.
Then the Gamma function formulas
$\Gamma\left(r+m(j-1)+\frac12\right)=\Gamma\left(\frac12\right)\left(\frac12\right)_{m(j-1)+r}$
and $\Gamma\left(\frac12\right)=\sqrt{\pi}$ lead to the desired identity by direct computation.
\end{demo}
\begin{demo}{Proof of Theorem~\ref{th:tildeN}}
Let $\INT_{A}=\INT_{A}(a,r)$, $\INT_{B}=\INT_{B}(a,r)$ and $\INT_{D}=\INT_{D}(a,r)$  
 be the left-hand side of \eqref{eq:ThA}, \eqref{eq:ThB}
and \eqref{eq:ThD}, respectively.
Write $s=r+\binom l2$. We first take $a>0$.
Using Corollary~\ref{cor:Hankel-hyperpfaff} with the shifted moments of
Lemma~\ref{lem:measure}, we obtain the following formulas in the ranges:
type $A$, $s\ge1$ (also $s=0$ when $a=1$); type $B$, $s\ge0$;
type $D$, $s\ge2$ (also $s=1$ when $a=1$).
For example, in type $D$ apply the corollary to the moment sequence
$N(D_{k+2},a)$ with shift $t=r-2$; the $a=1,s=1$ case uses
\eqref{eq:D-one-shift} instead. All these integrals are absolutely convergent.
\begin{align*}
\INT_{A}(a,r)
&=\frac{\ConstH{l,n}}{(2\pi)^{n}}\int_{I_{a}^n}
\prod_{i=1}^{n} x_i^{r+\binom{l}{2}-1}
\sqrt{4a-(x_{i}-a-1)^2}
\,\prod_{i<j}(x_{j}-x_{i})^{l^2}
\,d\x,
\\
\INT_{B}(a,r)
&=\frac{\ConstH{l,n}}{\pi^{n}}\int_{I_{a}^n}
\prod_{i=1}^{n}
\frac{ x_i^{r+\binom{l}{2}}}
{\sqrt{4a-(x_{i}-a-1)^2}}
\,\prod_{i<j}(x_{j}-x_{i})^{l^2}
\,d\x,
\\
\INT_{D}(a,r)
&=\frac{\ConstH{l,n}}{(2\pi)^{n}}\int_{I_{a}^n}
\prod_{i=1}^{n}
\frac{x_i^{r+\binom{l}{2}-2}\left\{2x_{i}^2-(a+1)x_{i}+(a-1)^2\right\}}
{\sqrt{4a-(x_{i}-a-1)^2}}
\,\prod_{i<j}(x_{j}-x_{i})^{l^2}
\,d\x,
\end{align*}
where $I_{a}=\left[(\sqrt{a}-1)^2,(\sqrt{a}+1)^2\right]$. 
If we set $x_{i}=4\sqrt{a}\,t_{i}+(\sqrt{a}-1)^2$ 
 in the above integrals,
then a straightforward computation leads to the following identities:
\begin{align}
\INT_{A}(a,r)
&=
C^{A}_{n,r}
\int_{\left[0,1\right]^n}
\prod_{i=1}^{n} \left\{ t_i+\frac{(\sqrt{a}-1)^{2}}{4\sqrt{a}}\right\}^{r+\binom{l}{2}-1}
\sqrt{t_{i}(1-t_{i})}
\,\prod_{i<j}(t_{j}-t_{i})^{l^2}
\,d\boldsymbol{t},
\label{eq:I_A}
\\
\INT_{B}(a,r)
&=C^{B}_{n,r}\int_{\left[0,1\right]^n}
\prod_{i=1}^{n}
\frac{ \left\{ t_i+\frac{(\sqrt{a}-1)^{2}}{4\sqrt{a}}\right\}^{r+\binom{l}{2}}}
{\sqrt{t_{i}(1-t_{i})}}
\prod_{i<j}(t_{j}-t_{i})^{l^2}
\,d\boldsymbol{t},
\label{eq:I_B}
\\
\INT_{D}(a,r)
&=
C^{D}_{n,r}
\int_{\left[0,1\right]^n}
\frac{
\prod_{i=1}^{n} \left\{ t_i+\frac{(\sqrt{a}-1)^{2}}{4\sqrt{a}}\right\}^{r+\binom{l}{2}-2}
 \left( t_i-\ConstD{a}^{+}\right) \left( t_i-\ConstD{a}^{-}\right)
\,\prod_{i<j}(t_{j}-t_{i})^{l^2}
}
{\prod_{i=1}^{n}\sqrt{t_{i}(1-t_{i})}}
\,d\boldsymbol{t},
\label{eq:I_D}
\end{align}
where 
$C^{A}_{n,r}=\frac{\ConstH{l,n}}{(2\pi)^{n}}\left(4\sqrt{a}\right)^{\left\{r+\binom{l}{2}+1\right\}n+l^2\binom{n}2}$,
$C^{B}_{n,r}=C^{D}_{n,r}=\frac{\ConstH{l,n}}{\pi^{n}}\left(4\sqrt{a}\right)^{\left\{r+\binom{l}{2}\right\}n+l^2\binom{n}2}$
and
$\ConstD{a}^{\pm}=\frac{-3(a+1)+8\sqrt{a}\pm\sqrt{4a-7(a-1)^2}}{16\sqrt{a}}$.
\par\smallskip
Equations \eqref{eq:I_A}--\eqref{eq:I_D} are initially asserted only
in the ranges just specified. For $s\ge1$ in type $A$, $s\ge0$ in
type $B$, and $s\ge2$ in type $D$, their remaining factors in the
$t_i$ are polynomials. Their integrals are therefore finite combinations
of beta-type integrals with coefficients analytic in $\sqrt a$ on
$\mathbb C\setminus\mathbb R_{\le0}$. In type $D$ the product
$(t-\eta_a^+)(t-\eta_a^-)$ is used as a quadratic polynomial, so its
individual roots introduce no extra branch choice. The identity theorem
extends the formulas from $a>0$ to this slit plane. The exceptional
$a=1$ endpoint shifts are evaluated directly using the convergent
representations of Lemma~\ref{lem:measure}.
\par\smallskip
Now we determine when these integrals are as in the form of the Selberg integral \eqref{eq:Selberg} or \eqref{eq:Aomoto2}.
\par\smallbreak
First, assume $(a,r)=(1,r)$. Then we obtain
\begin{align*}
\INT_{A}(1,r)
&=
C^{A}_{n,r}
\int_{\left[0,1\right]^n}
\prod_{i=1}^{n} 
t_i^{r+\binom{l}{2}-1/2}
(1-t_{i})^{1/2}
\,\prod_{i<j}(t_{j}-t_{i})^{l^2}
\,d\boldsymbol{t}
\\&
=C^{A}_{n,r}\cdot
S_{n}\left(r+\binom{l}{2}+\frac12,\frac32,\frac{l^2}2\right),
\\
\INT_{B}(1,r)
&=C^{B}_{n,r}\int_{\left[0,1\right]^n}
\prod_{i=1}^{n}
t_i^{r+\binom{l}{2}-1/2}(1-t_{i})^{-1/2}
\prod_{i<j}(t_{j}-t_{i})^{l^2}
\,d\boldsymbol{t}
\\&
=C^{B}_{n,r}\cdot
S_{n}\left(r+\binom{l}{2}+\frac12,\frac12,\frac{l^2}2\right),
\\
\INT_{D}(1,r)
&=
C^{D}_{n,r}
\int_{\left[0,1\right]^n}
\prod_{i=1}^{n} 
 \left( t_i-\frac14\right)
t_i^{r+\binom{l}{2}-3/2}(1-t_{i})^{-1/2}
\prod_{i<j}(t_{j}-t_{i})^{l^2}
\,d\boldsymbol{t}
\\&
=C^{D}_{n,r}\sum_{k=0}^{n}\left(-\frac14\right)^{n-k}
\int_{\left[0,1\right]^n}
 e_{k}(\boldsymbol{t})
\prod_{i=1}^{n} 
t_i^{r+\binom{l}{2}-3/2}(1-t_{i})^{-1/2}
\prod_{i<j}(t_{j}-t_{i})^{l^2}
\,d\boldsymbol{t}
\\&
=C^{D}_{n,r}\cdot
S_{n}\left(r+\binom{l}{2}-\frac12,\frac12,\frac{l^2}2\right)
\sum_{k=0}^{n}\binom{n}{k}\left(-\frac14\right)^{n-k}
\prod_{j=1}^{k}\frac{r+\binom{l}{2}-\frac12+\frac{l^2}2(n-j)}
{r+\binom{l}{2}+\frac{l^2}2(2n-j-1)}
\end{align*}
from \eqref{eq:Selberg} and \eqref{eq:Aomoto2}.
Hence,
using \eqref{eq:lem-s}, we obtain
 the desired identities,
 i.e., \eqref{eq:ThA},  \eqref{eq:ThB} and \eqref{eq:ThD}, when $(a,r)=(1,r)$.
\par\smallbreak
Hereafter we may assume $a\neq1$.
If $r=1-\binom{l}{2}$ then
\eqref{eq:I_A} reduces to \eqref{eq:Selberg}.
Hence we obtain
\begin{equation}
\INT_{A}\left(a,1-\binom{l}{2}\right)
=C^{A}_{n,r}
\,
S_{n}\left(\frac32,\frac32,\frac{l^2}2\right)
=2^{-n}\ConstH{l,n}a^{n+\frac{l^2}2\binom{n}2}
\Phi_{n}\left(1,1,\frac{l^2}2\right),
\end{equation}
which is equal to the second case of \eqref{eq:ThA}.
If we put $r=2-\binom{l}{2}$ in \eqref{eq:I_A} then
we obtain 
$
\prod_{i=1}^{n} (t_i+\ConstB{a})^{r+\binom{l}{2}-1}=\sum_{k=0}^{n}e_{k}(t)\ConstB{a}^{n-k}.
$
Hence, by \eqref{eq:Aomoto2} we obtain
\begin{align}
\INT_{A}\left(a,2-\binom{l}{2}\right)
=
C^{A}_{n,r}
\,
S_{n}\left(\frac32,\frac32,\frac{l^2}{2}\right)
\sum_{k=0}^{n}\binom{n}{k}\ConstB{a}^{n-k}
\prod_{j=1}^{k}\frac{\frac32+\frac{l^2}{2}(n-j)}
{3+\frac{l^2}{2}(2n-j-1)},
\end{align}
which proves the third case of \eqref{eq:ThA}.
\par\smallbreak
Next, for $X=B$,
if we put $r=-\binom{l}{2}$ in \eqref{eq:I_B} then we obtain
\begin{equation}
\INT_{B}\left(a,-\binom{l}{2}\right)
=C^{B}_{n,r}
\,
S_{n}\left(\frac12,\frac12,\frac{l^2}2\right)
=\ConstH{l,n}a^{\frac{l^2}2\binom{n}2}
\Phi_{n}\left(0,0,\frac{l^2}2\right).
\end{equation}
which proves the second case of \eqref{eq:ThB}.
If we put $r=1-\binom{l}{2}$ in \eqref{eq:I_B}, then we have
\begin{align}
\INT_{B}\left(a,1-\binom{l}{2}\right)
=
C^{B}_{n,r}
\,
S_{n}\left(\frac12,\frac12,\frac{l^2}{2}\right)
\sum_{k=0}^{n}\binom{n}{k}\ConstB{a}^{n-k}
\prod_{j=1}^{k}\frac{\frac12+\frac{l^2}{2}(n-j)}
{1+\frac{l^2}{2}(2n-j-1)},
\end{align}
which agrees with the third case of 
 \eqref{eq:ThB}.
\par\smallskip
Finally, we consider the $X=D$ case.
When $a\neq1$, \eqref{eq:I_D} can be in the form of
\eqref{eq:Selberg} or \eqref{eq:Aomoto2}
only if $r+\binom{l}{2}-2=0$ and one of $\ConstD{a}^{+}$ or $\ConstD{a}^{-}$
equals $0$ or $1$.
This actually happens if and only if $a=\omega$ or $a=\omega^{-1}$
in which case we have $\ConstD{a}^{-}=0$ and $\ConstD{a}^{+}=\frac58$.
Hence we obtain
\begin{align*}
&\INT_{D}\left(\omega^{\pm1},2-\binom{l}{2}\right)
=
C^{D}_{n,r}
S_{n}\left(\frac32,\frac12,\frac{l^2}2\right)
\sum_{k=0}^{n}\binom{n}{k}\left(-\frac58\right)^{n-k}
\prod_{j=1}^{k}\frac{\frac32+\frac{l^2}2(n-j)}{2+\frac{l^2}2(2n-j-1)}
\\&=
2^{2n}\omega^{n+\frac{l^2}2\binom{n}2}\ConstH{l,n}
\Phi_{n}\left(1,0,\frac{l^2}2\right)
\sum_{k=0}^{n}\binom{n}{k}\left(-\frac58\right)^{n-k}
\prod_{j=1}^{k}\frac{\frac32+\frac{l^2}2(n-j)}{2+\frac{l^2}2(2n-j-1)}
\end{align*}
by \eqref{eq:Aomoto2}.
This completes the proof of the theorem.
\end{demo}
\begin{remark}
If one writes the product in \eqref{eq:I_A}, \eqref{eq:I_B} or \eqref{eq:I_D}
 as a linear combination of the Jack polynomials 
 and uses \cite[Theorem~1,~2]{Kad2} or \cite[Corollary~1.2]{Wa},
then he/she can obtain more formulas of the hyperpfaffians 
for other values of $r$.
For an explicit finite insertion formula with all coefficients given by
products, see Theorem~\ref{app-thm:powers} and
Corollary~\ref{app-cor:quadratic}.
\end{remark}
In this section we studied only the Narayana polynomials
of type $A$, $B$ and $D$ since they cover the most of famous combinatorial numbers
as we noted.
The common feature of the combinatorial numbers of this section is that the generating functions
of the sequences
 are always solutions of quadratic equations as we saw in Lemma~\ref{lem:GF-Narayana}. 
This makes us easier to find the distribution  function and evaluate the related Selberg-Aomoto integral.
Of course,
we can choose other moment sequences.
However it is not always easy to find the distribution function.
Conjecture~\ref{conj:Gessel-Xin} in the last section is such an example.

\section[Hankel Pfaffians and little q-Jacobi polynomials]{Hankel Pfaffians and little $q$-Jacobi polynomials}
\label{sec:Askey}
The aim of this section is 
to give a new proof of \eqref{eq:ITZ},
i.e.,
 Theorem~\ref{conj:01}, \eqref{eq:pf-special}.
This $q$-identity first appeared in \cite[Corollary~3.2]{ITZ2},
and was also proved in \cite{GITZ1}, \cite[Corollary 3.4]{GITZ2}
by using a quadratic formula for the basic hypergeometric series related to 
the Askey-Wilson polynomials.
Here we use a different approach from the above two proofs,
that is,
we establish a key identity \eqref{eq:Pf-Delta2}
which enables us to compute $q$-Pfaffians of this type.
This key identity will be used to derive %the main theorem,
Theorem~\ref{conj:01},
and also used to derive Theorem~\ref{conj:03} in the next section.
\par\smallskip
As mentioned just before Corollary~\ref{cor:Hankel-hyperpfaff},
another case that
we can evaluate 
\eqref{eq:gen_formula1}
is when $l=2$.
Then we call the left-hand side of \eqref{eq:gen_formula1} \defterm{$q$-Hankel Pfaffian}.
We give examples of $q$-Hankel Pfaffians
in this section and the next one.
In this section we reduce the evaluation of the $q$-Hankel Pfaffian \eqref{eq:pf-special} 
to the $k=2$ case of the Askey-Habsieger-Kadell integral \eqref{eq:ask-conj1}.
\par\smallbreak
This section is composed as follows.
In \cite{Hab1,Kad} the $q$-Selberg integral has several expressions 
which include one of the $q$-difference products 
$\Delta_{k}(\x)$, $\Delta_{k}^{0}(\x)$, $\Delta_{k}^{1}(\x)$ and $\Delta_{k}^{2}(\x)$ 
(to be defined later).
To prove Theorem~\ref{conj:01},
we first state the relations between the integrals
which include the $q$-difference products
(see \eqref{eq:Delta1-Delta2}, \eqref{eq:Delta_k-Delta_k^1} and \eqref{eq:Delta1-Delta0}).
These identities hold for any $q$-measure $\omega$ and any $l=2k$.
Then we set $l=2$ to obtain the Pfaffian identity \eqref{eq:Pf-Delta2},
which will be an appropriate form to apply the Askey-Habsieger-Kadell formula
\eqref{eq:ask-conj1}.
Finally we take the measure which gives the little $q$-Jacobi polynomials,
and perform a straightforward computation to prove Theorem~\ref{conj:01}.
\par\smallbreak
If we put $l=2$ in \eqref{eq:gen_formula1},  
we obtain
\begin{equation}
\Pf\Bigl(
(q^{i-1}-q^{j-1})\mu_{i+j+r-2}
\Bigr)_{1\leq i<j\leq2n}
=\frac{q^{n(n-1)}(1-q)^{n}}{n!}
\int_{[0,a]^{n}}
\prod_{i=1}^{n}x_{i}^{r+1}
\cdot
\Delta_{2}^{1}(\x)
\,d_q\omega(\x).
\label{eq:r-mom-jacobi}
\end{equation}
In this section we use the $q$-analogue of the Selberg integral formula
\eqref{eq:ask-conj1}
to evaluate this Pfaffian when $\mu_n=\frac{(aq;q)_{n}}{(abq^2;q)_{n}}$ is 
the $n$th moment of the little $q$-Jacobi polynomials.
Our main result of this section is a new proof of
the following theorem,
 which was first proved in \cite{ITZ2}.
\begin{theorem}
\label{conj:01}
For integers $n\geq1$ and $r\geq0$, and generic parameters $a,b,q$ for which the displayed denominators are nonzero, we have
\begin{align}
&\Pf\left((q^{i-1}-q^{j-1})\frac{(aq;q)_{i+j+r-2}}{(abq^2;q)_{i+j+r-2}}\right)_{1\leq i<j\leq 2n}
\nonumber\\&=
a^{n(n-1)}q^{n(n-1)(4n+1)/3+n(n-1)r}
\prod_{k=1}^n\frac{(aq;q)_{2k+r-1}(bq;q)_{2(k-1)}(q;q)_{2k-1}}
{(abq^2;q)_{2(k+n)+r-3}}.
\label{eq:pf-special}
\end{align}
\end{theorem}
Let us use the notation (see \cite[p.1476]{Hab1})
\begin{align}
&
\Delta_{k}^{0}(\x)
=\prod_{i<j} \left( x_i/x_j;q\right)_{k} \left( q x_j/x_i;q\right)_{k},
\label{eq:Delta0}
\\&
\Delta_{k}(\x)=\frac1{n!}\sum_{\sigma\in\Sym_{n}}
\Delta_{k}^{0}(\sigma\x).
\label{Delta}
\end{align}
The $q$-gamma function is defined on $\C\setminus \Z_{\leq0}$ by
\[
\Gamma_q(x)=(1-q)^{1-x}\frac{(q;q)_\infty}{(q^{x};q)_\infty},
\qquad 0<q<1.
\]
The following lemma appeared in \cite[Lemma~3.2]{Wa} as a constant term identity.
\begin{lemma}%[\cite{Hab1}]
We have 
\begin{equation}\label{eq:Delta1-Delta2}
\int_{[0,a]^{n}}
\Delta_{k}^{1}(\x)
\prod_{i=1}^{n}x_{i}^{r+1}
\cdot
d_q\omega(\x)
=
\frac{n!}{\Gamma_{q^k}(n+1)}
\int_{[0,a]^{n}}
\Delta_{k}^{2}(\boldsymbol{x})
\prod_{i=1}^{n}x_{i}^{r+1}
\cdot
d_q\omega(\x),
\end{equation}
where $\Delta_{k}^{1}(\x)$ is as defined in \eqref{eq:Delta1} and $\Delta_{k}^{2}(\boldsymbol{x})$ is defined by
\begin{equation}
\Delta_{k}^{2}(\boldsymbol{x})
=\prod_{1\leq i<j\leq n}
x_i^{2k} \left( q^{1-k}x_j/x_i;q\right)_{2k}
=\prod_{1\leq i<j\leq n}\prod_{\nu=-k+1}^{k}
(x_{i}-q^{\nu}x_{j})
\label{eq:Delta2}
\end{equation}
(see \cite[(3.30)]{BF}).
\end{lemma}
\begin{demo}{Proof}
If we use 
\begin{equation}
\Delta_{k}(\x)
=\frac{\Gamma_{q^k}(n+1)}{n!}\prod_{i\neq j}
\left(x_i/x_j;q\right)_{k}
\label{eq:Hab}
\end{equation}
which is proved in \cite[(2.8)]{Hab1},
then a direct computation shows
\begin{equation}
\Delta_{k}(\x)
=(-1)^{k\binom{n}2}q^{\binom{k}2\binom{n}2}
\frac{\Gamma_{q^k}(n+1)}{n!}
\prod_{i=1}^{n}x_{i}^{-k(n-1)}
\cdot
\Delta_{k}^{1}(\x)
\end{equation}
(see \cite[(3.1)]{Hab1}).
This implies that
\begin{equation}
\int_{[0,a]^{n}}
\prod_{i=1}^{n}x_{i}^{r+1}
\cdot
\Delta_{k}^{1}(\x)
\,d_q\omega(\x)
=
\frac{(-1)^{k\binom{n}2}q^{-\binom{k}2\binom{n}2}n!}{\Gamma_{q^k}(n+1)}
\int_{[0,a]^{n}}
\prod_{i=1}^{n}x_{i}^{r+1+k(n-1)}
\cdot
\Delta_{k}(\x)
\,d_q\omega(\x).
\label{eq:Delta_k-Delta_k^1}
\end{equation}
Since
$\int_{[0,a]^n}f(\sigma \x)d_{q}\x
=\int_{[0,a]^n}f(\x)d_{q}\x$ for any $\sigma \in \Sym_{n}$, 
we obtain
\begin{equation}
\int_{[0,a]^{n}}
\prod_{i=1}^{n}x_{i}^{r+1}
\cdot
\Delta_{k}^{1}(\x)
\,d_q\omega(\x)
=
\frac{(-1)^{k\binom{n}2}q^{-\binom{k}2\binom{n}2}n!}{\Gamma_{q^k}(n+1)}
\int_{[0,a]^{n}}
\prod_{i=1}^{n}x_{i}^{r+1+k(n-1)}
\cdot
\Delta_{k}^{0}(\x)
\,d_q\omega(\x).
\label{eq:Delta1-Delta0}
\end{equation}
A direct computation shows
\[
x_i^{2k} \left( q^{1-k}x_j/x_i;q\right)_{2k}
=(-1)^k q^{-{k\choose 2}} (x_ix_j)^k  \left( x_i/x_j;q\right)_{k} \left( q x_j/x_i;q\right)_{k},
\]
which implies
\begin{equation}
\Delta_{k}^{0}(\x)
=(-1)^{k\binom{n}2}q^{\binom{k}2\binom{n}2}
\prod_{i=1}^{n}x_{i}^{-k(n-1)}
\prod_{1\leq i<j\leq n}x_i^{2k} \left( q^{1-k}x_j/x_i;q\right)_{2k}.
\label{eq:prod0}
\end{equation}
Hence \eqref{eq:Delta1-Delta0} and \eqref{eq:prod0} immediately
imply \eqref{eq:Delta1-Delta2}.
\end{demo}
From \eqref{eq:r-mom-jacobi} and \eqref{eq:Delta1-Delta2},
we derive the key identity 
\begin{equation}
\Pf\Bigl(
(q^{i-1}-q^{j-1})\mu_{i+j+r-2}
\Bigr)_{1\leq i<j\leq2n}
=\frac{q^{n(n-1)}(1-q)^{n}}{\Gamma_{q^2}(n+1)}
\int_{[0,a]^{n}}
\Delta_{2}^{2}(\boldsymbol{x})
\prod_{i=1}^{n}x_{i}^{r+1}
\cdot
d_q\omega(\x)
\label{eq:Pf-Delta2}
\end{equation}
to compute our $q$-Pfaffians.
\par\bigskip
Next we begin our proof by recalling the notation of the little $q$-Jacobi polynomials.
Let 
\begin{align*}
{}_{2}\phi_{1}\left[\,
{{a,b}\atop{c}};q,z
\,\right]
=\sum_{n=0}^{\infty}\frac{(a,b;q)_{n}}{(q,c;q)_{n}}z^{n}
\end{align*}
denote the \defterm{basic hypergeometric series}.
The \defterm{little $q$-Jacobi polynomials} \cite{GR,KLS} are defined by
\begin{equation}
p_{n}(x;a,b;q)=\frac{(aq;q)_{n}}{(abq^{n+1};q)_{n}}(-1)^{n}q^{\binom{n}2}
{}_{2}\phi_{1}\left[
{{q^{-n}, abq^{n+1}}
\atop{aq}}
\,;\, q, xq
\right],
\label{eq:little-q-Jacobi}
\end{equation}
which are orthogonal with respect to the inner product 
defined as
\begin{align}
&\int_{0}^{1}f(x)g(x)\,d_{q}\omega(x)
=\frac{(aq;q)_{\infty}}{(abq^2;q)_{\infty}}
\sum_{k=0}^{\infty}
\frac{(bq;q)_{k}}{(q;q)_{k}}(aq)^{k}f\left(q^{k}\right)g\left(q^{k}\right)
\nonumber\\&\qquad\qquad
=\frac{(aq,bq;q)_{\infty}}{(q,abq^2;q)_{\infty}}
\sum_{k=0}^{\infty}
\frac{(q^{k+1};q)_{\infty}}{(bq^{k+1};q)_{\infty}}
(aq)^kf\left(q^{k}\right)g\left(q^{k}\right).
\label{eq:inner-product}
\end{align}
Hence the measure is given by the weight function
\begin{equation}
w(x)
=\frac1{1-q}\cdot
\frac{(aq,bq;q)_{\infty}}{(q,abq^2;q)_{\infty}}\cdot
\frac{(q x;q)_{\infty}}{(q^{\beta+1} x;q)_{\infty}}
\, x^{\alpha},
\label{eq:weight-little-Jacobi}
\end{equation}
where $a=q^{\alpha}$ and $b=q^{\beta}$.
Theorem~\ref{conj:01} is equivalent to a polynomial identity in $a$ and $b$
if we multiply both sides by the denominator 
of the right-side of \eqref{eq:pf-special}. 
So it is sufficient to prove it
 for $a=q^\alpha$ and $b=q^\beta$.
By the $q$-binomial formula, 
the $n$th moment of the little $q$-Jacobi polynomials is
\begin{align}\label{eq:moment}
\mu_{n}=\int_{0}^{1} x^{n}w(x)\,d_{q}x
=\frac{(aq;q)_{n}}{(abq^2;q)_{n}}\quad(n=0,1,2,\dots).
\end{align}
Let 
\begin{equation}
A_n(x,y,k;q)=\prod_{j=1}^n
\frac{\Gamma_q(x+(j-1)k)\Gamma_q(y+(j-1)k) \Gamma_q(jk+1)}
{\Gamma_{q}(x+y+(n+j-2)k)\Gamma_q(k+1)}.
\label{eq:A_n}
\end{equation}
Askey proposed several $q$-analogues of the Selberg integral in \cite{Ask}.
The following is the well-known Askey-Habsieger-Kadell integral, 
conjectured by Askey~\cite[Conjecture 1]{Ask},
 and proved independently by Habsieger~\cite[(3.2), (3.3)]{Hab1} and Kadell \cite[Theorem 2;
$l=m=0$]{Kad},
\begin{align}\label{eq:ask-conj1}
\int_{[0,1]^n} \prod_{i<j} t_i^{2k}& 
\left( q^{1-k}t_j/t_i;q\right)_{2k} \prod_{i=1}^n t_i^{x-1} \frac{(t_iq;q)_\infty}{(t_iq^y;q)_\infty} d_q\boldsymbol{t}
=q^{kx{n\choose 2}+2k^2{n\choose 3}}A_n(x,y,k;q).
\end{align}
Recently Kim and Stanton \cite{KS} gave a combinatorial interpretation of a 
$q$-Selberg integral, which is equivalent to \eqref{eq:ask-conj1}.
\begin{lemma}
We have 
\begin{equation}\label{eq:int-Delta2}
\int_{[0,1]^n} 
\Delta_{k}^{2}(\boldsymbol{t})
\prod_{i=1}^{n}t_{i}^{r+1}w(t_{i})
\cdot
d_q\boldsymbol{t}
=
\left\{\frac{(aq,bq;q)_{\infty}}{(q,abq^2;q)_{\infty}}\right\}^{n}
\frac{
q^{k(\alpha+r+2){n\choose 2}+2k^2{n\choose 3}}
A_n(\alpha+r+2,\beta+1,k;q)
}
{(1-q)^n},
\end{equation}
where $w$ is the weight function \eqref{eq:weight-little-Jacobi}.
\end{lemma}
\begin{demo}{Proof}
From \eqref{eq:weight-little-Jacobi} and \eqref{eq:Delta2},
we obtain
\begin{align*}
I&=
\int_{[0,1]^n} 
\Delta_{k}^{2}(\boldsymbol{t})
\prod_{i=1}^{n}t_{i}^{r+1}w(t_{i})
\cdot
d_q\boldsymbol{t}
\\&=
\frac{1}{(1-q)^n}
\left\{\frac{(aq,bq;q)_{\infty}}{(q,abq^2;q)_{\infty}}\right\}^{n}\int_{[0,1]^n} 
\prod_{i<j}
t_i^{2k} \left( q^{1-k}t_j/t_i;q\right)_{2k}
\cdot
\prod_{i=1}^{n}t_{i}^{\alpha+r+1}
\frac{(t_{i} q;q)_{\infty}}{(t_{i} q^{\beta+1} ;q)_{\infty}}
d_q\boldsymbol{t}.
\end{align*}
It is easy to see that this equals the desired identity \eqref{eq:int-Delta2}
from \eqref{eq:ask-conj1}.
\end{demo}
\begin{demo}{Proof of Theorem~\ref{conj:01}} 
If we put $k=2$,
we obtain the following identity from \eqref{eq:Pf-Delta2} and \eqref{eq:int-Delta2} :
\begin{align}
&\Pf\Bigl(
(q^{i-1}-q^{j-1})\mu_{i+j+r-2}
\Bigr)_{1\leq i<j\leq2n}
\nonumber\\&
=
\left\{\frac{(aq,bq;q)_{\infty}}{(q,abq^2;q)_{\infty}}\right\}^{n}
\frac{
q^{n(n-1)+2(\alpha+r+2){n\choose 2}+8{n\choose 3}}
A_n(\alpha+r+2,\beta+1,2;q)
}
{\Gamma_{q^2}(n+1)}.
\label{eq:Pf-Askey}
\end{align}
From \eqref{eq:A_n} we have
\[
A_n(\alpha+r+2,\beta+1,2;q)
=\prod_{j=1}^{n}
\frac{\Gamma_q(\alpha+r+2+2(j-1))\Gamma_q(\beta+1+2(j-1)) \Gamma_q(2j+1)}
{\Gamma_{q}(\alpha+\beta+r+3+2(n+j-2))\Gamma_q(3)}.
\]
Using $\Gamma_{q}(x+k)=\Gamma_{q}(x)\frac{(q^x:q)_{k}}{(1-q)^{k}}$
when $k$ is nonnegative integer,
we see that this equals
\begin{align*}
&\left\{\frac{\Gamma_q(\alpha+1)\Gamma_q(\beta+1)}{\Gamma_{q}(\alpha+\beta+2)}\right\}^{n}
\prod_{j=1}^{n}\frac{(q^{\alpha+1};q)_{r+1+2(j-1)}(q^{\beta+1};q)_{2(j-1)}(q;q)_{2j}}
{(q^{\alpha+\beta+2}:q)_{r+1+2(n+j-2)}(q;q)_{2}}
\\&
=\left\{\frac{(abq^{2},q;q)_{\infty}}{(1-q^2)(aq,bq;q)_{\infty}}\right\}^{n}
\prod_{j=1}^{n}\frac{(aq;q)_{r+1+2(j-1)}(bq;q)_{2(j-1)}(q;q)_{2j}}
{(abq^{2}:q)_{r+1+2(n+j-2)}}.
\end{align*}
If we substitute $A_n(\alpha+r+2,\beta+1,2;q)$ into \eqref{eq:Pf-Askey} and use 
$(1-q^2)^{n}\Gamma_{q^2}(n+1)=\prod_{j=1}^{n}(1-q^{2j})$,
then 
we immediately obtain the desired identity \eqref{eq:pf-special}
by a straightforward computation.
\end{demo}

\section{Hankel Pfaffians for Rogers-Szeg\"o polynomials}
\label{sec:Al-Salam-Carlitz}
The aim of this section is to prove Theorem~\ref{conj:03}
in which the first two identities are conjectured in  \cite[Conjecture~6.1]{ITZ2}.
In \cite{ITZ2} we treated the Al-Salam-Carlitz I polynomials only,
and only
the product formulas are conjectured.
But here we extend the conjectured identities to the Al-Salam-Carlitz I and II polynomials,
and two more additive identities,
i.e.,
\eqref{eq:U:r=1} and \eqref{eq:V:r=1},
are added.
The key tool for the proof is \eqref{eq:Pf-Delta2} as well,
and we apply the results in \cite{BF}.
\par
This section is composed as follows.
First we recall the Al-Salam-Carlitz I, II polynomials
and their weight functions.
We state our main result in Theorem~\ref{conj:03},
then we define the multivariate Al-Salam-Carlitz polynomials
$\{U^{(a)}_{\lambda}(\x;q,t)\}$ and $\{V^{(a)}_{\lambda}(\x;q,t)\}$,
which are both basis of the vector space of the symmetric functions 
in the variables $\x=(x_{1},\dots,x_{n})$
over $\C(q,t)$.
The orthogonality \cite[(4.27), (4.29), (4.34), (4.35)]{BF}
 of $\{U^{(a)}_{\lambda}(\x;q,t)\}$ and $\{V^{(a)}_{\lambda}(\x;q,t)\}$
plays an important role to prove the theorem.
We use a result \cite[(A4), (A5), (A6)]{BF} which tells us how to express the elementary symmetric function
$e_{r}(\x)$ by means of $\{U^{(a)}_{\lambda}(\x;q,t)\}$ or $\{V^{(a)}_{\lambda}(\x;q,t)\}$
and apply the above orthogonality.
\par\smallskip
In this section we use the notation:
\begin{align*}
&
e_{q}(x)=\sum_{n=0}^{\infty}\frac{x^n}{(q;q)_{n}}=\frac1{(x;q)_{\infty}},
\qquad
E_{q}(x)=\sum_{n=0}^{\infty}\frac{q^{\frac{n(n-1)}2}x^n}{(q;q)_{n}}=(-x;q)_{\infty}.
\end{align*}
These functions are the classical $q$-analogues of the exponential function.
Al-Salam and Carlitz have introduced two families $\{U^{(a)}_{n}(x;q)\}_{n=0}^{\infty}$ 
and $\{V^{(a)}_{n} (x; q)\}_{n=0}^{\infty}$ of polynomials by means of
the generating functions
\begin{align}
&
\rho_{a}(x;q)e_{q}(xy)=\sum_{n=0}^{\infty}U^{(a)}_{n}(y;q)\frac{x^n}{(q;q)_{n}},
\\&
\frac1{\rho_a(x;q)}E_{q}(-xy) 
=\sum_{n=0}^{\infty}
V^{(a)}_{n}(y;q)
\frac{(-1)^{n}q^{\frac{n(n-1)}2}x^n}{(q;q)_{n}},
\end{align}
where $\rho_{a}(x;q)=(x;q)_{\infty}(ax;q)_{\infty}=E_{q}(-x)E_{q}(-ax)$.
The families
$\{U^{(a)}_{n}(x;q)\}_{n=0}^{\infty}$ 
and
$\{V^{(a)}_{n} (x; q)\}_{n=0}^{\infty}$
are called the
\defterm{Al-Salam-Carlitz I polynomials}
and
\defterm{Al-Salam-Carlitz II polynomials},
respectively.
Al-Salam and Carlitz obtain the explicit orthogonal relations (see \cite{AC,Ch})
\begin{align}
&
\int_{a}^{1}U^{(a)}_{m}(x;q)U^{(a)}_{n}(x;q)w_{U}^{(a)}(x;q)\,d_qx
=(1-q)(-a)^{n}q^{\frac{n(n-1)}2}(q;q)_{n}\delta_{m,n},
\\&
\int_{1}^{\infty}V^{(a)}_{m}(x;q)V^{(a)}_{n}(x;q)w_{V}^{(a)}(x;q)\,d_qx
=(1-q)a^{n}q^{-n^2}(q;q)_{n}\delta_{m,n},
\end{align}
where
\begin{align*}
&
w_{U}^{(a)}(x;q)=\frac{\left(qx;q\right)_{\infty}\left(\frac{qx}{a};q\right)_{\infty}}
{\left(q;q\right)_{\infty}\left(aq;q\right)_{\infty}\left(\frac{q}{a};q\right)_{\infty}},
\qquad
w_{V}^{(a)}(x;q)=\frac
{\left(q;q\right)_{\infty}\left(aq;q\right)_{\infty}\left(\frac{q}{a};q\right)_{\infty}}
{\left(x;q\right)_{\infty}'\left(\frac{x}{a};q\right)_{\infty}}.
\end{align*}
Here $\left(x;q\right)_{\infty}'$ stands for the product except zeros,
and the above Jackson integrals are defined by
\begin{align*}
&
\int_{a}^{1}f(x)d_qx
=(1-q)\left\{\sum_{n=0}^{\infty}f(q^n)q^n-a\sum_{n=0}^{\infty}f(aq^n)q^n\right\},
\\&
\int_{1}^{\infty}f(x)d_qx
=(1-q)\sum_{n=0}^{\infty}f(q^{-n})q^{-n}.
\end{align*}
Recall that $\{U^{(a)}_{n}(x;q)\}_{n=0}^{\infty}$ is positive-definite for $a<0$ and $0<q<1$,
and
$\{V^{(a)}_{n} (x; q)\}_{n=0}^{\infty}$ is positive-definite for $a>0$ and $0<q<1$.
Then the moments $F_{n}(a;q)$ and $G_{n}(a;q)$ are defined by
\begin{align}
&
\int_{a}^{1}x^nw_{U}^{(a)}(x;q)d_qx=(1-q)\,F_{n}(a;q),
\quad
\int_{1}^{\infty}x^nw_{V}^{(a)}(x;q)d_qx=(1-q)\,G_{n}(a;q)
\label{eq:BF-moments}
\end{align}
and have the expressions \cite{GR,KLS}
\begin{align}
&
F_{n}(a;q)=\sum_{k=0}^{n}\left[{{n}\atop{k}}\right]_{q}a^{k},
\qquad\qquad
G_{n}(a;q)=\sum_{k=0}^{n}\left[{{n}\atop{k}}\right]_{q}a^{k}q^{k(k-n)},
\end{align}
where $\left[{{n}\atop{k}}\right]_{q}=\frac{(q;q)_{n}}{(q;q)_{k}(q;q)_{n-k}}$.
We call $F_{n}(a;q)$ the Rogers-Szeg\"o polynomials 
\cite[Ex.~7.40]{GR}.
The main result of this section is the following:
\begin{theorem}
\label{conj:03}
Let $F_{n}(a;q)$ and $G_{n}(a;q)$ be as above.
Then the following identities hold:
\begin{align}
&\Pf\biggl((q^{i-1}-q^{j-1})F_{i+j-3}(a;q)\biggr)_{1\leq i<j\leq 2n}
=
a^{n(n-1)}q^{\frac16n(n-1)(4n-5)}\prod_{k=1}^{n}(q;q)_{2k-1},
\label{eq:U:r=-1}
\\
&\Pf\biggl((q^{i-1}-q^{j-1})F_{i+j-2}(a;q)\biggr)_{1\leq i<j\leq 2n}
=a^{n(n-1)}q^{\frac16n(n-1)(4n+1)}
\prod_{k=1}^{n}(q;q)_{2k-1}
\cdot G_{n}(a;q^2),
\label{eq:U:r=0}
\\
&\Pf\biggl((q^{i-1}-q^{j-1})F_{i+j-1}(a;q)\biggr)_{1\leq i<j\leq 2n}
=
(-a)^{n(n-1)}
q^{n(n-1)(4n+7)/6}
\prod_{i=1}^{n}(q;q)_{2i-1}
\nonumber\\&
\qquad\qquad\qquad\qquad\qquad\qquad
\times
\sum_{i=0}^{n}
(-a)^{i}q^{i(i+1)-2ni}
(q;q^{2})_{i}\left[{{n}\atop{i}}\right]_{q^2}
G_{n-i}(a;q^2)^2,
\label{eq:U:r=1}
\\
&\Pf\biggl((q^{i-1}-q^{j-1})G_{i+j-3}(a;q)\biggr)_{1\leq i<j\leq 2n}
=
a^{n(n-1)}q^{-n(n-1)(4n-5)/3}\prod_{k=1}^{n}(q;q)_{2k-1},
\label{eq:V:r=-1}
\\&
\Pf\biggl((q^{i-1}-q^{j-1})G_{i+j-2}(a;q)\biggr)_{1\leq i<j\leq 2n}
=a^{n(n-1)}q^{-\frac23n(n-1)(2n-1)}\prod_{k=1}^{n}(q;q)_{2k-1}
\cdot F_{n}(a;q^2),
\label{eq:V:r=0}
\\&
\Pf\biggl((q^{i-1}-q^{j-1})G_{i+j-1}(a;q)\biggr)_{1\leq i<j\leq 2n}
=
a^{n(n-1)}q^{-n(n-1)(4n+1)/3}
\prod_{i=1}^{n}(q;q)_{2i-1}
\nonumber\\&
\qquad\qquad\qquad\qquad\qquad\qquad\qquad\qquad
\times
\sum_{i=0}^{n}
a^{i}q^{-i}
(q;q^{2})_{i}\left[{{n}\atop{i}}\right]_{q^2}
F_{n-i}(a;q^{2})^2.
\label{eq:V:r=1}
\end{align}
\end{theorem}
\begin{remark}
The first (resp. the second) identity in \cite[Conjecture~6.1]{ITZ2}
obviously corresponds to \eqref{eq:U:r=-1} (resp. \eqref{eq:U:r=0}). 
The powers of $q$ in conjectured identities look complicated and different.
However \eqref{eq:U:r=-1} is equivalent to the conjectured identity and
\eqref{eq:U:r=0} corrects the mistaken conjectured formula.
\end{remark}
\par\bigskip
Let $P_{\lambda}(\x;q,t)$ denote Macdonald's $P$-function with respect to a partition $\lambda$
in $n$-tuple $\x=(x_{1},\dots,x_{n})$ of variables
(see \cite[Chap. VI, (4.7)]{Mac}).
In what follows, we use the notation
\begin{equation}
\begin{array}{lll}
a_{\lambda}(s)=\lambda_{i}-j,&\qquad&a_{\lambda}'(s)=j-1,\\
l_{\lambda}(s)=\lambda_{j}'-i,&\qquad&l_{\lambda}'(s)=i-1
\end{array}
\end{equation}
for each cell $s=(i,j)\in\lambda$ (\cite[VI, (6.14)]{Mac}),
where $\lambda'$ denotes the conjugate of $\lambda$.
We also use the notation
\begin{equation}
h_{\lambda}(q,t)=\prod_{s\in\lambda}(1-q^{a_{\lambda}(s)}t^{l_{\lambda}(s)+1})
,\qquad
h_{\lambda}'(q,t)=\prod_{s\in\lambda}(1-q^{a_{\lambda}(s)+1}t^{l_{\lambda}(s)})
\end{equation}
(see \cite[(1.3), (1.4)]{BF}).
The hypergeometric functions ${}_{0}\mathcal{F}_{0}(\x;\boldsymbol{y};q,t)$ 
and
${}_{0}\varphi_{0}(\x;\boldsymbol{y};q,t)$ are defined by
\begin{align}
&
{}_{0}\mathcal{F}_{0}(\x;\boldsymbol{y};q,t)
=\sum_{\lambda}\frac{t^{n(\lambda)}}{h_{\lambda}'(q,t)P_{\lambda}(t^{\bar\delta};q,t)}
P_{\lambda}(\x;q,t)P_{\lambda}(\boldsymbol{y};q,t),
\\&
{}_{0}\varphi_{0}(\x;\boldsymbol{y};q,t)
=\sum_{\lambda}\frac{(-1)^{|\lambda|}q^{n(\lambda')}}
{h_{\lambda}'(q,t)P_{\lambda}(t^{\bar\delta};q,t)}
P_{\lambda}(\x;q,t)P_{\lambda}(\boldsymbol{y};q,t),
\end{align}
where
$n(\lambda)=\sum_{i\geq1}(i-1)\lambda_{i}$
(\cite[(1.3), (1.4)]{BF}, \cite[I, (1.5)]{Mac})
and $t^{\bar\delta}=(1,t,\dots,t^{n-1})$.
In \cite[VI, (6.11')]{Mac} (also see \cite[(1.4)]{BF}) 
an explicit description of $P_{\lambda}(t^{\bar\delta};q,t)$ is given that
\begin{equation}
P_{\lambda}(t^{\bar\delta};q,t)
=t^{n(\lambda)}\prod_{s\in\lambda}\frac{1-q^{a_{\lambda}'(s)}t^{n-l_{\lambda}'(s)}}
{1-q^{a_{\lambda}(s)}t^{l_{\lambda}(s)+1}}
=t^{n(\lambda)}\frac{\prod_{(i,j)\in\lambda}(1-q^{j-1}t^{n-i+1})}{h_{\lambda}(q,t)}.
\end{equation}
Two families $\{U^{(a)}_{\lambda}(\x;q,t)\}$ and $\{V^{(a)}_{\lambda}(\x;q,t)\}$ 
of multivariate polynomials are introduced in \cite[(2.32), (2.33)]{BF}
using the following generating functions:
\begin{align}
\prod_{i=1}^{n}\rho_{a}(y_{i};q)\cdot
{}_{0}\mathcal{F}_{0}(\x;\boldsymbol{y};q,t)
&=\sum_{\lambda}\frac{t^{n(\lambda)}P_{\lambda}(\boldsymbol{y};q,t)}{h_{\lambda}'(q,t)P_{\lambda}(t^{\bar\delta};q,t)}
\,U^{(a)}_{\lambda}(\x;q,t),
\label{eq:gen-U}
\\
\frac1{\prod_{i=1}^{n}\rho_{a}(t^{-(n-1)}y_{i};q)}\cdot
{}_{0}\varphi_{0}(\x;\boldsymbol{y};q,t)
&=\sum_{\lambda}\frac{(-1)^{|\lambda|}q^{n(\lambda')}P_{\lambda}(\boldsymbol{y};q,t)}
{h_{\lambda}'(q,t)P_{\lambda}(t^{\bar\delta};q,t)}
\,
V^{(a)}_{\lambda}(\x;q,t).
\label{eq:gen-V}
\end{align}
It is shown in \cite[(2.34)]{BF} that the following relation holds:
\begin{equation}
V^{(a)}_{\lambda}(\x;q,t)=U^{(a)}_{\lambda}(\x;q^{-1},t^{-1}).
\label{eq:V-U}
\end{equation}
Hereafter we write 
$d_{q}\mu^{(U)}(\x)=\prod_{i=1}^nw^{(a)}_{U}(x_i;q)\,d_q\x$
and
$d_{q}\mu^{(V)}(\x)=\prod_{i=1}^nw^{(a)}_{V}(x_i;q)\,d_q\x$.
Baker and Forrester proved the following identities:
\begin{align}
\Norm_{\emptyset}^{(U)}(a;q,t)
&=
\int_{[a,1]^n}\Delta_{k}^{2}(\x)
d_{q}\mu^{(U)}(\x)
=
(1-q)^{n}(-a)^{kn(n-1)/2}t^{k\binom{n}3-\frac{k-1}2\binom{n}2}
\prod_{i=1}^{n}\frac{(q;q)_{ki}}{(q;q)_{k}},
\label{eq:BF-U}
\\
\Norm_{\emptyset}^{(V)}(a;q,t)
&=
\int_{[1,\infty]^n}\Delta_{k}^{2}(\x)
d_{q}\mu^{(V)}(\x)
=
(1-q)^na^{kn(n-1)/2}t^{-2k\binom{n}3-k\binom{n}2}
\prod_{i=1}^{n}\frac{(q;q)_{ki}}{(q;q)_{k}}
\label{eq:BF-V}
\end{align}
for $t=q^{k}$ with any nonnegative integer $k$,
where $\Delta_{k}^{2}(\x)$ is as in \eqref{eq:Delta2}
(see \cite[(4.27), (4.29)]{BF}).
They also prove the orthogonality
\begin{align}
&\int_{[a,1]^n}
U^{(a)}_{\lambda}(\x;q,t)
U^{(a)}_{\mu}(\x;q,t)
\Delta_{k}^{2}(\x)
d_{q}\mu^{(U)}(\x)
=\Norm_{\lambda}^{(U)}(a;q,t)\,\delta_{\lambda \mu},
\label{eq:orth-U}
\\&
\int_{[1,\infty]^n}
V^{(a)}_{\lambda}(\x;q,t)
V^{(a)}_{\mu}(\x;q,t)
\Delta_{k}^{2}(\x)
d_{q}\mu^{(V)}(\x)
=\Norm_{\lambda}^{(V)}(a;q,t)\,\delta_{\lambda \mu}
\label{eq:orth-V}
\end{align}
of $\{U^{(a)}_{\lambda}(\x;q,t)\}$ and $\{V^{(a)}_{\lambda}(\x;q,t)\}$ 
(\cite[(3.29), (3.37)]{BF})
and obtain the normalization integrals 
\begin{align}
\Norm_{\lambda}^{(U)}(a;q,t)
&=
(-at^{n-1})^{|\lambda|}q^{n(\lambda')}t^{-2n(\lambda)}
h_{\lambda}'(q,t) P_{\lambda}(t^{\bar\delta};q,t)
\Norm_{\emptyset}^{(U)}(a;q,t),
\label{eq:BF-U-lambda}
\\
\Norm_{\lambda}^{(V)}(a;q,t)
&=
(aq^{-1}t^{-2(n-1)})^{|\lambda|}q^{-2n(\lambda')}t^{n(\lambda)}
h_{\lambda}'(q,t) P_{\lambda}(t^{\bar\delta};q,t)
\Norm_{\emptyset}^{(V)}(a;q,t)
\label{eq:BF-V-lambda}
\end{align}
for any $\lambda$
(\cite[(4.34), (4.35)]{BF})
when $t=q^{k}$, $k=0,1,2,\dots$.
\begin{demo}{Proof of Theorem~\ref{conj:03}}
In this proof we always assume $t=q^k$, where $k$ is a nonnegative integer.
However, we basically use only the case when $k=2$ to prove Theorem~\ref{conj:03}.
Substituting \eqref{eq:BF-moments}, \eqref{eq:BF-U} and \eqref{eq:BF-V}
into \eqref{eq:Pf-Delta2} with $r=-1$ and $k=2$,
we obtain
\begin{align*}
&\Pf\biggl((q^{i-1}-q^{j-1})F_{i+j-3}(a;q)\biggr)_{1\leq i<j\leq 2n}
=\frac{q^{n(n-1)}(1-q)^{n}}{\Gamma_{q^2}(n+1)}
a^{n(n-1)}q^{4\binom{n}3-\binom{n}2}
\prod_{i=1}^{n}\frac{(q;q)_{2i}}{(q;q)_{2}},
\\
&\Pf\biggl((q^{i-1}-q^{j-1})G_{i+j-3}(a;q)\biggr)_{1\leq i<j\leq 2n}
=\frac{q^{n(n-1)}(1-q)^{n}}{\Gamma_{q^2}(n+1)}
a^{n(n-1)}q^{-8\binom{n}3-4\binom{n}2}
\prod_{i=1}^{n}\frac{(q;q)_{2i}}{(q;q)_{2}},
\end{align*}
which proves \eqref{eq:U:r=-1} and \eqref{eq:V:r=-1} by direct computation.
\par\smallskip
Next 
we need to expand the $n$th elementary symmetric function
$\prod_{i=1}^{n}x_i=e_n(\x)$ in terms of 
 $U^{(a)}_{\lambda}(\x;q,t)$
(resp. $V^{(a)}_{\lambda}(\x;q,t)$) to prove \eqref{eq:U:r=0} 
(resp. \eqref{eq:V:r=0}).
Then we use the orthogonality \eqref{eq:orth-U} or \eqref{eq:orth-V}
to calculate the right-hand side of \eqref{eq:Pf-Delta2} when $r=0$.
First we evaluate the integrals
\begin{equation}
\int_{[a,1]^{n}} \prod_{i=1}^{n}x_{i}
\cdot \Delta_{k}^{2}(\x) d_{q}\mu^{(U)}(\x)
\quad\text{ and }\quad
\int_{[1,\infty]^{n}} \prod_{i=1}^{n}x_{i} 
\cdot \Delta_{k}^{2}(\x) d_{q}\mu^{(V)}(\x)
\label{eq:integrals-BF}
\end{equation}
for any nonnegative integer $k$.
Baker~and~Forrester \cite[(A.4), (A.5), (A.6)]{BF} show the $r$th elementary symmetric function
$e_r(\x)$ can be expanded in $U^{(a)}_{\lambda}(\x;q,t)$ as
\begin{equation}
e_r(\x)
=\sum_{i=0}^{r}\tilde{f}_{r-i}(a)\left[{{n-i}\atop{r-i}}\right]_{t}
U^{(a)}_{(1^i)}(\x;q,t),
\label{eq:e-U}
\end{equation}
where $\tilde{f}_{i}(a)$ is defined by the initial condition $\tilde{f}_{0}(a)=1$,
 $\tilde{f}_{1}(a)=1+a$ 
and the recurrence equation
\begin{align*}
&
\tilde{f}_{i}(a)
=(1+a)t^{i-1}\tilde{f}_{i-1}(a)
+at^{i-2}(1-t^{i-1})\tilde{f}_{i-2}(a).
\end{align*}
Using this recurrence, we can show that
\begin{equation*}
\tilde{f}_{i}(a)
=\sum_{j=0}^{i}\left[{{i}\atop{j}}\right]_{t}
t^{\frac{j(j-1)}2+\frac{(i-j)(i-j-1)}2}a^j
=t^{\binom{i}2}G_{i}(a;t)
\end{equation*}
by induction.
Since we have $U^{(a)}_{\emptyset}(\x;q,t)=1$ from \eqref{eq:gen-U},
the first integral of \eqref{eq:integrals-BF} is written as
\begin{equation}
\int_{[a,1]^n}
\prod_{i=1}^{n}x_i
\cdot
\Delta_{k}^{2}(\x)
d_{q}\mu^{(U)}(\x)
=
\int_{[a,1]^n}
e_{n}(\x)U^{(a)}_{\emptyset}(\x;q,t)
\Delta_{k}^{2}(\x)
 d_{q}\mu^{(U)}(\x)
\label{eq:rhs:F:r=0}
\end{equation}
Substituting \eqref{eq:e-U} with $r=n$ into \eqref{eq:rhs:F:r=0},
we see this becomes
\begin{equation*}
\sum_{i=0}^{n}\tilde{f}_{n-i}(a)\int_{[a,1]^n}
U^{(a)}_{(1^i)}(\x;q,t)U^{(a)}_{\emptyset}(\x;q,t)
\Delta_{k}^{2}(\x)
d_{q}\mu^{(U)}(\x)
\end{equation*}
which equals
\begin{equation*}
\tilde{f}_{n}(a)
\int_{[a,1]^n}
\Delta_{k}^{2}(\x)
d_{q}\mu^{(U)}(\x)
\end{equation*}
from the orthogonality \eqref{eq:orth-U}.
Hence, by \eqref{eq:BF-U}, we conclude that \eqref{eq:rhs:F:r=0} equals
\begin{equation*}
\tilde{f}_{n}(a)
(1-q)^n(-a)^{kn(n-1)/2}q^{k^2\binom{n}3-\frac{k(k-1)}2\binom{n}2}
\prod_{i=1}^{n}\frac{(q;q)_{ki}}{(q;q)_{k}},
\end{equation*}
Finally,
by substituting $k=2$ into this result and using
the identity \eqref{eq:Pf-Delta2} with $r=0$, we obtain
\begin{align*}
\Pf\biggl((q^{i-1}-q^{j-1})F_{i+j-2}(a;q)\biggr)_{1\leq i<j\leq 2n}
=\frac{q^{n(n-1)}(1-q)^{n}}{\Gamma_{q^2}(n+1)}
\tilde{f}_{n}(a)
a^{n(n-1)}q^{4\binom{n}3-\binom{n}2}
\prod_{i=1}^{n}\frac{(q;q)_{2i}}{(q;q)_{2}}.
\end{align*}
One can prove \eqref{eq:U:r=0} from this identity by straight computation.
\par\bigbreak
Next we use \eqref{eq:V-U} to prove \eqref{eq:V:r=0}.
The identities \eqref{eq:e-U} and \eqref{eq:V-U} imply
\begin{equation}
e_r(\x)
=\sum_{i=0}^{r}\tilde{g}_{r-i}(a)\left[{{n-i}\atop{r-i}}\right]_{t}
V^{(a)}_{(1^i)}(\x;q,t),
\label{eq:e-V}
\end{equation}
where $\tilde{g}_{i}(a)$ is given by
$
\tilde{g}_{i}(a)
=t^{-\frac{i(i-1)}2}\sum_{j=0}^{i}\left[{{i}\atop{j}}\right]_{t}
a^j
=t^{-\binom{i}2}F_{i}(a;t).
$
A similar argument as above shows that
we can rewrite the second integral of \eqref{eq:integrals-BF} as
\begin{equation*}
\int_{[1,\infty]^{n}} 
\prod_{i=1}^{n}x_{i} 
\cdot
\Delta_{k}^{2}(\x)
d_{q}\mu^{(V)}(\x)
=\tilde{g}_{n}(a)
(1-q)^na^{kn(n-1)/2}q^{-2k^2\binom{n}3-k^2\binom{n}2}
\prod_{i=1}^{n}\frac{(q;q)_{ki}}{(q;q)_{k}}
\end{equation*}
by \eqref{eq:BF-V}, \eqref{eq:orth-V} and \eqref{eq:e-V}.
Again the $k=2$ case of this identity 
and \eqref{eq:Pf-Delta2} with $r=0$ lead to
\begin{align*}
&\Pf\biggl((q^{i-1}-q^{j-1})G_{i+j-2}(a;q)\biggr)_{1\leq i<j\leq 2n}
=\frac{q^{n(n-1)}(1-q)^{n}}{\Gamma_{q^2}(n+1)}
\tilde{g}_{n}(a)\,
a^{n(n-1)}q^{-8\binom{n}3-4\binom{n}2}
\prod_{i=1}^{n}\frac{(q;q)_{2i}}{(q;q)_{2}},
\end{align*}
which proves \eqref{eq:V:r=0} immediately.
\par\smallskip
Next 
we can use \eqref{eq:e-U}, \eqref{eq:e-V} and the orthogonality \eqref{eq:orth-U}, \eqref{eq:orth-V} to obtain
\begin{align}
\int_{[a,1]^n}
e_{n}(\x)^2
\Delta_{k}^{2}(\x)
d_{q}\mu^{(U)}(\x)
&=\sum_{i=0}^{n}\tilde{f}_{n-i}(a)^2
\Norm_{(1^{i})}^{(U)}(a;q,t),
\label{eq:rhs:F:r=1}\\
\int_{[1,\infty]^n}
e_{n}(\x)^2
\Delta_{k}^{2}(\x)
d_{q}\mu^{(V)}(\x)
&=\sum_{i=0}^{n}\tilde{g}_{n-i}(a)^2
\Norm_{(1^{i})}^{(V)}(a;q,t).
\label{eq:rhs:G:r=1}
\end{align}
Substituting
$|\lambda|=i$, 
$n(\lambda)=\binom{i}2$,
$n(\lambda')=0$,
$h_{\lambda}'(q,t)=(q;t)_{i}$
and $P_{\lambda}(t^{\bar\delta};q,t)=t^{\binom{i}2}\frac{(t^{n-i+1};t)_{i}}{(t;t)_{i}}$ 
for $\lambda=(1^{i})$
into \eqref{eq:BF-U-lambda} and \eqref{eq:BF-V-lambda},
we obtain
\begin{align}
\Norm_{\lambda}^{(U)}(a;q,t)
&=
(-a)^{i}t^{(n-1)i-\binom{i}2}
\frac{(q;t)_{i}(t^{n-i+1};t)_{i}}{(t;t)_{i}}
\Norm_{\emptyset}^{(U)}(a;q,t),
\label{eq:BF-U-lambda-con}
\\
\Norm_{\lambda}^{(V)}(a;q,t)
&=
a^{i}q^{-i}t^{-2(n-1)i+2\binom{i}2}
\frac{(q;t)_{i}(t^{n-i+1};t)_{i}}{(t;t)_{i}}
\Norm_{\emptyset}^{(V)}(a;q,t).
\label{eq:BF-V-lambda-con}
\end{align}
Substituting \eqref{eq:BF-U-lambda-con} (resp. \eqref{eq:BF-V-lambda-con})
into \eqref{eq:rhs:F:r=1} (resp. \eqref{eq:rhs:G:r=1}),
we obtain
\begin{align}
\int_{[a,1]^n}
e_{n}(\x)^2
\Delta_{k}^{2}(\x)
d_{q}\mu^{(U)}(\x)
&=(1-q)^{n}
(-a)^{kn(n-1)/2}
t^{k\binom{n}3-\frac{k-1}2\binom{n}2+n(n-1)}
\prod_{i=1}^{n}\frac{(q;q)_{ki}}{(q;q)_{k}}
\nonumber\\&\times
\sum_{i=0}^{n}
(-a)^{i}t^{-ni+i(i+1)/2}
\frac{(q;t)_{i}(t^{n-i+1};t)_{i}}{(t;t)_{i}}
G_{n-i}(a;q^2)^2,
\label{eq:rhs:F:r=1b}\\
\int_{[1,\infty]^n}
e_{n}(\x)^2
\Delta_{k}^{2}(\x)
d_{q}\mu^{(V)}(\x)
&=(1-q)^{n}
a^{kn(n-1)/2}t^{-2k\binom{n}3-k\binom{n}2-n(n-1)}
\prod_{i=1}^{n}\frac{(q;q)_{ki}}{(q;q)_{k}}
\nonumber\\&\times
\sum_{i=0}^{n}
a^{i}q^{-i}
\frac{(q;t)_{i}(t^{n-i+1};t)_{i}}{(t;t)_{i}}
F_{n-i}(a;t)^2.
\label{eq:rhs:G:r=1b}
\end{align}
By replacing $t$ by $q^{k}$, substituting $k=2$ into \eqref{eq:rhs:F:r=1b} and \eqref{eq:rhs:G:r=1b},
and using the identity \eqref{eq:Pf-Delta2} with $r=1$,
we obtain \eqref{eq:U:r=1} and \eqref{eq:V:r=1}.
This completes the proof of our theorem.
\end{demo}
The reader may notice that \eqref{eq:V:r=-1} (resp. \eqref{eq:V:r=0}, \eqref{eq:V:r=1})
can be directly derived from \eqref{eq:U:r=-1} (resp. \eqref{eq:U:r=0}, \eqref{eq:U:r=1}) 
by replacing $q$ by $q^{-1}$.
Here we give an explicit form \eqref{eq:e-V},
which tells us how to express 
$e_{r}(\x)$ by means of $\{V^{(a)}_{\lambda}(\x;q,t)\}$.
We would like to add the following remark at the end of this section.
\begin{remark}
\rm
We can evaluate the Pfaffian
$\Pf\biggl((q^{i-1}-q^{j-1})F_{i+j+r-2}(a;q)\biggr)_{1\leq i<j\leq 2n}$
for a nonnegative integer $r$ if we can find a formula which expands
the elementary symmetric function $e_{n^{r+1}}(\x)=e_{n}(x_1,\dots,x_n)^{r+1}$ by means of
$U^{(a)}_{\lambda}(\x;q,t)$.
\end{remark}

\section{Polynomial insertions and applications to beta ensembles}
\label{sec:applications}
\subsection{Moment hyperpfaffians and polynomial modifications}
Fix an even positive integer $\ell$, and put
\[
 \beta=\ell^2,\qquad \gamma=\ell^2/2,\qquad
 C=\binom\ell2,\qquad c_\ell=\prod_{r=0}^{\ell-1}r!.
\]
For an increasing $\ell$-tuple $I=(i_1,\ldots,i_\ell)$ of nonnegative integers, write
\[
 d(I)=\sum_{p=1}^{\ell}i_p-C,\qquad
 \Delta(I)=\prod_{p<q}(i_q-i_p).
\]
In particular $d(I)\ge0$. For a moment sequence $\nu$, define
\begin{equation}\label{app-eq:H}
 \HH_n[\nu]=\Pf^{[\ell,1]}
 \left(\frac{\Delta(I)}{c_\ell}\nu_{d(I)}\right)_{I\subset\{0,\ldots,\ell n-1\},\ |I|=\ell},
 \qquad \HH_0[\nu]=1.
\end{equation}
We use the hyperpfaffian convention of Section~\ref{sec:Pfaffian}: sum over ordered blocks, increasing within each block, divided by $n!$.

Let $w\ge0$ on a real interval $J$, and suppose all moments needed below are finite. Set
\begin{equation}\label{app-eq:Z}
 \mu_k=\int_Jx^kw(x)\dd x,\qquad
 Z_n[w]=\frac1{n!}\int_{J^n}\Delta(x)^{\ell^2}
                    \prod_{j=1}^nw(x_j)\dd x_j.
\end{equation}
Assume $Z_n[w]>0$. The labelled joint probability density is the integrand in \eqref{app-eq:Z} divided by $n!Z_n[w]$.

\begin{prop}\label{app-prop:insert}
We have $Z_n[w]=\HH_n[\mu]$. If $P(x)=\sum_{r=0}^Dp_rx^r$ and
$\mu_k^{[P]}=\sum_{r=0}^Dp_r\mu_{k+r}$, then
\begin{equation}\label{app-eq:insert}
 \EE\prod_{j=1}^nP(X_j)=\frac{\HH_n[\mu^{[P]}]}{\HH_n[\mu]}.
\end{equation}
The polynomial $P$ need not be nonnegative.
\end{prop}
\begin{proof}
Use normalized derivatives $D^r=(r!)^{-1}(\mathrm d/\mathrm dx)^r$. The monomial Wronskian is
\[
 \det\bigl(D^{q-1}x^{i_p}\bigr)_{p,q=1}^{\ell}
 =\frac{\Delta(I)}{c_\ell}x^{d(I)}.
\]
The full confluent Vandermonde determinant, with columns ordered first by particle and then by derivative order, is
\[
 \det\bigl(D^{r}x_j^i\bigr)_{\substack{0\le i<\ell n\\1\le j\le n,\ 0\le r<\ell}}
 =\Delta(x)^{\ell^2}.
\]
Its Laplace expansion into column blocks, followed by integration, gives $Z_n[w]=\HH_n[\mu]$. Replacing $w$ by $Pw$ proves \eqref{app-eq:insert}. Absolute convergence justifies the finite expansions and iterated integrals.
\end{proof}

\subsection{Jacobi weights and characteristic-polynomial powers}
Take $J=[0,1]$, $a,b>0$, and $w(x)=x^{a-1}(1-x)^{b-1}$. Thus
\[
 \mu_k=B(a+k,b),\qquad
 Z_n[w]=\frac{S_n(a,b,\gamma)}{n!},
\]
where $S_n$ is the Selberg integral in \eqref{eq:Selberg}.
For a partition $\lambda$, let $\lambda'$ be its conjugate and define
\begin{align*}
 \gp{u}{\lambda}&=\prod_{i\ge1}(u-(i-1)\gamma)_{\lambda_i},\\
 h_\lambda(\gamma)&=\prod_{(i,j)\in\lambda}
  \bigl(\lambda_i-j+\gamma(\lambda'_j-i+1)\bigr),\\
 h'_\lambda(\gamma)&=\prod_{(i,j)\in\lambda}
  \bigl(\lambda_i-j+1+\gamma(\lambda'_j-i)\bigr).
\end{align*}
Empty products are one. Put $A=a+(n-1)\gamma$ and $B=a+b+2(n-1)\gamma$; the italic parameter $B$ here is distinct from the beta function $B(u,v)$.

\begin{theorem}\label{app-thm:powers}
For every integer $m\ge0$, set
\[
 \mu_k^{(m,z)}=\sum_{r=0}^m(-1)^r\binom mr z^{m-r}B(a+k+r,b).
\]
Then the following polynomial identity holds for all $z\in\mathbb C$:
\begin{equation}\label{app-eq:powers}
 \frac{\HH_n[\mu^{(m,z)}]}{\HH_n[\mu]}
 =\sum_{\lambda\subseteq(m^n)}
 \frac{\gp{-m}{\lambda}\gp{n\gamma}{\lambda}\gp{A}{\lambda}}
      {h_\lambda(\gamma)h'_\lambda(\gamma)\gp{B}{\lambda}}
 z^{mn-|\lambda|}.
\end{equation}
Equivalently, the right-hand side evaluates $\EE\prod_{j=1}^n(z-X_j)^m$.
\end{theorem}
\begin{proof}
Let $P_\lambda^{(1/\gamma)}$ denote the monic Jack polynomial. In the present normalization, the Jack binomial identity reads
\begin{equation}\label{app-eq:binomial}
 \prod_{i=1}^n(1-tx_i)^m
 =\sum_{\lambda\subseteq(m^n)}
 \frac{\gp{-m}{\lambda}}{h'_\lambda(\gamma)}
 t^{|\lambda|}P_\lambda^{(1/\gamma)}(x).
\end{equation}
The evaluation at the all-ones alphabet and Kadell's integral give
\begin{align}
 P_\lambda^{(1/\gamma)}(1^n)
 &=\frac{\gp{n\gamma}{\lambda}}{h_\lambda(\gamma)},\label{app-eq:jackones}\\
 \EE P_\lambda^{(1/\gamma)}(X)
 &=P_\lambda^{(1/\gamma)}(1^n)
   \frac{\gp{A}{\lambda}}{\gp{B}{\lambda}}.\label{app-eq:kadell}
\end{align}
For \eqref{app-eq:binomial} and \eqref{app-eq:jackones}, see the standard Jack identities in \cite{Mac}; for \eqref{app-eq:kadell}, see \cite{Kad2} and the normalization in \cite[Eq.~(4.1)]{ARW}. Averaging \eqref{app-eq:binomial}, setting $t=1/z$, and multiplying by $z^{mn}$ proves the assertion for $z\ne0$. Polynomiality handles $z=0$. Proposition~\ref{app-prop:insert} identifies the average with the hyperpfaffian ratio.
\end{proof}

\begin{remark}
The Selberg-average side of \eqref{app-eq:powers} is classical Jack-hypergeometric material. The purpose here is its explicit realization as a finite evaluation of polynomially modified Hankel hyperpfaffians, extending the elementary-symmetric insertions in Section~\ref{sec:Motzkin}. The theorem is presented as a consequence, without a novelty claim for the integral evaluation.
\end{remark}

\subsubsection{An explicit quadratic insertion}
For $r,s\ge0$ with $r+s\le n$, set $\lambda(r,s)=(2^r,1^s)$ and
\[
 T_{r,s}(u)=\prod_{i=0}^{r-1}(u-i\gamma)(u+1-i\gamma)
             \prod_{i=r}^{r+s-1}(u-i\gamma).
\]
The two hook products specialize to
\begin{align*}
 h_{r,s}&=\prod_{i=1}^{r}
  \bigl(1+\gamma(r+s-i+1)\bigr)\gamma(r-i+1)
  \prod_{i=r+1}^{r+s}\gamma(r+s-i+1),\\
 h'_{r,s}&=\prod_{i=1}^{r}
  \bigl(2+\gamma(r+s-i)\bigr)\bigl(1+\gamma(r-i)\bigr)
  \prod_{i=r+1}^{r+s}\bigl(1+\gamma(r+s-i)\bigr).
\end{align*}
\begin{corollary}\label{app-cor:quadratic}
With $\widehat\mu_k(z)=z^2\mu_k-2z\mu_{k+1}+\mu_{k+2}$, one has
\begin{equation}\label{app-eq:quadratic}
 \frac{\HH_n[\widehat\mu(z)]}{\HH_n[\mu]}
 =\sum_{r=0}^n\sum_{s=0}^{n-r}
 \frac{T_{r,s}(-2)T_{r,s}(n\gamma)T_{r,s}(A)}
      {h_{r,s}h'_{r,s}T_{r,s}(B)}z^{2n-2r-s}.
\end{equation}
\end{corollary}
\begin{proof}
In Theorem~\ref{app-thm:powers}, take $m=2$. Every partition inside $(2^n)$ has the unique form $(2^r,1^s)$. Substitute the displayed products.
\end{proof}
Thus the quadratic average is evaluated by $(n+1)(n+2)/2$ explicit terms, without unevaluated Jack polynomials.
For comparison, $m=1$ gives the Aomoto formula
\[
 \EE\prod_{j=1}^n(z-X_j)
 =\sum_{k=0}^n(-1)^k\binom nk z^{n-k}
 \prod_{j=1}^k\frac{a+(n-j)\gamma}{a+b+(2n-j-1)\gamma}.
\]
At the other extreme, evaluation at $z=0$ gives the check
\[
 \EE\prod_j(-X_j)^m=(-1)^{mn}\frac{S_n(a+m,b,\gamma)}{S_n(a,b,\gamma)}.
\]

\subsubsection{A fourth-order hyperpfaffian example}
For $\ell=4$, $n=2$ and $a=b=1$, we have $\gamma=8$, $\mu_k=1/(k+1)$ and $c_4=12$. Directly,
\[
 \HH_2[\mu]=\frac12\int_0^1\!\int_0^1(x-y)^{16}\dd x\dd y=\frac1{306}.
\]
Corollary~\ref{app-cor:quadratic} yields the explicit identity
\begin{equation}\label{app-eq:example}
 \HH_2[\widehat\mu(z)]
 =\frac1{306}\left(z^4-2z^3+\frac{21}{19}z^2
                         -\frac2{19}z+\frac1{209}\right),
 \quad
 \widehat\mu_k(z)=\frac{z^2}{k+1}-\frac{2z}{k+2}+\frac1{k+3}.
\end{equation}
The entries of this fourth-order hyperpfaffian, indexed by four-element subsets of $\{0,\ldots,7\}$, are $\Delta(I)\widehat\mu_{d(I)}(z)/12$. The corresponding first-power average is $z^2-z+1/20$.

\subsection{Correlation functions and gap probabilities}
For $0\le q\le n$, let $R_q$ denote the factorial correlation density, so that $\int R_q=n!/(n-q)!$. Set
\[
 P_{\mathbf y}(x)=\prod_{r=1}^q(y_r-x)^{\ell^2},\qquad
 \nu_k(\mathbf y)=\int_Jx^kP_{\mathbf y}(x)w(x)\dd x.
\]
\begin{prop}\label{app-prop:corr}
For $q\ge1$,
\begin{equation}\label{app-eq:correlations}
 R_q(\mathbf y)=\Delta(\mathbf y)^{\ell^2}
 \prod_{r=1}^qw(y_r)
 \frac{\HH_{n-q}[\nu(\mathbf y)]}{\HH_n[\mu]}.
\end{equation}
\end{prop}
\begin{proof}
Fix $q$ variables in the labelled joint density and multiply its marginal by $n!/(n-q)!$. Factoring the Vandermonde separates the fixed-fixed, fixed-free and free-free interactions. The free integral is $(n-q)!Z_{n-q}[P_{\mathbf y}w]$. Proposition~\ref{app-prop:insert} gives \eqref{app-eq:correlations}.
\end{proof}
This is the monomial-moment specialization of the hyperpfaffian correlation mechanism studied by Sinclair and Wells \cite{SW}. Their work also treats arbitrary monic polynomial bases, circular ensembles, and the special reduction to Pfaffian kernels when $\ell=2$. We do not assert a kernel representation for $\ell>2$.

For the Jacobi weight, the one-point formula has completely explicit entries:
\begin{align*}
 R_1(y)&=y^{a-1}(1-y)^{b-1}
            \frac{\HH_{n-1}[\nu(y)]}{\HH_n[\mu]},\\
 \nu_k(y)&=\sum_{r=0}^{\ell^2}(-1)^r\binom{\ell^2}{r}
                        y^{\ell^2-r}B(a+k+r,b).
\end{align*}
Moreover, for $n\ge2$, Theorem~\ref{app-thm:powers} with particle number $n-1$ and insertion power $m=\ell^2$ evaluates this ratio as a finite partition sum, with prefactor
\[
 \frac{Z_{n-1}[w]}{Z_n[w]}=
 n\frac{S_{n-1}(a,b,\gamma)}{S_n(a,b,\gamma)}.
\]
For the preceding $\beta=16$, two-particle uniform example,
\[
 R_1(y)=18\bigl(y^{17}+(1-y)^{17}\bigr),\qquad \int_0^1R_1(y)\dd y=2.
\]

A further direct application concerns a measurable set $E\subseteq J$. Let $\mu_k^{J\setminus E}=\int_{J\setminus E}x^kw(x)\dd x$. Then
\[
 \mathbb P\{X_j\notin E\text{ for all }j\}
 =\frac{\HH_n[\mu^{J\setminus E}]}{\HH_n[\mu]}.
\]
For example, with $B_t(u,v)=\int_0^t x^{u-1}(1-x)^{v-1}\dd x$,
\[
 \mathbb P\{\max_jX_j\le t\}
 =\frac{\HH_n[(B_t(a+k,b))_{k\ge0}]}{\HH_n[(B(a+k,b))_{k\ge0}]},\quad 0\le t\le1.
\]
This formula is exact, but is not claimed to be a product evaluation or an efficient large-$n$ algorithm.

\subsection{Recent connections and scope of the extension}
Sinclair and Wells \cite{SW}, first posted in 2025 and revised in July 2026, discuss our moment hyperpfaffians in their treatment of hyperpfaffian evaluations. The application developed here emphasizes explicit Jacobi moments and the evaluated insertion formula \eqref{app-eq:powers}; partition functions and the general correlation mechanism are already established in the literature.

Introduce formal times by
\[
 \mu_k(\mathbf t)=\int_0^1x^k
       \exp\left(\sum_{r\ge1}t_rx^r\right)w(x)\dd x,
 \qquad \tau_n(\mathbf t)=\HH_n[\mu(\mathbf t)].
\]
Coefficientwise in the formal times,
\[
 \partial_{t_r}\mu_k=\mu_{k+r},\qquad
 \tau_n(\mathbf0)=\frac{S_n(a,b,\gamma)}{n!}.
\]
Thus the Selberg product supplies explicit initial values for time-deformed hyperpfaffians. Sinclair's March 2026 preprint \cite{Sinclair2026} develops a Hirota-type framework for charge-$L$ ensembles. The displayed identities give the moment dictionary relevant to comparing the two settings. No particular KP or BKP equation is asserted here without matching its conventions and hypotheses.

Wolff's constellation ensembles \cite{Wolff} provide another related viewpoint. In the proof of Corollary~\ref{cor:q-Hankel-hyperpfaff}, the Vandermonde variables occur in geometric clusters $q^{\ell-1}x_j,\ldots,qx_j,x_j$. Their coalescence explains the confluent determinant and the exponent $\ell^2$. A multiplicative-constellation interpretation is a possible further development; positivity of the interaction and the role of the Jackson measure must be treated before proposing a probabilistic model.

\section{Open Problems}
In \cite[Conjecture~6.3]{ITZ2}
we presented another conjecture.
We retain this conjecture as an open problem here,
but instead we add more conjectures of this type
in this last section.
Let us take $a_{n}=\frac1{2n+1}\binom{3n}{n}=\frac1{3n+1}\binom{3n+1}{n}$.
This number is famous because
Tamm \cite{Ta} proved that the Hankel determinant
$\det(a_{i+j-1})_{1\leq i,j\leq n}$ equals
the number of $(2n+1)\times(2n+1)$ alternating sign matrices
that are invariant under vertical reflection.
He evaluated the determinant by showing that the
generating function of the sequence $\{a_{n}\}_{n\geq0}$ 
is a ratio of hypergeometric series,
i.e.,
\[
g=\sum_{n=0}^{\infty}a_{n}x^{n}
={}_{2}F_{1}\left(\frac23,\frac43;\frac32;\frac{27}4x\right)
\Big/{}_{2}F_{1}\left(\frac23,\frac13;\frac12;\frac{27}4x\right),
\]
and by using the continued fraction expansion of $g$.
In \cite[Section~3]{GX} 
Gessel and Xin made a systematic
application of this continued fraction method to a number of similar Hankel
determinants.
They observed empirically that
there are five pairs of generating functions of the form
\[
\sum_{n=0}^{\infty}p_{n}x^{n}
={}_{2}F_{1}\left(a,b+1;c+1;\frac{27}4 x\right)
\Big/{}_{2}F_{1}\left(a,b;c;\frac{27}4 x\right),
\]
that can be expressed as polynomials in $g-1$.
Let $\{a^{(i)}_{n}\}_{n\geq0}$ ($i=1,\dots,5$) 
be  Gessel and Xin's five sequences (see below).
Gessel and Xin studied the Hankel determinants of these sequences. %using the continued fraction expansion of the generating functions.
By experiments %of ``Hankel Pfaffians'' for these sequences,
 we observe that the Hankel Pfaffian transforms of these sequences have mysteriously nice product formulas as in the following.
\begin{conjecture}
\label{conj:Gessel-Xin}
Let $n$ be a nonnegative integer.
The following identities would hold.
\par\noindent
If $a^{(1)}_{n}
=\frac1{3n+1}\binom{3n+1}{n}
=[x^n]\frac{{}_{2}F_{1}\left(\frac23,\frac43;\frac32;\frac{27}4x\right)}
{{}_{2}F_{1}\left(\frac23,\frac13;\frac12;\frac{27}4x\right)}$ (A001764),
then 
\begin{equation}
\Pf\left((j-i)a^{(1)}_{i+j-1}\right)_{1\leq i,j\leq2n}
=2^{-n}\prod_{k=0}^{n-1}\frac
{(12k+6)!(4k+1)!(3k+2)!}{(8k+2)!(8k+5)!(3k+1)!}.
\label{eq:Tamm}
\end{equation}
If $a^{(2)}_{n}
=\frac1{3n+2}\binom{3n+2}{n+1}
=[x^n]\frac{{}_{2}F_{1}\left(\frac43,\frac53;\frac52;\frac{27}4x\right)}
{{}_{2}F_{1}\left(\frac43,\frac23;\frac32;\frac{27}4x\right)}$ (A006013),
then 
\begin{equation}
\Pf\left((j-i)a^{(2)}_{i+j-2}\right)_{1\leq i,j\leq2n}
=12^{-n}\prod_{k=0}^{n-1}\frac{(12k+10)!(4k+2)!(4k+1)}{(8k+3)!(8k+7)!(3k+2)(12k+5)}.
\end{equation}
If $a^{(3)}_{n}
=\frac2{3n+1}\binom{3n+1}{n+1}
=2\,[x^n]\frac{{}_{2}F_{1}\left(\frac53,\frac73;\frac72;\frac{27}4x\right)}
{{}_{2}F_{1}\left(\frac53,\frac43;\frac52;\frac{27}4x\right)}$
(A007226),
then 
\begin{equation}
\Pf\left((j-i)a^{(3)}_{i+j-1}\right)_{1\leq i,j\leq2n}
=\left(\frac43\right)^{n}\prod_{k=0}^{n-1}\frac{(12k+15)!(4k+5)!(2k+1)}{(8k+8)!(8k+11)!(12k+13)}.
\end{equation}
If $a^{(4)}_{n}
=\frac2{(3n+1)(3n+2)}\binom{3n+2}{n+1}
=2\,[x^n]\frac{{}_{2}F_{1}\left(\frac53,\frac73;\frac52;\frac{27}4x\right)}
{{}_{2}F_{1}\left(\frac53,\frac43;\frac32;\frac{27}4x\right)}$
(A000139),
then 
\begin{equation}
\Pf\left((j-i)a^{(4)}_{i+j-1}\right)_{1\leq i,j\leq2n}
=\left(\frac23\right)^{n}(6n+1)!\prod_{k=0}^{n-1}\frac{(12k+6)!(4k+5)!(4k+3)}{(8k+5)!(8k+10)!(k+1)(3k+1)}.
\end{equation}
If $a^{(5)}_{n}
=\frac{9n+5}{(3n+1)(3n+2)}\binom{3n+2}{n+1}
=5\,[x^n]\frac{{}_{2}F_{1}\left(\frac23,\frac43;\frac52;\frac{27}4x\right)}
{{}_{2}F_{1}\left(\frac23,\frac13;\frac32;\frac{27}4x\right)}$
,
then 
\begin{equation}
\Pf\left((j-i)a^{(5)}_{i+j-2}\right)_{1\leq i,j\leq2n}
=3^{-n}\prod_{k=0}^{n-1}\frac{(6k+6)!(2k)!}{(4k+1)!(4k+4)!(3k+2)}.
\label{eq:A_UU}
\end{equation}
The first four of the five sequences  $\{a^{(i)}_{n}\}_{n\geq0}$ ($i=1,\dots,5$) are identified in OEIS except the last one.
\end{conjecture}
\begin{remark}
Further we identified the Hankel Pfaffian transform of $\{a^{(5)}_{n}\}_{n\geq0}$
 as A059489 in OEIS.
Namely we observe that
\begin{align}
\Pf\left((j-i)a^{(5)}_{i+j-2}\right)_{1\leq i,j\leq2n}
=A^{(2)}_{UU}(4n;1,1,1)
\end{align}
would hold for $n\geq0$.
Kuperberg \cite[Theorem~4]{Ku2} showed that the generating function
$A_{UU}(4n;x,y,z)$ of UUASMs factors as
\[
A_{UU}(4n;x,y,z)=A_{V}(2n+1;z)A^{(2)}_{UU}(4n;x,y,z),
\]
which characterizes $A^{(2)}_{UU}(4n;x,y,z)$.
Recall that a UUASM is an alternating sign matrix
with U-turn boundary on the right and on the top.
Here we define the $x$-weight to be $x^k$ if $k$ is the number of $-1$s,
$y$-weight to be $y^k$ if $k$ of the U-turns on the right are oriented upward,
 and $z$-weight to be $z^k$ if $k$ of the U-turns on the top are oriented
to the right.
Further
$A_{V}(2n+1;x)$ denotes the generating function of $(2n+1)\times(2n+1)$
vertically symmetric alternating sign matrices,
where we define the $x$-weight be the number of $-1$s, as well.
It is also a challenging problem to find the reason
which explains $A^{(2)}_{UU}(4n;1,1,1)$ appears as
the Hankel Pfaffian sequence of $\{a^{(5)}_{n}\}_{n\geq0}$.
\end{remark}
The conjectured identity \eqref{eq:Tamm} first appeared in \cite[Conjecture~6.3]{ITZ2},
but the others are new.
We note that the generating functions of these sequences are roots of certain cubic equations, 
and so the evaluation of the Pfaffians would be more difficult than the Pfaffians in Section~\ref{sec:Motzkin}.
We also note that we can find more Hankel determinants of this type 
in \cite[Theorem~31]{K2}.
We hope to address this topic in a forthcoming paper.

\par\bigskip
\noindent
{\Large\bf Acknowledgment}
\par\smallskip
We are grateful to an anonymous reviewer of our previous version
for valuable and helpful comments.

\end{document}